\documentclass[10pt]{article}
\usepackage[english]{babel}
\usepackage{amssymb,amsmath,amsthm,amsfonts,mathrsfs}
\usepackage[pdftex,pagebackref,colorlinks=true,linktocpage,citecolor=black,linkcolor=black]{hyperref}
\usepackage{graphicx}
\usepackage{epstopdf}
\usepackage{graphics}
\usepackage{subfigure}
\usepackage{diagbox}
\usepackage{bm}
\usepackage{multirow,multicol}
\usepackage{comment}
\usepackage{booktabs}
\usepackage{threeparttable}
\usepackage[ruled]{algorithm2e}
\graphicspath{{./pictures/}}
\usepackage{verbatim}
\usepackage{indentfirst}
\usepackage{mathtools}
\usepackage[usenames,dvipsnames,svgnames,table]{xcolor}
\usepackage{cite}
\usepackage[active]{srcltx}
\numberwithin{equation}{section}
\newtheorem{theorem}{Theorem}[section]
\newtheorem{lemma}{Lemma}[section]
\theoremstyle{remark}
\newtheorem{remark}{Remark}[section]
\theoremstyle{definition}

\newtheorem{example}{Example}[section]
\def\cT{\mathcal{T}}
\def\bold {\boldsymbol}

\def\J {\mathcal{J}}
\def\F {\mathcal{F}}
\def\U {\mathcal{U}}
\def\to{\rightarrow}
\def\eps {\varepsilon}
\def\A {\mathcal{A}}
\def\R {\mathbb{R}}

\def\p {\partial}

\newcommand{\dx}{\,{\rm d}x}

\newcommand{\dd}{\,{\rm d}}

\begin{document}
	
\title{Adaptive Computation of Elliptic Eigenvalue Topology Optimization with a Phase-Field Approach}
	
\author{Jing Li\footnote{School of Mathematical Sciences, East China Normal University, Shanghai 200241, China. (betterljing@163.com)} \and Yifeng Xu\footnote{Department of Mathematics \& Scientific Computing Key Laboratory of Shanghai Universities, Shanghai Normal University, Shanghai 200234, China. (corresponding author, yfxu@shnu.edu.cn)} \and Shengfeng Zhu\footnote{Key Laboratory of MEA (Ministry of Education) \& Shanghai Key Laboratory of PMMP, School of Mathematical Sciences, East China Normal University, Shanghai 200241, China. (sfzhu@math.ecnu.edu.cn)}}
	
\date{}
	
\maketitle
	
	%----------------------------------------------------------------
	
\begin{abstract}
In this paper, we study adaptive approximations of an elliptic eigenvalue optimization problem in a phase-field setting by a conforming finite element method. An adaptive algorithm is proposed and implemented in several two-dimensional numerical examples for illustration of efficiency and accuracy. Theoretical findings consist in the vanishing limit of a subsequence of estimators and the convergence of the relevant subsequence of adaptively-generated solutions to a solution to the continuous optimality system.
\end{abstract}

%\medskip

	
\medskip
	
	
	%----------------------------------------------------------------
	
\section{Introduction}\label{sec:intro}
	
Topology optimization features an optimal distribution of materials (optimum design) blue for certain objective subject to design constraints in a given domain \cite{BensoeSigmund:2003}. It has a wide range of applications in solid mechanics, fluid mechanics, acoustics, etc. In this field, due to its intriguing mathematical nature, significance has been attached to eigenvalue optimization, which, constrained by an eigenvalue problem, seeks the best possible design to minimize/maximize an objective involving the spectrum of the underlying differential operator. Moreover, this topic is useful in practical engineering. For instance, the optimization process helps to separate vibration frequencies and modes of a running engine in a machine from those of other components of the machine in order to avoid resonance (see e.g., \cite{AllaireAubryJouve, AllaireJouve:2001, ButtazzoDalMaso1993, GarckeHuttlKahleKnopfLaux:2023, GarckeHuttlKnopf:2022, Henrot:2006, HenrotPierre:2018, Pedersen:2000, QianHuZhu:2022, ZhangZhuLiuShen:2021} for theoretical analysis and numerical methods on eigenvalue optimization).
	
Though much progress has been made in the study of eigenvalue optimization, the cost of the numerical simulations remains high. As mentioned before, eigenvalue optimization may be viewed as a PDE-constrained optimization problem. This indicates that a lot of computer resources are devoted to solving the underlying eigenvalue problem by finite element methods in numerically searching for the optimal domain, especially in 3D. This paper focuses on the computing efficiency in numerical treatment of eigenvalue optimization.  
In particular, we plan to incorporate an adaptive strategy into the usual finite element method. Utilizing this idea, the authors have recently investigated an adaptive algorithm for efficient numerical solution of the minimum compliance problem in structural topology optimization \cite{JinLiXuZhu}. However, compared with linear elasticity for structural optimization, the eigenvalue problem is more challenging in numerical simulation and consequently complicates the resulting optimization model. To illustrate this point, our discussion starts with the following elliptic eigenvalue problem, one application of which lies in design of an inhomogeneous drum \cite{OsherSantosa:2001, ZhangZhuLiuShen:2021}: 
\begin{equation}\label{eigen_vp}
    \text{Find}~(\lambda, u) \in \R \times H_0^1(D)\setminus 0~\text{such that}~\int_{D} \bold{\nabla} u \cdot \bold{\nabla} v \dx + \alpha\int_{D}\chi_\Omega u v \dx = \lambda \int_{D} u v \dx    \quad \forall v \in  H_0^1(D), 
\end{equation}
where $D\subset\R^d$ $(d=2,3)$ is an open bounded polyhedral/polygonal domain,  $\Omega\subset\subset D$ is an open set,  $\chi_\Omega$ denotes the characteristic function of $\Omega$, and $\alpha>0$ is a given constant. Then the eigenvalue optimization problem of interest is
\begin{equation}\label{eigen-opt}
	\begin{aligned} \inf_{\Omega\in\mathcal{A}}\left\{\J(\Omega):=\Psi\left(\lambda_{i_1},\lambda_{i_2},\cdots,\lambda_{i_l}\right)\right\}\quad
\text{subject to different eigenvalues}~\lambda_{i_1},\lambda_{i_2},\cdots,\lambda_{i_l}~\text{of}~\eqref{eigen_vp},
	\end{aligned}
\end{equation}
where $\Psi:(\R_{>0})^l\to\R$ is a $C^1$ function with $\mathbb{R}_{>0}$ denoting the set of all positive real numbers,  $\A:=\{\text{open set}~\Omega\subset D~|~|\Omega| = V\}$ with $V\in (0,|D|)$ as a given volume  and $|\cdot|$ denoting the Lebesgue measure of some domain, and $l\in\mathbb{N}$. Since the linear operator associated with \eqref{eigen_vp} is self-adjoint, positive and compact from $L^2(D)$ to $L^2(D)$, the sequence $\{\lambda_i\}_{i\geq 1}$ of eigenvalues can be arranged in a nondecreasing order
$
0< \lambda_1 \leq \lambda_2 \leq \lambda_3 \leq \cdots \leq  \lambda_i \to \infty$ as $i\to \infty$
and corresponding eigenfunctions $\{u_i\}_{i\geq 1}\subset H_0^1(D)$ form an $L^2$-orthonormal basis of $L^2(D)$, where each eigenvalue is repeated according to its geometric multiplicity. Therefore, without loss of generality we assume that $\lambda_{i_1}, \lambda_{i_2},\cdots,\lambda_{i_l}$ are the $i_j$-th $(1\leq j\leq l)$ eigenvalue of \eqref{eigen_vp} each and put in a strictly increasing order.

When $\Psi\left(\lambda_{i_1},\lambda_{i_2},\cdots,\lambda_{i_l}\right)=\lambda_1$, it is proved in \cite{Chanillo:2000} that problem \eqref{eigen-opt} is mathematically equivalent to the composite drum problem \cite{Cox:1990} (see e.g., \cite{LiangLuYang:2015, ZhuWuLiu:2010,ZhangLiangCheng:2011} for numerical investigations). The model can be directly solved by the boundary variation technique \cite{BucurButtazzo:2005, Henrot:2006, SokolowskiZolesio:1992}. However, this approach results in high computational costs in numerical approximations using finite element methods, as the mesh must be continually updated with each boundary change. Moreover this classical technique takes no topological change into account. To overcome these two drawbacks, level set methods were used in \cite{HeKaoOsher:2007, OsherSantosa:2001, ZhangZhuLiuShen:2021}.

In this paper, we consider this problem from another point of view and adopt a surrogate model
using the phase-field approach. The main idea is the parametrization of a material distributed in the fixed domain $D$, also called the design domain, by a phase-field function $\phi\in [0,1]$ and then the variable domain $\Omega$ is approximated by $\phi\approx 1$. Now the resulting constrained minimization problem in the phase-field setting reads
	\begin{equation}\label{eigen-opt_phasefield}
		\begin{aligned}
			&\qquad\inf_{\phi\in\U}\left\{\J^{\eps}(\phi):=\Psi(\lambda_{i_1}^\eps,\lambda_{i_2}^\eps, \cdots, \lambda_{i_l}^\eps ) +\gamma \F_{\eps}(\phi)\right\}\\
			&\text{different eigenvalues}~\lambda_{i_1}^\eps,\lambda_{i_2}^\eps, \cdots, \lambda_{i_l}^\eps~\text{given by}
		\end{aligned}
	\end{equation}
	\begin{equation}\label{eigen_vp_phasefield}
		\int_{D} \bold{\nabla} w \cdot \bold{\nabla} v \dx +  \alpha\int_D \phi w v \dx = \lambda^\eps(\phi)\int_{D} w v \dx \quad \text{for}~(\lambda^\eps,w)\in \mathbb{R}\times H_0^1(D)\setminus \{0\},~\forall v \in H_0^1(D).
	\end{equation}
	Here $\eps$ and $\gamma$ are both positive and fixed constants,
	\[
	   \U:=\bigg\{\phi\in H^1(D)~|~\phi\in[0,1]~\text{a.e. in}~D, \int_{D}\phi\dx = V \bigg\},
	\]
	and $\F_{\eps}$ is a Ginzburg-Landau energy functional
	\begin{equation}\label{GL-energy}
		\F_{\eps}(\phi) : = \frac{\eps}{2}\|\bold{\nabla}\phi\|_{L^2(D)}^2 + \frac{1}{\eps} \int_{D}  f(\phi)  \dx,\quad \eps > 0
	\end{equation}
	with $f(\phi)$ being a double obstacle potential given by $\frac{1}{2}\phi(1-\phi)$ or a double well potential by $\frac{1}{4}\phi^2(1-\phi)^2$. Problem \eqref{eigen-opt_phasefield}-\eqref{eigen_vp_phasefield} may be viewed as a regularized formulation of \eqref{eigen-opt} as the regularization parameter $\gamma$ reflects the influence of Ginzburg-Landau energy on the configuration determined by the finite $l$ different eigenvalues $\lambda^\eps_{i_1},\cdots,\lambda^\eps_{i_l}$. As in \eqref{eigen_vp}, the eigenvalues $\{\lambda_i^\eps\}_{i\geq 1}$ of \eqref{eigen_vp_phasefield} are also a positive nondecreasing sequence
	\begin{equation}\label{eigenvalue_phase-field}
		0< \lambda^\eps_1 \leq \lambda^\eps_2 \leq \lambda^\eps_3 \leq \cdots \leq  \lambda^\eps_i \to \infty\quad  \text{as}~i\to \infty
	\end{equation}
	with the geometric multiplicity counted for each eigenvalue and the sequence of corresponding eigenfunctions $\{w_i\}_{i\geq1}\subset H_0^1(D)$ is an $L^2$-orthonormal basis of $L^2(D)$. As before, it is supposed that each of $\{\lambda_{i_j}^\eps\}_{j=1}^l$ in the functional $\J^\eps$ is the $i_j$-th eigenvalue of \eqref{eigen_vp_phasefield} and strictly ascends as $j$ increases. The existence of a minimizer to \eqref{eigen-opt_phasefield}-\eqref{eigen_vp_phasefield} can be proved as in \cite{GarckeHuttlKahleKnopfLaux:2023, GarckeHuttlKnopf:2022}. With $\phi^\ast\in\U$ and $\{(\lambda_{i_j}^{\eps,\ast},w^\ast_{i_j})\}_{j=1}^l\subset \mathbb{R} \times H_0^1(D)$ denoting the minimizer and the $l$ corresponding eigenpairs respectively, if  $\{\lambda_{i_j}^{\eps,\ast}\}_{j=1}^l$ are all simple, then the necessary optimality system \cite{GarckeHuttlKahleKnopfLaux:2023} reads 
	\begin{equation}\label{optsys_phase-field}
		\left\{
		\begin{array}{ll}
			\gamma\mathcal{F}_\eps'(\phi^\ast)(\phi-\phi^\ast) + \displaystyle\alpha\sum_{j=1}^l\frac{\p \Psi}{\p \lambda_{i_j}^\eps}(\lambda^{\eps,\ast}_{i_1},\lambda^{\eps,\ast}_{i_2},\cdots,\lambda^{\eps,\ast}_{i_l})\int_{D}(w_{i_j}^{\ast})^2 (\phi -\phi^\ast)  \dx  \geq 0 \quad \forall\phi \in \U,    \\[1ex]
			\displaystyle\int_{D} \bold{\nabla} w_{i_j}^\ast \cdot \bold{\nabla} v \dx +  \alpha \int_D \phi^\ast w^\ast_{i_j} v \dx = \lambda^{\eps,\ast}_{i_j}\int_{D} w_{i_j}^\ast v \dx \quad  \forall v \in  H^{1}_0(D),~\forall j \in \{1,2,\cdots, l\},
		\end{array}
		\right.
	\end{equation}
	where $$\mathcal{F}_\eps'(\phi^\ast)(\phi-\phi^\ast) = \eps\int_{D} \bold{\nabla}\phi^\ast\cdot\bold{\nabla}(\phi-\phi^\ast)\dd x+\frac{1}{\eps}\int_{D}f'(\phi^\ast)(\phi-\phi^\ast)\dd x.$$
	
With \eqref{eigen-opt_phasefield}-\eqref{eigen_vp_phasefield} in hand, the corresponding numerical implementation by finite element methods is performed over a fixed mesh and topological changes, even nucleation of holes in $D$, are allowed in the optimization process. This approach has been widely applied in structural optimization (see e.g., \cite{BurgerStainko:2006, WangZhou:2004}) since it was first used in \cite{BourdinChambolle:2003, BourdinChambolle:2006}. The reason for the success is that each sequence of minimizers for $\{\mathcal{F}_{\eps}\}_{\eps>0}$ has a cluster point $\chi_{\Omega^\ast}$ with respect to the $L^1(D)$-topology and the resulting set $\Omega^\ast$ has the minimal perimeter when $\eps\to0^+$  \cite{BloweyElliot:1991, Braides:2002, Modica:1987}. It is also evident that for small $\eps$ a minimizer $\phi^\ast$ to \eqref{eigen-opt_phasefield}-\eqref{eigen_vp_phasefield} takes values close to $0$ and $1$ mainly in $D$ and transits from 0 to 1 smoothly only in a thin layer with thickness of order $\eps$. As a result, a lot of degrees of freedom, especially distributed in the region where $\phi^\ast\approx0$ or $1$, are redundant for the improvement of numerical resolutions if one seeks finite element approximations of $\phi^\ast$ over uniformly refined meshes. In addition, the commonly used uniform mesh refinement strategy causes very high computational cost in solving the eigenvalue problem \eqref{eigen_vp_phasefield}. To remedy this situation, we incorporate the adaptive technique based on an \textit{a posteriori} error estimator and a local mesh refinement strategy into the numerical optimization process, which yields a so-called adaptive finite element method (AFEM). In contrast to traditional finite element methods, AFEM attempts to equi-distribute the error of numerical solutions at the expense of minimum degrees of freedoms. Therefore, the aim of this paper is the study of AFEM for \eqref{eigen-opt_phasefield}-\eqref{eigen_vp_phasefield}. We shall propose an adaptive algorithm to approximate three variables $(\phi, \lambda^\eps, w)$ over the same mesh and study its convergence. 
With an initial mesh $\cT_0$ given and a counter $k$ set to be zero, it performs successive loops of the form (see Algorithm \ref{alg_afem_eigenvalue} in Section \ref{sec:alg} for more detail)
	\begin{equation}\label{afem_loop}
		\textsf{SOLVE} \rightarrow \textsf{ESTIMATE} \rightarrow \textsf{MARK} \rightarrow \textsf{REFINE}.
	\end{equation}
	\textsf{SOLVE} outputs a discrete minimizer $\phi_k^\ast$ and $l$  corresponding discrete eigenpairs $\left\{\left(\lambda^{\eps,\ast}_{k, i_j}, w_{k,i_j}^\ast\right)\right\}_{j=1}^l$ over the current mesh $\cT_k$; \textsf{ESTIMATE} provides $l+1$ estimators of the residual type as in \cite{AinsworthOden:2000, Verfurth:2013}, written as $\eta_{k,0}$, $\eta_{k,1}$, $\cdots$, $\eta_{k,l}$ for brevity in our case, measuring the discretization error for optimality system \eqref{optsys_phase-field}; \textsf{MARK} selects some elements within $\cT_k$ to be refined according to some strategy using the $l+1$ estimators; \textsf{REFINE} generates a finer new mesh $\cT_{k+1}$ by locally refining all marked elements and their neighbors (if necessary) for conformity.
	
The major force driving our algorithm is clearly $l+1$ estimators $\eta_{k,j}$ $(0\leq j \leq l)$ in the module \textsf{ESTIMATE}. Each \textit{a posteriori} error estimator is a computable quantity involving the numerical minimizer, associated discrete eigenpairs, the coefficients in \eqref{eigen-opt_phasefield}-\eqref{eigen_vp_phasefield} and the mesh size (see \eqref{def:err-est_0}-\eqref{def:err-est_j}). We note that our attempt to prove the reliability of estimators in \eqref{def:err-est_0}-\eqref{def:err-est_j} or an upper error bound is thwarted by the nonlinear relation between $\phi$ and eigenpair $(\lambda^\eps,w)$ through \eqref{eigen_vp_phasefield}. Nevertheless, numerical results in Section \ref{sec:numerics} show that even if the thin layer for the phase transition is unknown, with the help of these estimators the free interface is identified correctly by performing more local mesh refinements as the algorithm proceeds. Further, they play an important role in the convergence analysis of Algorithm \ref{alg_afem_eigenvalue} (see Sections \ref{sec:estimator}-\ref{sec:conv}). Theorem \ref{thm:estimator->zero} asserts that each of $\{\eta_{k,j}\}_{k\geq0}$ $(0\leq j \leq l)$ has a subsequence convergent to zero. Then this allows us to prove that the relevant subsequence $\left\{\phi_{k_n}^\ast,(\lambda^{\eps,\ast}_{k_n, i_1}, w_{k_n,i_1}^\ast), \cdots, (\lambda^{\eps,\ast}_{k_n, i_l}, w_{k_n,i_l}^\ast)\right\}_{n \geq 0}$ of discrete solutions generated by Algorithm \ref{alg_afem_eigenvalue} converges strongly to a solution of optimality system \eqref{optsys_phase-field} (see Theorem \ref{thm:alg_conv}). It is worth pointing out that although some techniques in our convergence analysis are borrowed from the routine for forward problems \cite{CasconKreuzerNochetto:2008, GanterPraetorius:2022, GarauMorin:2011, GarauMorinZuppa:2009, MorinSiebertVeeser:2008, NochettoSiebertVeeser:2009, Siebert:2011}, our approach should be modified accordingly due to the nonlinear nature of optimization problem \eqref{eigen-opt_phasefield}-\eqref{eigen_vp_phasefield}. We begin with the discussion on how to obtain $\eta_{k,j}$ $(0\leq j\leq l)$ without reliability and efficiency from the optimality system (see \eqref{eqn:est_derive01}-\eqref{eqn:est_derive03}). The standard residual approach in the \textit{a posteriori} error analysis \cite{AinsworthOden:2000} is used here. More importantly, our argument in \eqref{eqn:est_derive01}-\eqref{eqn:est_derive03} motivates the subsequent convergence analysis although no reliability is made available. Then the key step of our analysis consists in a limiting constrained minimization problem posed over a subset defined as the limit of the sequence of discrete admissible sets given by Algorithm \ref{alg_afem_eigenvalue} (see \eqref{eigen-opt_phasefield_limit}-\eqref{eigen_vp_phasefield_limit}). After establishing some preliminary results for the limiting problem and the discrete problem (see Lemmas \ref{lem:Courant-Fisher_limit}-\ref{lem:eigen_weak_cont}), we shall prove the adaptively generated sequence of discrete solutions contains a subsequence strongly convergent to a minimizer to the limiting problem (see Theorem \ref{thm:alg_conv_limit}). This crucial finding, together with a separate-collective marking featuring a reasonable assumption \eqref{marking} in Algorithm \ref{alg_afem_eigenvalue}, and the approach in \cite{GanterPraetorius:2022} as well as arguments in \cite{MorinSiebertVeeser:2008, Siebert:2011} for forward problems, leads to the vanishing limit of  estimators up to a subsequence in Theorem \ref{thm:estimator->zero}, from which, \eqref{eqn:est_derive01}-\eqref{eqn:est_derive03} and the auxiliary convergence in Theorem \ref{thm:alg_conv_limit} again, we conclude the desired convergence in Theorem \ref{thm:alg_conv}. In contrast to our recent work \cite{JinLiXuZhu} on AFEM for structural optimization, where only D\"{o}rfler's marking strategy is considered, the analysis summarized above is quite different and more involved since a general yet practical condition \eqref{marking} is imposed in the module \textsf{MARK}.

	The layout of this paper is as follows. %In section 2, we establish the existence of a minimizer to \eqref{eigen-opt_phasefield}-\eqref{eigen_vp_phasefield} and present the related necessary optimality system. Then
	An adaptive finite element method is presented in Section \ref{sec:alg}, which is followed by several two-dimensional numerical examples showing the effectiveness of our adaptive algorithm in Section \ref{sec:numerics}. Sections \ref{sec:estimator}-\ref{sec:conv} are concerned with the theoretical aspect of our proposed algorithm. Section \ref{sec:estimator} relates to the motivation behind the $l+1$ estimators serving in the module \textsf{ESTIMATE} while an auxiliary limiting problem characterizing the limiting behavior of the algorithm is carefully examined in Section \ref{sec:limit}.  Finally, Section \ref{sec:conv} deals with the convergence analysis of the adaptive algorithm. Throughout the paper, we adopt standard notation for $L^p$ ($1\leq p\leq \infty$) spaces, Sobolev spaces and related (semi)-norms \cite{AdamsFournier:2003}. Moreover the lower-case letter $c$ with/without a subscript stands for a generic constant possibly different at each occurrence but independent of the mesh size and the function under consideration. %For any given integer $r\geq0$, $P_r(\omega)$ denotes the set of all polynomials of total degree no more than $r$ on a given bounded set $\omega\subset\mathbb{R}^d$ while
	
	%\section{Optimization problem in the phase-field %setting}\label{sec:phase-field}
	
	%In this section, we establish the well-posedness of \eqref{eigen-opt_phasefield}-\eqref{eigen_vp_phasefield} and present the related necessary optimality system.
	
	%Since two differential operators associated with \eqref{eigen_vp_phasefield} and \eqref{eigen_vp} are self-adjoint, compact and positive, two sequences of eigenvalues $\{\lambda_i\}_{i\geq 1}$ and $\{\lambda_i^\eps\}_{i\geq 1}$ can be both arranged in an nondecreasing order
	%\begin{equation}\label{eigenvalue}
	%  0< \lambda_1 \leq \lambda_2 \leq \lambda_3 \leq \cdots \leq  \lambda_i \to \infty,\quad
	%  0< \lambda^\eps_1 \leq \lambda^\eps_2 \leq \lambda^\eps_3 \leq \cdots \leq  \lambda^\eps_i \to \infty\quad  \text{as}~i\to \infty,
	%\end{equation}
	%and each sequence of corresponding eigenfunctions   $\{u_i\}_{i\geq 1}, \{w_i\}_{k\geq1}\subset H_0^1(D)$ forms an $L^2$-orthonormal basis of $L^2(D)$, where each eigenvalue is repeated according its geometric multiplicity.
	
	\section{Adaptive algorithm}\label{sec:alg}
	This section presents an adaptive finite element method of form \eqref{afem_loop} to approximate \eqref{eigen-opt_phasefield}-\eqref{eigen_vp_phasefield}. Let $\cT_0$ be a shape regular conforming triangulation of the domain
	$\overline{\Omega}$ into closed triangles/tetrahedra and $\mathbb{T}$ be the set of all possible conforming triangulations of $\overline{\Omega}$
	obtained from $\cT_0$ by the successive use of bisection \cite{Kossaczky:1995} \cite{NochettoSiebertVeeser:2009}. Then the set $\mathbb{T}$ is uniformly shape regular, i.e., the
shape regularity of any mesh $\mathcal{T}\in\mathbb{T}$ is bounded by a constant depending only on $\cT_0$ \cite{NochettoSiebertVeeser:2009} \cite{Traxler:1997}. Over any $\cT\in\mathbb{T}$, let $V_\cT$ be the usual $H^1$-conforming space of continuous piecewise polynomials of fixed order $r\geq1$. With $\U_\cT:=\U\cap V_\cT$ and $V_\cT^0:=V_\cT\cap H_0^1(D)$ coming into play, the finite element approximation of  \eqref{eigen-opt_phasefield}-\eqref{eigen_vp_phasefield} is formulated as
	
	\begin{equation}\label{eigen-opt_phasefield_disc}
		\begin{aligned}
			&\qquad\min_{\phi_\cT\in\U_\cT}\left\{\J_\cT^{\eps}(\phi_\cT):=\Psi(\lambda_{\cT,i_1}^\eps,\lambda_{\cT,i_2}^\eps, \cdots, \lambda_{\cT,i_l}^\eps ) +\gamma \F_{\eps}(\phi_\cT)\right\}\\
			&\text{subject to different eigenvalues}~\lambda_{\cT,i_1}^\eps,\lambda_{\cT,i_2}^\eps, \cdots, \lambda_{\cT,i_l}^\eps~\text{given by}
		\end{aligned}
	\end{equation}
	\begin{equation}\label{eigen_vp_phasefield_disc}
		%\left\{\begin{array}{l}
			\begin{aligned}
				\int_{D} \bold{\nabla} w_\cT & \cdot \bold{\nabla} v_\cT \dx +  \alpha \int_D \phi_\cT w_\cT v_\cT \dx = \lambda_\cT^\eps(\phi_\cT)\int_{D} w_\cT v_\cT \dx\quad  \forall v_\cT \in V^0_\cT \\  &\text{for}~(\lambda^\eps_\cT,w_\cT) \in \mathbb{R} \times V_\cT^0~\text{with} ~\|w_\cT\|_{L^2(D)} = 1.
			\end{aligned}
			%\end{array}
			%\right.
		\end{equation}
		Similar to the continuous case, \eqref{eigen_vp_phasefield_disc} has a finite sequence of eigenvalues counting geometric multiplicities
		\[
		0 < \lambda_{\cT,1}^\eps \leq \lambda_{\cT,2}^\eps \leq \lambda_{\cT,2}^\eps \leq \cdots \leq \lambda_{\cT,\dim V_{\cT}^0}^\eps
		\]
		and an $L^2(D)$-orthonormal sequence $\{w_{\cT,i}\}_{i=1}^{\dim V^0_{\cT}}$ of corresponding eigenfunctions. We assume that the $l$, not greater than $\mathrm{dim}V_{\cT}^0$, eigenvalues $\{\lambda_{\cT,i_j}^\eps\}_{j=1}^l$ are the $i_j$-th $(1\leq j\leq l)$ eigenvalue of \eqref{eigen_vp_phasefield_disc} each and enumerated in a strictly increasing order. The existence of a minimizer to \eqref{eigen-opt_phasefield_disc}-\eqref{eigen_vp_phasefield_disc} can be established by similar arguments for \eqref{eigen-opt_phasefield}-\eqref{eigen_vp_phasefield} (cf. \cite{GarckeHuttlKahleKnopfLaux:2023, GarckeHuttlKnopf:2022}) and Lemma \ref{lem:eigenvalue_disc_estimate} below.
		
Before getting down to our adaptive algorithm, we introduce   more notation. The collection of all edges/faces (resp. all interior edges/faces) on the boundary of elements in $\cT\in\mathbb{T}$ is denoted by $\mathcal{F}_\cT$ (resp. $\mathcal{F}_\cT(D)$). Associated to each $F\in \mathcal{F}_\cT$ is a fixed unit normal vector $\bold{n}_F$, which is specially taken as the unit outward normal on $\partial D$. Over $\cT$, we introduce a piecewise constant mesh-size function $h_\cT:\overline{D}
		\rightarrow \mathbb{R}_{>0}$ as
		\begin{equation}\label{meshsize_def}
			h_\cT |_T := h_T = |T|^{1/d}\quad \forall T\in\cT.
		\end{equation}
		Since $\cT$ is shape regular, $h_T$ is equivalent to the diameter of any $T\in\cT$. The union of elements in $\cT$ neighbouring some $T\in\cT$ (resp. sharing a common edge/face $F\in \mathcal{F}_\cT(D)$) is denoted by $\omega_\cT(T)$ (resp. $\omega_\cT(F)$). For the solution $\phi_{\cT}^\ast$ of \eqref{eigen-opt_phasefield_disc}-\eqref{eigen_vp_phasefield_disc} and the $l$ corresponding eigenpairs $\{(\lambda_{\cT,i_j}^{\eps,\ast}, w_{\cT,i_j}^\ast)\}_{j=1}^l$, we define element residuals on each $T\in\cT$ and jumps across each $F\in\mathcal{F}_\cT$ respectively by
		\begin{equation*}\label{def:err-est_element}
			\begin{aligned}
				&R_{0}(\phi_\cT^\ast)|_{T}:= \alpha\sum_{j=1}^{l}\frac{\p\Psi(\lambda_{\cT,i_1}^{\eps,\ast},\cdots,\lambda_{\cT,i_l}^{\eps,\ast})}{\p\lambda_{i_j}^{\eps}}(w_{\cT,i_j}^\ast)^2-\gamma\eps\Delta\phi^\ast_{\cT}+\frac{\gamma}{\eps}f'(\phi^\ast_\cT),\\
				& R_{j}(\phi_\cT^\ast,\lambda_{\cT,i_j}^{\eps,\ast},w_{\cT,i_j}^{\ast})|_{T}:= -\Delta w_{\cT,i_j}^\ast + \alpha\phi_{\cT}^\ast w_{\cT,i_j}^\ast - \lambda_{\cT,i_j}^{\eps,\ast}w_{\cT,i_j}^\ast\quad 1\leq j\leq l,
			\end{aligned}
		\end{equation*}
		\begin{equation*}\label{def:err-est_jump}
			J_0(\phi^\ast_\cT)|_{F}:= \left\{\begin{array}{lll}
				\gamma\eps[\bold{\nabla}\phi_{\cT}^{\ast}\cdot\boldsymbol{n}_{F}]\quad&
				F\in\mathcal{F}_{\cT}(D),\\ [1ex]
				\gamma\eps\bold{\nabla}\phi_{\cT}^{\ast}\cdot\boldsymbol{n}\quad&
				F\in \mathcal{F}_\cT \setminus \mathcal{F}_\cT(D),
			\end{array}\right. ~~~ J_{j}(w_{\cT,i_j}^{\ast})|_{F}:= [\bold{\nabla}w_{\cT,i_j}^{\ast}\cdot\boldsymbol{n}_{F}] \quad F\in \mathcal{F}_\cT(D)~1\leq j \leq l,
		\end{equation*}
		where $[\cdot]$ denotes the jumps across an inter-element face/edge. Then for any collection of elements
		$\mathcal{M}_{\cT}\subseteq\cT$, we define the following local error indicators and global error estimators
\begin{equation}\label{def:err-est_0}
    \eta_{\cT,0}^{2}(\phi^\ast_\cT;T):=h_T^2\|R_0(\phi^\ast_\cT)\|^2_{L^2(T)}+h_T\sum_{F\subset\p T}\|J_0(\phi^\ast_\cT)\|^2_{L^2(F)} \quad \eta^2_{\cT,0}(\phi^\ast_\cT;\mathcal{M}_{\cT})	:=\sum_{T\in\mathcal{M}_\cT}\eta_{\cT,0}^{2}(\phi^\ast_\cT;T),
\end{equation}
\begin{equation}\label{def:err-est_j}
    \begin{aligned}
        &\eta_{\cT,j}^{2}(\phi_\cT^\ast,\lambda_{\cT,i_j}^{\eps,\ast},w_{\cT,i_j}^{\ast};T):= h_T^2\|R_j(\phi_\cT^\ast,\lambda_{\cT,i_j}^{\eps,\ast},w_{\cT,i_j}^{\ast})\|^2_{L^2(T)}+h_T\sum_{F\subset\p T}\|J_j(w_{\cT,i_j}^{\ast})\|^2_{L^2(F)},\\ &\eta_{\cT,j}^{2}(\phi_\cT^\ast,\lambda_{\cT,i_j}^{\eps,\ast},w_{\cT,i_j}^{\ast};\mathcal{M}_{\cT})
        :=\sum_{T\in\mathcal{M}_{\cT}}\eta_{\cT,j}^{2}(\phi_\cT^\ast,\lambda_{\cT,i_j}^{\eps,\ast},w_{\cT,i_j}^{\ast};T)\quad 1\leq j\leq l.
	\end{aligned}
\end{equation}
When $\mathcal{M}_{\cT}=\cT$, the notation $\mathcal{M}_{\cT}$ will be suppressed and finite element functions involved in $\eta_{k,j}$ $(0 \leq j \leq l)$ are omitted if no confusion arises in the context. Now we are ready to present an adaptive algorithm for problem \eqref{eigen-opt_phasefield}-\eqref{eigen_vp_phasefield}. To avoid the abuse of notation, all dependence on a triangulation $\cT_k$ is indicated by the iteration number $k$ in the subscript, e.g., $V_{k}^0:=V_{\cT_k}^0$.

\begin{algorithm}
    \caption{AFEM for a shape design associated with elliptic eigenvalues as in \eqref{eigen-opt_phasefield}-\eqref{eigen_vp_phasefield}}\label{alg_afem_eigenvalue}%
    \LinesNumbered
    \KwIn{Specify an initial mesh $\cT_{0}$, fix $\eps>0$ and $\gamma>0$ and set the maximum number $K$ of refinements.}
    \KwOut{$\phi_k^\ast$ and $(\lambda_{k,i_j}^{			\eps,\ast},w_{k,i_j}^\ast)$ $(1\leq j\leq l)$.}
    \For {$k=0:K-1$} {
        {\textsf{(SOLVE)}} Solve \eqref{eigen-opt_phasefield_disc}-\eqref{eigen_vp_phasefield_disc} over $\cT_{k}$ for a minimizer $\phi_k^\ast\in \U_k$ and $l$ associated eigenpairs $(\lambda_{k,i_j}^{\eps,\ast},w_{k,i_j}^\ast)\in \mathbb{R}\times V_{k}^0$ $~(1\leq j\leq l)$\;
	
        {\textsf{(ESTIMATE)}} Compute $l+1$ error estimators $\eta_{k,0}(\phi_{k}^{\ast})$ and $\eta_{k,j}(\phi_{k}^{\ast},\lambda_{k,i_j}^{\eps,\ast},w_{k,i_j}^{\ast})$~$(1\leq j\leq l)$\;
	
        {\textsf{(MARK)}} Mark $l+1$ subsets $\mathcal{M}_{k}^{j}\subseteq\cT_{k}$~$(0 \leq j\leq l)$, each of which contains
    at least one element $T_{k}^j$~$(0 \leq j\leq l)$ holding the largest error indicator among all $\eta_{k,j}(T)$, i.e., there exists one element $T_{k}^j\in \mathcal{M}_k^j$ such that
    \begin{equation}\label{marking}		\eta_{k,j}(T_k^{j})=\max_{T\in\cT_k}\eta_{k,j}(T).%\widetilde{\eta}_{k}(\mathcal{M}_k)\geq\theta\widetilde{\eta}_k%=\max_{T\in\cT_k}\widetilde{\eta}_{k}(T).
	\end{equation}
        Then $\mathcal{M}_k:=\bigcup_{j=0}^{l}\mathcal{M}_{k}^j$\;
    %Let $\widetilde{\eta}_{k}:= \max_{0\leq j\leq l}\eta_{k,j}$, and mark a subset $\mathcal{M}_{k}\subseteq\cT_k$ such that $\mathcal{M}_k$ contains at least one element $T_{k}$ holding the largest error indicator among all $\widetilde{\eta}_{k}(T)$, i.e.,
    %\begin{equation}\label{marking}		\widetilde{\eta}_{k}(T_k)=\max_{T\in\cT_k}\widetilde{\eta}_{k}(T).%\widetilde{\eta}_{k}(\mathcal{M}_k)\geq\theta\widetilde{\eta}_k%=\max_{T\in\cT_k}\widetilde{\eta}_{k}(T).
    %\end{equation}
	
        {\textsf{(REFINE)}} Refine each element $T$ in $\mathcal{M}_{k}$ by bisection to get $\cT_{k+1}$\;

        Check the stopping criterion;
    }
\end{algorithm}

It should be noted that we use a separate-collective scheme in the module \textsf{MARK} of Algorithm \ref{alg_afem_eigenvalue}. First for each global estimator $\eta_{k,j}$, elements of $\cT_k$ are collected in a subset $\mathcal{M}_k^j$ according to a practical assumption \eqref{marking}, which is fulfilled by several popular marking strategies, e.g., the maximum strategy, the equi-distribution strategy, the modified equi-distribution strategy and D\"{o}rfler's strategy \cite{Siebert:2011}. Then the union of the $l+1$ resultant $\mathcal{M}_{k}^j$ yields the desired marked set $\mathcal{M}_k$ for refinement. For example, in our numerical experiments we implement D\"{o}rfler's strategy $l+1$ times to determine $l+1$ subsets $\mathcal{M}_{k}^j$ ($0\leq j\leq l$) separately, each of which satisfies 
\[
    \eta_{k,j}(\mathcal{M}_{k}^j) \geq \theta_j  \eta_{k,j}   
\]
with a given $\theta_j \in (0,1]$ and then the marked set for refinements is given by $\mathcal{M}_k:=\bigcup_{j=0}^{l}\mathcal{M}_{k}^j$.  A maximum iteration number $K$ is set in Algorithm \ref{alg_afem_eigenvalue} for termination while a stopping criterion, such as a prescribed tolerance for $\sum_{j=0}^l\eta^{2}_{k,j}$ or a given upper bound for the number of vertices over $\cT_k$, is also included in the adaptive loop. As mentioned in Section \ref{sec:intro}, even if we fail to establish any error bound in terms of $\eta_{k,j}$  ($0\leq j\leq l$), the relevant approach to obtain these computable quantities in Section \ref{sec:estimator} helps us to study the convergence of Algorithm \ref{alg_afem_eigenvalue} without $K$ and a stopping criterion in Sections \ref{sec:limit} and \ref{sec:conv}. More specifically, each $\{\eta_{k,j}\}_{k\geq0}$ ($0\leq j\leq l$) generated by Algorithm \ref{alg_afem_eigenvalue} is shown to have a subsequence converging to zero and the associated subsequence of discrete minimizers and $l$ discrete eigenpairs converge strongly to some solution of \eqref{optsys_phase-field}.

\section{Numerical example}\label{sec:numerics}
In this section, we validate the effectiveness and efficiency of the proposed adaptive phase-field algorithm for elliptic eigenvalue optimization through several numerical examples.  In Section \ref{subsect:num_detail}, we describe in detail the numerical realization in the module \textsf{SOLVE} of Algorithm \ref{alg_afem_eigenvalue} for discrete minimizers and associated eigenpairs of \eqref{eigen-opt_phasefield_disc}-\eqref{eigen_vp_phasefield_disc}
while numerical results are provided in Section \ref{subsect:num_discussion}.

\subsection{Detail of module SOLVE} \label{subsect:num_detail}
We employ an augmented Lagrangian method to relax the volume constraint. Let
\[\mathcal{L}(\phi)=\J^{\eps}(\phi)+ \mu G(\phi)+\frac{1}{2\beta} G^{2}(\phi),\]
where $\mu>0$ and $\beta>0$ are a Lagrange multiplier and a penalty parameter respectively and $G(\phi):=\int_{D} \phi \mathrm{d} x-V $. The volume constraint error is denoted by $|G(\phi)|$. The first-order necessary optimality condition for a minimizer $\phi^\ast$ can be derived by the 
G\^{a}teaux derivative as the following nonlinear system over $\widetilde{\U}:=\{\phi\in H^1(D)~|~\phi\in[0,1]~\text{a.e. in}~D\}$
\begin{equation}\label{LagOptCondtion}
    \begin{aligned}
        0&\leq \mathcal{L}^{\prime}(\phi^\ast)(\phi-\phi^\ast)=
        %\frac{\partial \mathcal{L}}{\partial \phi}(\phi^\ast-\phi) & =\J^{\eps\prime}(\phi^\ast)(\phi-\phi^\ast)+l G^{\prime}(\phi)+\frac{1}{\beta} G(\phi) G^{\prime}(\phi)  =
        \J^{\eps\prime}(\phi^\ast)(\phi-\phi^\ast)+\int_\Omega(\mu + \frac{1}{\beta}G(\phi^\ast))
        (\phi-\phi^\ast)\dx\quad \forall \phi\in \U\\%\left(\int_{D} \phi \dx -V\right).
        & = \gamma\mathcal{F}_\eps'(\phi^\ast)(\phi-\phi^\ast) + \displaystyle\alpha\sum_{j=1}^l\frac{\p \Psi}{\p \lambda_{i_j}^\eps}(\lambda^{\eps,\ast}_{i_1},\lambda^{\eps,\ast}_{i_2},\cdots,\lambda^{\eps,\ast}_{i_l})\int_{D}(w_{i_j}^{\ast})^2 (\phi -\phi^\ast)  \dx  \\
        &\quad+\int_\Omega(\mu + \frac{1}{\beta}G(\phi^\ast))
        (\phi-\phi^\ast)\dx\quad \forall \phi\in \widetilde{\U}.
    \end{aligned}
\end{equation}
% It is worth mentioning that all examples in the next subsection involve $\Psi(\lambda_{1}) = \pm\lambda_{1}$, $\Psi(\lambda_{1},\lambda_{2})=\lambda_{1}-\lambda_{2}$ or $\Psi(\lambda_{2},\lambda_{3})=\lambda_{2}-\lambda_3$ respectively. And the sensitivity of $\Psi$ reads \cite{QianHuZhu:2022}
%        \begin{equation}\label{ObjSensitivity}
%        \Psi^{\prime}(\phi)=\left\{\begin{array}{ll}
%        \alpha w_{1}^{2} \chi_{\Omega}^{\prime}(\phi) & \text { if } \Psi=\pm\lambda_{1}, \\
%        % -\alpha w_{1}^{2} \chi_{\Omega}^{\prime}(\phi) & \text { if } \Psi=-\lambda_{1}, \\
%        \alpha\left(w_{i}^{2}-w_{i+1}^{2}\right) \chi_{\Omega}^{\prime}(\phi) & \text { if } \Psi=\lambda_{i}-\lambda_{i+1},~~i=1,2,
%        \end{array}\right.
%        \end{equation}
%        where $\Omega=\{x \in D \mid \chi_\Omega(x)=1\}$ and $\chi_{\Omega}$ is the indicator function of $\Omega$. %$(\lambda_{i},w_{i})$ is the $i$-th eigenpair.
%         In order to enhance the smoothness of the descent direction, we solve a gradient smoothing flow: find  $\Theta \in H_{0}^{1}(D)$  such that
%         \begin{equation}\label{GradientFlow}
%        \int_{D}(\tilde{\alpha} \nabla \Theta \cdot \nabla \mathcal{V}+\Theta \mathcal{V}) \dx=\int_{D} \mathcal{L}^{\prime}(\phi) \mathcal{V} \dx \quad \forall \mathcal{V} \in H_{0}^{1}(D),
%        \end{equation}
%        where  $\tilde{\alpha}>0$ is specified in each example. And let $\Theta$ replace $\mathcal{L}^{\prime}(\phi)$.

The three variables $\phi$, $\mu$ and $\beta$ are calculated in an alternating manner. For $n=0,1,\ldots,N-1$ where $N$ is specified in Table \ref{tab:hyper} below for each example, we denote by $\phi_n$, $\mu_n$, $\beta_n$ and $(\lambda_{i,n}, w_{i,n})$ ($i=1,2,3$) the approximations at the $n$-th optimization iteration over a fixed mesh generated by Algorithm \ref{alg_afem_eigenvalue}. Our numerical scheme in \textsf{SOLVE} first produces the eigenpair $(\lambda_{i,n}, w_{i,n})$, followed by $\phi_{n+1}$, $\mu_{n+1}$ and $\beta_n$ in succession. As mentioned above, the computation of the minimizing phase-field function depends on the first-order necessary condition. To address this, we introduce the following gradient flow for $\phi=\phi(t,\bm{x})$ %phase-field equation is formulated as
\begin{equation}\label{PhaseFieldEquation1}
    \left\{\begin{array}{ll}
        \dfrac{\partial \phi}{\partial t}=\kappa \Delta \phi-F^{\prime}(\phi) & \text { in }(0, \infty) \times D, \\ [1ex]
        \dfrac{\partial \phi}{\partial \boldsymbol{n}}=0 & \text { on } \partial D, \\ [1ex]
        \phi(0, \cdot)=\phi_{0} & \text { in } D,
    \end{array}\right.
\end{equation}
where $t>0$ is the pseudo-time, $\bold{n}$ is the unit outward normal on $\p D$, and $\phi_0\in L^2(D)$ is the initial guess. Here, $\kappa=\gamma\eps$ is a coefficient of the diffusion term. 
% is a coefficient of the diffusion term. we set
We set $F(\phi)=W(\boldsymbol{x}) w(\phi)+\mathcal{G}(\boldsymbol{x}) g(\phi)$ as in \cite{AkihiroShinjiMitsuru:2010}, which is given by
\begin{equation}\label{PhaseFieldEquation1_aux}
    \begin{array}{l}
        w(\phi)=\phi^{2}(1-\phi)^{2}, \quad g(\phi)=\phi^{3}\left(6 \phi^{2}-15 \phi+10\right),\quad
        W(\boldsymbol{x})=\dfrac{1}{4}, \quad \mathcal{G}(\boldsymbol{x})=\Tilde{\gamma} \dfrac{\mathcal{L}^{\prime}\left(\phi_{0}\right)}{\left\|\mathcal{L}^{\prime}\left(\phi_{0}\right)\right\|_{L^{2}(D)}}.
    \end{array}
\end{equation}
The parameter $\tilde{\gamma}>0$, related to the scaling of $\mathcal{L}^{\prime}\left(\phi_{0}\right)$, is specified in Table \ref{tab:hyper} below for each example. With a fixed $T>0$ introduced, then \eqref{PhaseFieldEquation1} can be rewritten as
\begin{equation}\label{PhaseFieldEquation2}
   \begin{aligned}
       \frac{\partial \phi}{\partial t} & =\kappa \Delta \phi-\frac{\partial}{\partial \phi}(W(\boldsymbol{x}) w(\phi)+\mathcal{G}(\boldsymbol{x}) g(\phi))  =\kappa \Delta \phi-\left(\frac{1}{4} w^{\prime}(\phi)+\Tilde{\gamma} \frac{\mathcal{L}^{\prime}\left(\phi_{0}\right)}{\left\|\mathcal{L}^{\prime}\left(\phi_{0}\right)\right\|_{L^{2}(D)}} g^{\prime}(\phi)\right) \\
       & =\kappa \Delta \phi+\phi(1-\phi)\left[\phi-\frac{1}{2}-30 \Tilde{\gamma} \frac{\mathcal{L}^{\prime}\left(\phi_{0}\right)}{\left\|\mathcal{L}^{\prime}\left(\phi_{0}\right)\right\|_{L^{2}(D)}}(1-\phi) \phi\right] \quad \text { in }(0, T] \times D .
   \end{aligned}
\end{equation}
% Numerical treatment of \eqref{PhaseFieldEquation2} relies on a finite difference method for the time variable and a finite element method for the space variables respectively.
For better stability in temporal discretization, we adopt a semi-implicit variational formulation of \eqref{PhaseFieldEquation2} \cite{QianHuZhu:2022} to find $\phi^{[m+1]}_{n} \in H^{1}(D)$ such that for any $\varphi \in H^{1}(D)$, 
\begin{equation}\label{PhaseFieldEquation3}
    \tau^{-1}(\phi^{[m+1]}_n-\phi^{[m]}_n, \varphi) + \kappa (\nabla \phi^{[m+1]}_n , \nabla \varphi) = \left\{\begin{array}{ll}
    \left(\phi^{[m+1]}_n\left(1 - \phi^{[m]}_n \right) r\left(\phi^{[m]}_n\right), \varphi \right), & r\left( \phi^{[m]}_n \right) \leqslant 0, \\
    \left( \phi^{[m]}_n \left( 1 - \phi^{[m+1]}_n \right) r\left( \phi^{[m]}_n \right), \varphi \right), & r\left( \phi^{[m]}_n \right)>0,
    \end{array}\right.
\end{equation} %scheme for the diffusion term and the reaction term , namely %based on implicit and semi-implicit discretizations , respectively.
%\begin{equation}
%    \frac{\phi^{[n, \ell+1]}-\phi^{[n, \ell]}}{\tau}=\kappa \Delta \phi^{[n, \ell+1]}+\left\{\begin{array}{ll}
%    \phi^{[n, \ell+1]}\left(1-\phi^{[n, \ell]}\right) r\left(\phi^{[n, \ell]}\right) & \text { for } r\left(\phi^{[n, \ell]}\right) \leqslant 0 \\
%    \phi^{[n, \ell]}\left(1-\phi^{[n, \ell+1]}\right) r\left(\phi^{[n, \ell]}\right) & \text { for } r\left(\phi^{[n, \ell]}\right)>0,
%    \end{array}\right.
%\end{equation}
where $\tau>0$ is a suitable time step, $(\cdot,\cdot)$ denotes the $L^{2}$ inner product, $m=0,1,\ldots, M-1$ with $M=[T/\tau]$,
$\phi^{[0]}_0=\phi_0$ and
\begin{equation*}
    r\left(\phi^{[m]}_n \right)=\phi^{[m]}_n-\frac{1}{2}-30 \Tilde{\gamma} \frac{\mathcal{L}^{\prime}(\phi_{n}^{[0]})}{\|\mathcal{L}^{\prime}(\phi_{n}^{[0]})\|_{L^{2}(D)}}\left(1-\phi^{[m]}_n\right) \phi^{[m]}_n.
\end{equation*}
%while to deal with both regular and irregular design regions finite element method is utilized in spatial discretization instead of finite difference method. Then 
In the numerical experiments, the above semi-implicit scheme is run $M$ (specified in Table \ref{tab:hyper} below for each example) times for each fixed $n$, with a conforming linear finite element method over each fixed mesh generated by Algorithm \ref{alg_afem_eigenvalue} utilized in spatial discretization.
%and $\phi^{[n, \ell]}$ denotes the $\ell$th inner iteration to solve the phase field equation (\ref{PhaseFieldEquation2}) at the $n$th optimization iteration. The semi-implicit variational formulation reads \cite{QianHuZhu:2022}: given $\phi^{[n, \ell]}$ $(\ell=0,1, \ldots, L-1)$, find  $\phi^{[n, \ell+1]} \in H^{1}(D)$ such that for any $\varphi \in H^{1}(D)$
% \begin{equation}\label{PhaseFieldVariational}
%     \left\{\begin{array}{l}
%     \left(\phi^{[n, \ell+1]}, \varphi\right)+\tau \kappa\left(\nabla \phi^{[n, \ell+1]}, \nabla \varphi\right)-\tau\left(\left(1-\phi^{[n, \ell]}\right) r\left(\phi^{[n, \ell]}\right) \phi^{[n, \ell+1]}, \varphi\right)=\left(\phi^{[n, \ell]}, \varphi\right)~
%     \text { for } r\left(\phi^{[n, \ell]}\right) \leqslant 0, \\
%     \left(\phi^{[n, \ell+1]}, \varphi\right)+\tau \kappa\left(\nabla \phi^{[n, \ell+1]}, \nabla \varphi\right)+\tau\left(\phi^{[n, \ell]} r\left(\phi^{[n, \ell]}\right) \phi^{[n, \ell+1]}, \varphi\right)=\left(\phi^{[n, \ell]}, \varphi\right)+\tau\left(\phi^{[n, \ell]} r\left(\phi^{[n, \ell]}\right), \varphi\right) \\
%     \text { for } r\left(\phi^{[n, \ell]}\right)>0,
%     \end{array}\right.
% \end{equation}
%where $(\cdot,\cdot)$ denotes the $L^{2}$ inner product.
        %Numerically we observe that at some positions the density function does not preserve the box constraint  [0,1].
Noting that the box condition $[0,1]$ might not be satisfied by $\phi^{[M]}_n$ in numerical computations, we resort to  %The box constraint on $\phi$ is enforced by
a projection step
\[
    \phi^{[m+1]}_{n} \leftarrow \min \left\{\max \left\{0, \phi^{[m+1]}_n \right\}, 1\right\}, 
\] and set $\phi_{n+1}^{[0]} = \phi_{n+1} = \phi_{n}^{[M]}$ after $M$ steps of \eqref{PhaseFieldEquation3}. The Lagrange multiplier $\mu$ is updated by %the scheme:
\[
    \mu_{n+1}=\mu_{n}+\frac{1}{\beta_{n+1}}G\left(\phi_{n+1}\right).
\]
The value of $\mu$ increases (resp. decreases) as the current volume is greater (resp. less) than the target volume. During iterations, the penalty parameter is updated by
$\beta_{n+1}=\xi \beta_{n}$ with $0<\xi<1$. Finally, we choose \[
    \tau=\zeta \max_{T\in\mathcal{T}_k}h_T/\max\lvert\mathcal{L}^{\prime}(\phi_n^{[0]}) \rvert
\]
%where $h:=\max_{T\in\mathcal{T}_k}h_T$ is the current maximal mesh size
with $\zeta>0$ to be specified in each example. We summarize the whole process on a fixed mesh in Algorithm \ref{Alga}, which specifically implements \textsf{SOLVE} of Algorithm \ref{alg_afem_eigenvalue}.

\begin{algorithm}
	\caption{Phase-field topology optimization on a fixed mesh}\label{Alga}
        \LinesNumbered
    \KwIn{parameters $\mu_0$, $\beta_0$, initial guess $\phi_0^{[0]}$, and integers $N$ and $M$.}
    \KwOut{the phase-field function $\phi_N$.}% and the associated eigenpairs}
		\For {$n=0:N-1$} {
            Solve the eigenvalue problem \eqref{eigen_vp_phasefield_disc} for $(\lambda_{n,i},w_{n,i})$\; %$l$\;
		  Compute the gradient \eqref{LagOptCondtion}\;
		  %\STATE Adaptive mesh moving based on residual error estimate\;
		  \For{$m=0:M-1$}{
                Semi-implicit updating step $\phi_n^{[m+1]}$ via \eqref{PhaseFieldEquation3}\;
		      Project $\phi_n^{[m+1]}$ to $[0,1]$\;
		      %\STATE  Check the stopping criterion.
          }
          Update $\phi_{n+1}^{[0]} = \phi_{n+1}= \phi_{n}^{[M]}$\;
          Update the Lagrange multiplier $\mu_{n+1}$\;
          Update the penalty parameter $\beta_{n+1}$\;
        }
\end{algorithm}

% \begin{algorithm}
% 	\caption{Phase-field topology optimization on fixed mesh}\label{Alga}
% 	\begin{algorithmic}[1]
% 		\STATE Set $\ell$, $\Tilde{\gamma}$, $\Tilde{\alpha}$, $\xi$, initial guess $\phi_0^{[0]}$, and integers $N$ and $M$
% 		\FOR {$n=1,\ldots,N$}
% 		\STATE Solve the state problem $\bm{w}_n$;
% 		\STATE Update the Lagrange multiplier $\ell_n$; %$l$\;
% 		\STATE Compute the gradient;
% 		%\STATE Adaptive mesh moving based on residual error estimate\;
% 		\FOR{$m=1,\ldots,M$}
% 		\STATE  Semi-implicit updating step $\phi_n^{[m]}$;
% 		\STATE  Projection $\phi_n^{[m]}$ to $[0,1]$;
% 		%\STATE  Check the stopping criterion.
% 		\ENDFOR
% 		\STATE Update the penalty parameter $\beta_{n+1}$;
% 		\ENDFOR
% 	\end{algorithmic}
% \end{algorithm}

\subsection{Discussion on numerical result} \label{subsect:num_discussion}
In this subsection, we evaluate Algorithm \ref{alg_afem_eigenvalue} through seven examples, with all numerical results implemented using MATLAB R2020a in a laptop with an Intel(R) Core(TM) i5-7500 CPU and 8 GB memory. The volume fraction is denoted by $C={V}/{|D|}$. For all experiments, $V_k$ is consistently chosen as the conforming linear finite element space. D\"{o}rfler's strategy is employed to represent the general assumption \eqref{marking} for each error estimator $\eta_{k,j}$ ($0\leq j\leq l$ and $l=1$ or $2$ in each example), which selects a subset $\mathcal{M}_k^j\subseteq\mathcal{T}_k$ such that
\begin{equation}\label{marking_num}
    \eta_{k,j}(\mathcal{M}_k^j)\geq \theta_j \eta_{k,j}\quad  0\leq j\leq l
\end{equation}
and defines $\mathcal{M}_k:=\bigcup_{j=0}^{l}\mathcal{M}_k^j$
with each threshold $\theta_j\in (0,1]$ to be specified in each example. In the design configurations below, red and blue stand for phase-field function values of $1$ and $0$, respectively. %In order to show the advantage of our proposed adaptive algorithm, we also list results over uniformly refined meshes in Examples \ref{example1}-\ref{example7} for comparison.

\begin{example}\label{example1}
Maximization of the first eigenvalue $\lambda_{1}$ with $\alpha=10.21$ and $C=0.5$ in a unit square $D=\left(0,1\right)^{2}$.
% In a unit square $D=\left(0,1\right)^{2}$, (a) minimization of the first eigenvalue $\lambda_{1}$ with $\alpha=10.21$ and $C=0.5$; (b) maximization of the first eigenvalue $\lambda_{1}$ with $\alpha=10.21$ and $C=0.5$. We note that part (a) is equivalent to a two-density drum problem \cite{Chanillo:2000}.

%Set $K=6$, $N=20$, $\theta_{1}=0.8$, $\beta=100$, $l=-5$, $\gamma=5\times 10^{-3}$, $\epsilon=5\times10^{-2}$ and $C=0.5$.

%The uniform result is refined by four times, and the eigenvalue problem is solved twenty-four times on each mesh.
\end{example}

\begin{example}\label{example2}
Minimization of $\lambda_{1}$ in a unit disk $D=\{ \left(x_1,x_2\right)\mid x_1^{2}+x_2^{2} < 1 \}$ with $\alpha=1$ and $C=0.5$. The optimal shape is a concentric annulus $\{\left(x_1,x_2\right) \mid R^2 < x_{1}^2+x_{2} < 1\}$ with $R=1/\sqrt{2}$ \cite{Chanillo:2000, LiangLuYang:2015}. %is given by $\pi$the value of $R>0$ is chosen to satisfy the volume constraint with $C=0.5$. 

%Set $\theta_{1}=0.7$, $\theta_{2}=0.2$, $C=0.5$, $\beta=50$, $\Tilde{\gamma}=20$, $\gamma=10^{-3}$, $\epsilon=10^{-2}$ and $\zeta=0.1$. For the initial design, choose $\phi_{0}=\frac{1}{2}(1-\cos(2\pi x_1)\cos(2\pi x_2))$, cf. Fig. \ref{fig:InitialPhi}(c).
\end{example}

\begin{example}\label{example4}
Minimization of $\lambda_{1}$ with $\alpha=1$ and $C=0.75$ in an L-shaped domain $\left(0,2\right)^{2}\backslash \left[1,2\right]^{2}$.
% (b) maximization of $\lambda_1$ with $\alpha=1$ and $C=0.2$.

%with $\theta_{1}=0.8$ and $\theta_{2}=0.3$. Set refinement parameter $K=6$ and $N=20$. Let $\epsilon=8\times10^{-3}$, $\gamma=10^{-3}$ and $\zeta=0.1$.
\end{example}

\begin{figure}[htbp]
	\centering
	\begin{tabular}{ccccc}
	\includegraphics[width=1in]{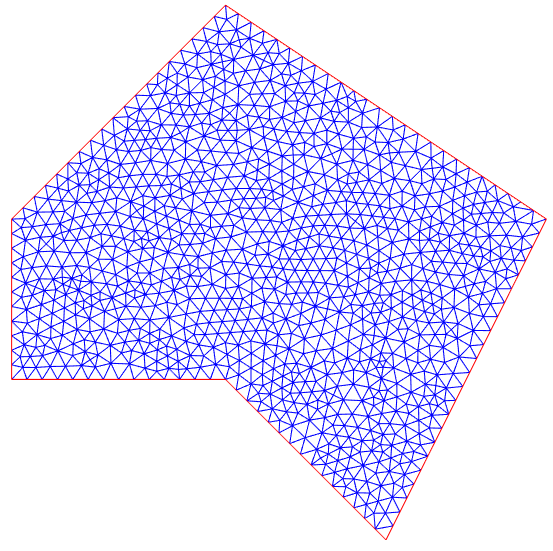}
        &\includegraphics[width=1in]{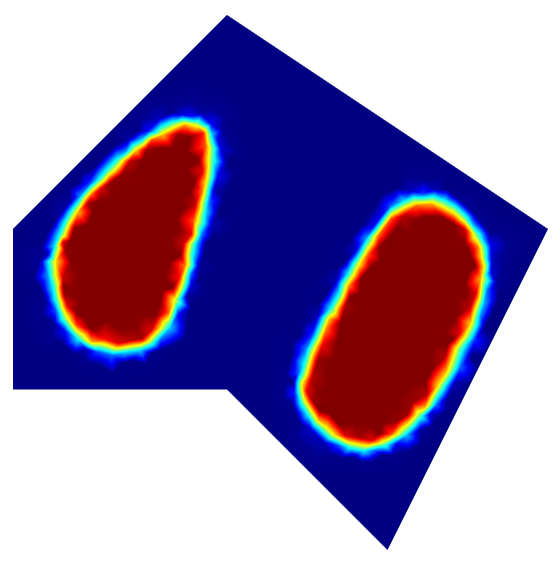}
        &\includegraphics[width=1in]{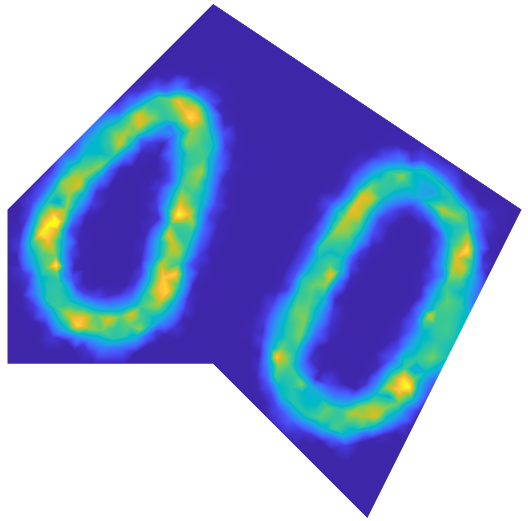}
        &\includegraphics[width=1in]{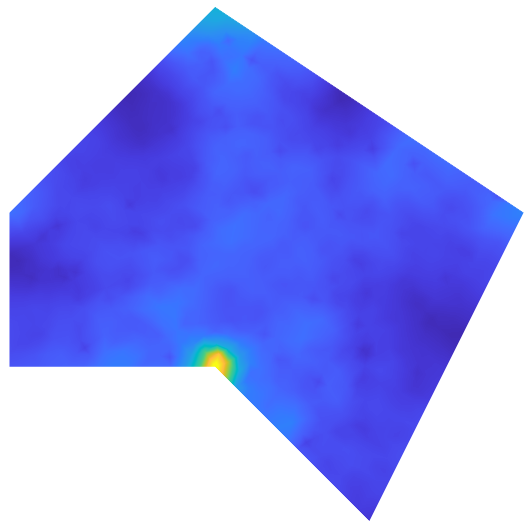}
        &\includegraphics[width=1in]{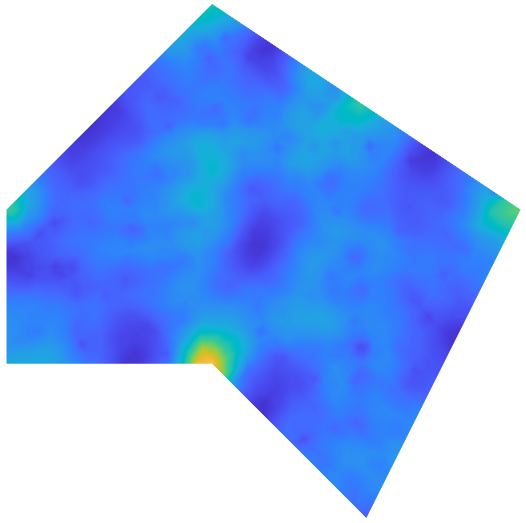}\\
        \includegraphics[width=1in]{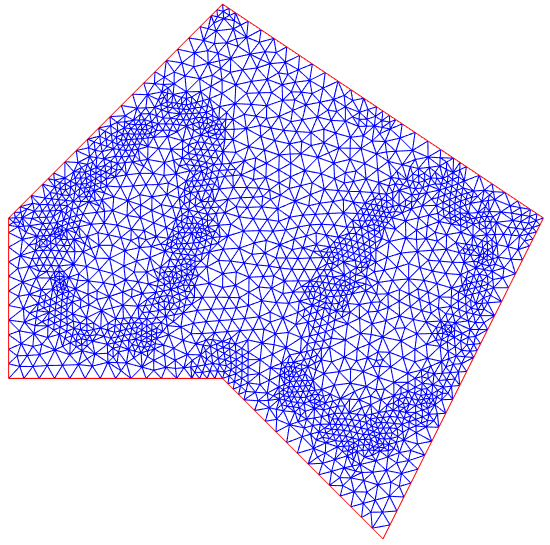}
        &\includegraphics[width=1in]{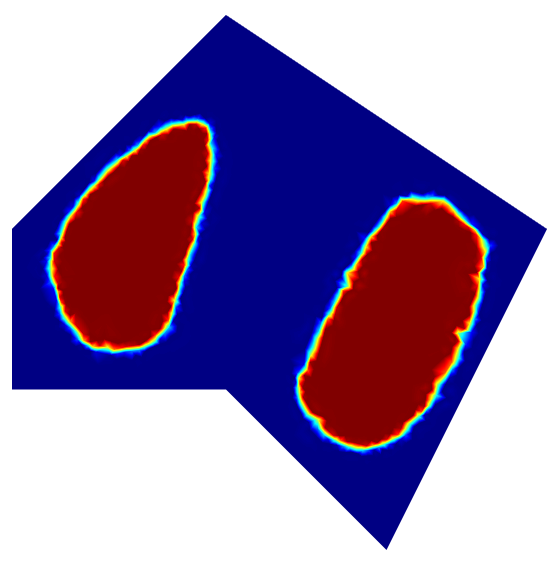}
        &\includegraphics[width=1in]{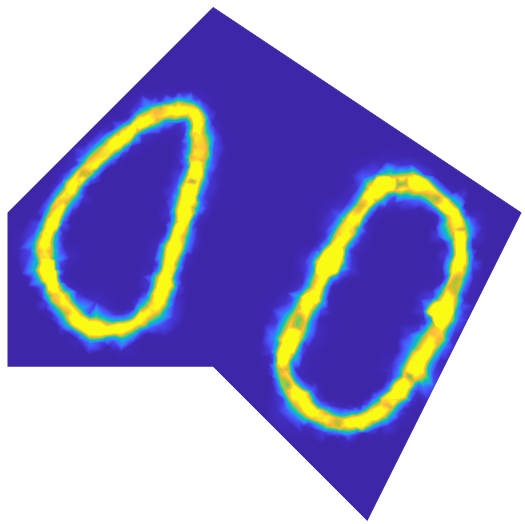}
        &\includegraphics[width=1in]{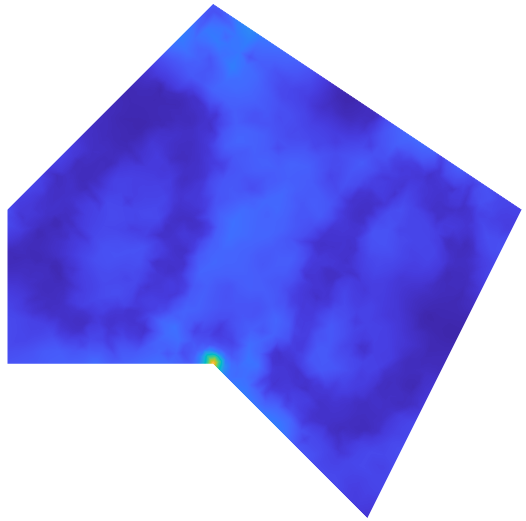}
        &\includegraphics[width=1in]{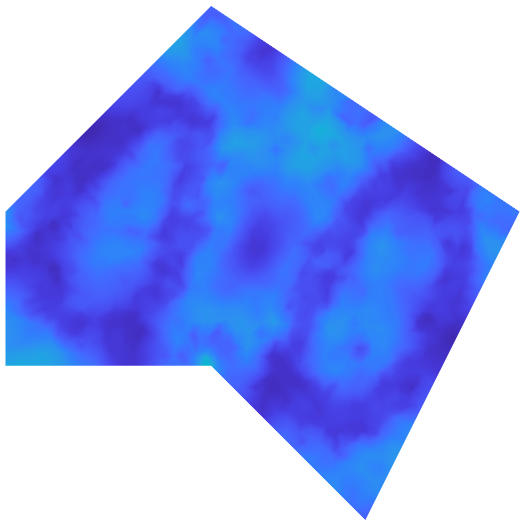}\\
        % \includegraphics[width=1in]{IrrPMax212TH.png}
        % &\includegraphics[width=1in]{IrrPMax21FinalDesign3TH.png}
        % &\includegraphics[width=1in]{IrrPMax213th0EstimatorS.png}
        % &\includegraphics[width=1in]{IrrPMax213th1EstimatorS.png}
        % &\includegraphics[width=1in]{IrrPMax213th2EstimatorS.png}\\
        \includegraphics[width=1in]{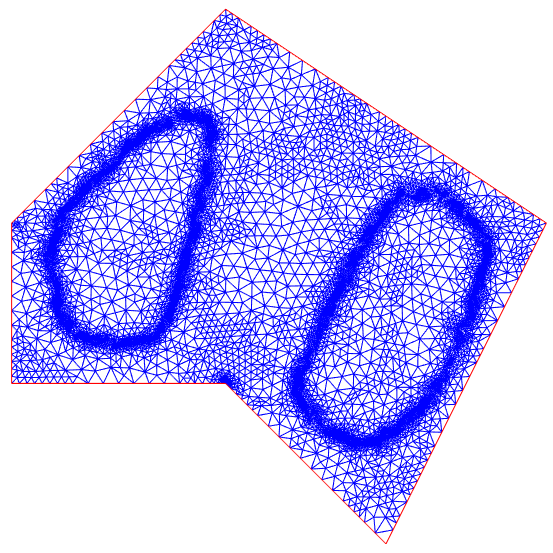}	
	  &\includegraphics[width=1in]{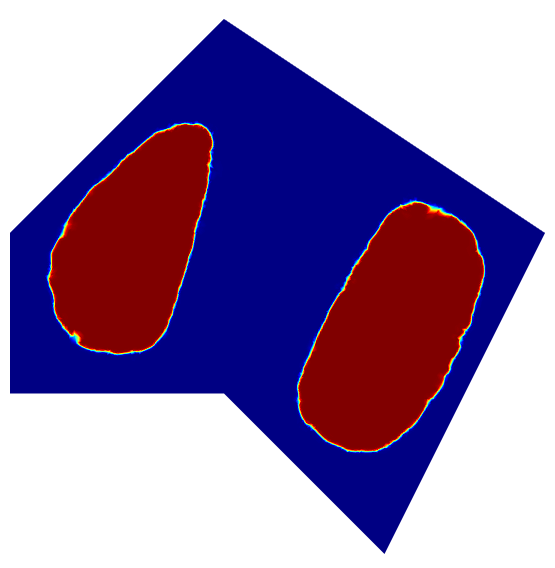}
        &\includegraphics[width=1in]{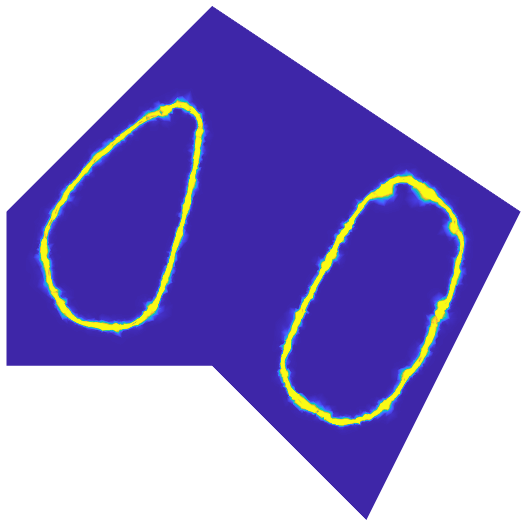}
        &\includegraphics[width=1in]{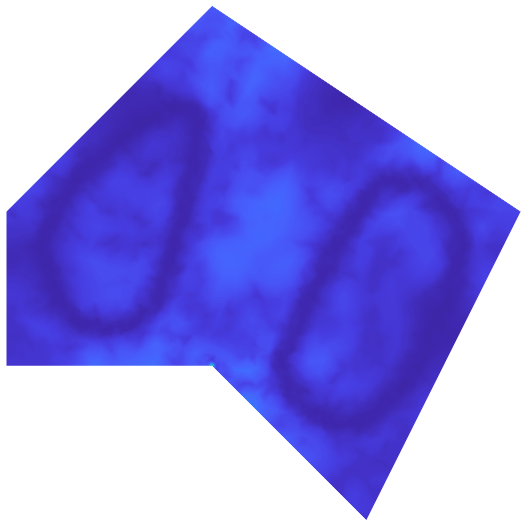}
        &\includegraphics[width=1in]{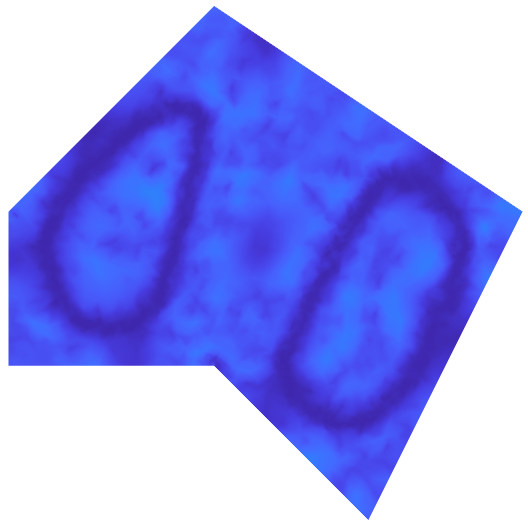}\\
        &
        &\includegraphics[width=1in]{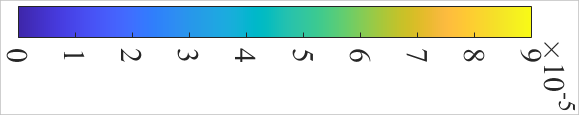}
        &\includegraphics[width=1in]{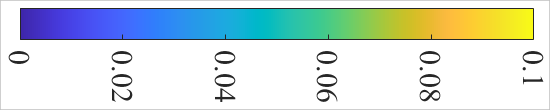}
        &\includegraphics[width=1in]{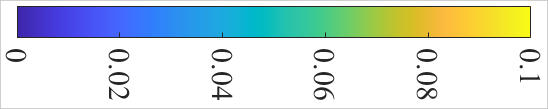}
	\end{tabular}
\caption{Evolution of adaptive mesh levels $k = 0,1,3$ for Example \ref{example7}, with the number of vertices on each mesh being 1087, 2133 and 7681. The 2nd, 3rd, 4th and 5th columns show optimized designs $\phi_k^\ast$ with a numerical optimal shape highlighted in red over each $\cT_k$ and three error indicators $\eta_{k,0}$, $\eta_{k,1}$ and $\eta_{k,2}$ respectively.}\label{AdaptiveExp7Max21}
\end{figure}

\begin{example}\label{example7}
Maximization of $\lambda_{2}-\lambda_{1}$ with $\alpha=1$ and $C=0.35$ in a nonconvex polygonal domain (see Figure \ref{AdaptiveExp7Max21}).
% In a nonconvex polygonal domain (see e.g., Figure \ref{AdaptiveExp7Min1}), (a) minimization of $\lambda_{1}$ with $\alpha=1$ and $C=0.65$; (b) maximization of $\lambda_{2}-\lambda_{1}$ with $\alpha=1$ and $C=0.35$.
%Set $K=4$, $N=30$, $\epsilon=10^{-2}$, $\theta_{1}=0.8$ and $\theta_{2}=\theta_{3}=0.1$ for minimizing $\lambda_{1}$ and maximizing $\lambda_{2}-\lambda_{1}$.
\end{example}

\begin{example}\label{example3}
Maximization of $\lambda_{3}-\lambda_{2}$ with $\alpha=9$ and $C=0.35$ in a rectangle $\left(-1,2\right)\times\left(-0.5,1\right)$.
% In a rectangle $\left(-1,2\right)\times\left(-0.5,1\right)$, (a) maximization of $\lambda_{2}-\lambda_{1}$ with $\alpha=9$ and $C=0.45$; (b) maximization of $\lambda_{3}-\lambda_{2}$ with $\alpha=9$ and $C=0.35$.
%and $\theta_{1}=0.8$.
\end{example}

\begin{example}\label{example5}
In a dumbbell-shaped domain $D_{\mathcal{H}}=B(-1, 0) \cup((-1, 1) \times(-\mathcal{H}, \mathcal{H})) \cup B(1,0)$, where $B(\bold{y})$ represents a unit disk centered at $\bold{y}$, (a) minimization of $\lambda_{1}$ with the handle width $2\mathcal{H}=1.2$, $\alpha=1$ and $C=0.5$; (b) minimization of $\lambda_{1}$ with the handle width $2\mathcal{H}=0.6$, $\alpha=1$ and $C=0.85$.
%		\[
%		  B(\bold{y})=\left\{\bold{x} \in \mathbb{R}^{2}:|\mathbf{x}-\mathbf{y}|<1\right\} .
%		\]
		%When $\mathcal{H}=0.6$  and  $0.3$, $\bar{a}=1$.
%And set $K=4$, $N=20$ and $\epsilon=10^{-2}$.
\end{example}

\begin{example}\label{example6}
In an annulus $D=\left\{(x_1, x_2) \mid r_{1} < \sqrt{x_1^{2}+x_2^{2}} < r_{2}\right\}$, (a) minimization of $\lambda_1$ with $r_{1}=1$, $r_{2}=3.5$, $\alpha=10$ and $C=0.85$; (b) minimization of $\lambda_1$ with $r_{1}=1$, $r_{2}=1.5$, $\alpha=1$ and $C=0.64$.
% In an annulus $D=\left\{(x_1, x_2) \mid r_{1} < \sqrt{x_1^{2}+x_2^{2}} < r_{2}\right\}$, (a) maximization  of $\lambda_1$ with $r_{1}=1$, $r_{2}=3.5$, $\alpha=1$ and $C=0.36$; (b) minimization of $\lambda_1$ with $r_{1}=1$, $r_{2}=3.5$, $\alpha=10$ and $C=0.85$; (c) minimization of $\lambda_1$ with $r_{1}=1$, $r_{2}=1.5$, $\alpha=1$ and $C=0.64$.
%And we denote $\eta=r_{2}/r_{1}$. We consider $r_{1}=1$, $r_{2}=3.5$ and  $r_{1}=1$, $r_{2}=1.5$, respectively.
\end{example}

%We set the Lagrangian multiplier $\mu=0$, $\alpha=1$, $\Tilde{\gamma}=20$, $\Tilde{\alpha}=10^{-4}$ and $\xi=0.9$ unless otherwise specified. The phase field equation is solved $M=10$ times on each mesh. And $C$ indicates volume fraction ${V}/{|D|}$.
%The uniform refinement is implemented three times, and the eigenvalue problem is solved in the same way with the corresponding adaptive algorithm on each mesh unless otherwise stated.

\begin{figure}[htbp]
    \centering\setlength{\tabcolsep}{0pt}
    \begin{tabular}{cccc}
		\includegraphics[width=1in]{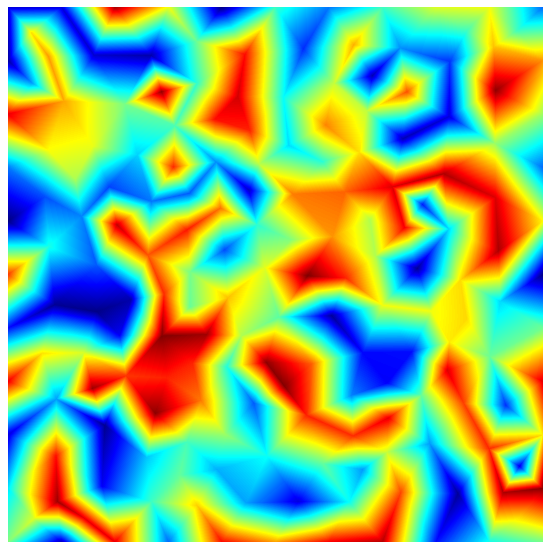} \hspace{0.1cm}
        &\includegraphics[width=1in]{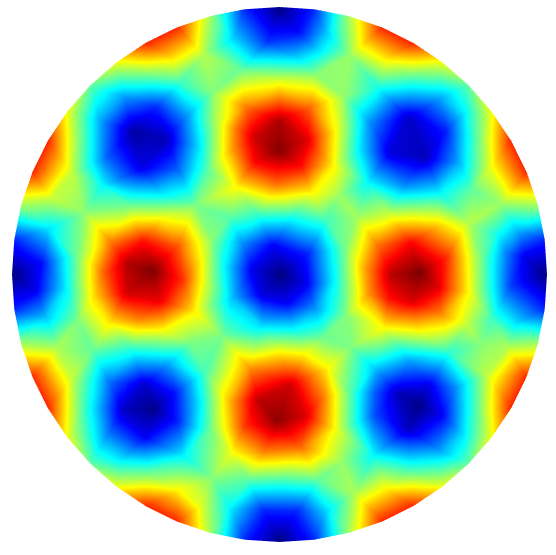} \hspace{0.1cm}
        % &\includegraphics[width=1in]{IrrPMin1Initialphi.png} \hspace{0.1cm}
        &\includegraphics[width=1in]{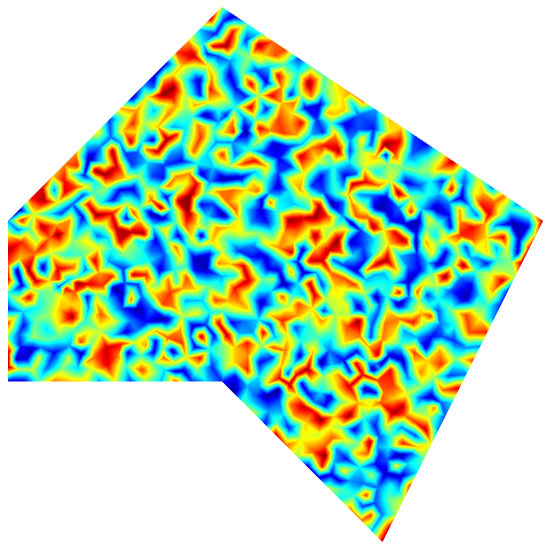} \hspace{0.05cm}
        &\includegraphics[width=0.185in]{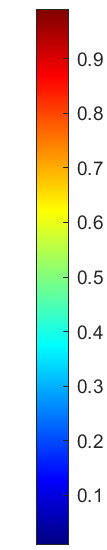}
    \end{tabular}
    \caption{Initial phase-field functions for Examples \ref{example1} (left), \ref{example2} (middle) and \ref{example7} (right).} %(a) for Example \ref{example1}: minimizing $\lambda_{1}$. (b) for Example \ref{example1}: maximizing $\lambda_{1}$. (c) for Example \ref{example2}. (d) for Example \ref{example7}: minimizing $\lambda_{1}$. (e) for Example \ref{example7}: maximizing $\lambda_{2}-\lambda_{1}$.}
    \label{fig:InitialPhi}
\end{figure}

The hyperparameters in Algorithm \ref{alg_afem_eigenvalue} and Algorithm \ref{Alga} for the numerical experiments are summarized in Table \ref{tab:hyper}. In the table, $\eps$ and $\gamma$ are positive constants appearing in \eqref{eigen-opt_phasefield}. $K$ refers to the maximum iteration number in Algorithm \ref{alg_afem_eigenvalue}. Each marking parameter $\theta_j$ ($0\leq j\leq 2$) is utilized in \eqref{marking_num} as a separate D\"{o}rfler's strategy in the module \textsf{MARK}. $N$ denotes the number of iterations to update the eigenpair and $M$ represents the number of pseudo-time steps in the inner iteration for the gradient flow to update the phase-field function. The parameters $\mu_0$ and $\beta_0$ are the initial Lagrange multiplier and the initial penalty weight for the augmented Lagrangian term while the factors $\tilde{\gamma}$, $\xi$ and $\zeta$ help to update $\mathcal{G}(\bold{x})$ in \eqref{PhaseFieldEquation1_aux}, the penalty weight $\beta_{n+1}$ and the phase-field pseudo-time step size $\tau$, respectively.
%The parameter $\theta_{2}$ equals to $0.1$ in Example \ref{example3}. 
Initial phase-field functions are set to be random values in $[0,1]$ in Examples \ref{example1} and \ref{example7}, and $\frac{1}{2}\cos(1-\cos(2 \pi x_1)\cos(2 \pi x_2))$ in Example \ref{example2} as displayed in Figure \ref{fig:InitialPhi}. In the remaining four examples, initial guesses are constants, i.e., $0.75$ in Example \ref{example4}, $0.35$ in Example \ref{example3}, $0.8$ and $0.5$ in Example \ref{example5}, and $0.5$ in Example \ref{example6}.

\begin{table}[htbp]
\centering
\begin{threeparttable}
\caption{Hyper-parameters for \eqref{eigen-opt_phasefield}, Algorithm \ref{alg_afem_eigenvalue} and Algorithm \ref{Alga} used in the experiments. \label{tab:hyper}}
\begin{tabular}{|c|cccccccccccc|}
\toprule
Example & $\eps$ & $\gamma$ & $K$ & $\theta_0$ & $\theta_1$ & $\theta_2$ & $(N,M)$ & $\mu_0$ & $\beta_0$ & $\tilde{\gamma}$ & $\xi$ & $\zeta$ \\
\midrule
    % \ref{example1}(a) & $5\times10^{-2}$ & $ 5\times10^{-3}$ & $6$ & $0.8$ & $0.1$ &  & $(20,10)$ & $-5$ & $100$ & $20$  & $0.9$ & $2$ \\
    \ref{example1} & $5\times10^{-2}$ & $5\times 10^{-3}$ & $6$ & $0.8$ & $0.3$ & & $(20,10)$ & $-5$ & $100$ & $20$ & $0.7$ & $2.5$\\
 \hline
    \ref{example2}  & $10^{-2}$ & $10^{-3}$ & $6$ & $0.7$ & $0.2$ &  & $(20,10)$ & $0$ & $50$ & $20$  & $0.9$ & $0.1$ \\
\hline
    \ref{example4} & $8\times10^{-3}$ & $10^{-3}$ & $6$ & $0.8$ & $0.3$ &  & $(20,10)$ & $0$ & $100$ & $20$  & $0.9$ & $0.1$ \\
    % \ref{example4}(b) & $8\times10^{-3}$ & $10^{-3}$ & $6$ & $0.8$ & $0.3$ &  & $(20,10)$ & $0$ & $80$ & $20$  & $0.9$ & $0.1$ \\
\hline
    %\ref{example7} & $10^{-2}$ & $10^{-3}$ & $4$ & $0.8$ & $0.1$ &  & $(30,10)$ & $0$ & $100$ & $20$  & $0.9$ & $0.05$ \\
    \ref{example7} & $10^{-2}$ & $5 \times 10^{-3}$ & $4$ & $0.8$ & $0.1$ & $0.1$  & $(30,10)$ & $0$ & $300$ & $20$  & $0.9$ & $0.05$ \\
\hline
    % \ref{example3}(a) &  $10^{-2}$ & $5\times10^{-3}$ & $6$ & $0.8$ & $0.2$ & $0.2$  & $(30,10)$ & $0$ & $80$ & $20$  & $0.9$ & $5$ \\
    \ref{example3}  & $5 \times 10^{-3}$ & $ 10^{-3}$ & $4$ & $0.8$ & $0.1$ & $0.1$  & $(40,10)$ & $0$ & $110$ & $20$  & $0.9$ & $0.6$ \\
\hline
    \ref{example5}(a) & $10^{-2}$ & $10^{-3}$ & $4$ & $0.7$ & $0.1$ & & $(20,10)$ & $0$ & $80$ & $20$  & $0.9$ & $0.03$ \\
    \ref{example5}(b) & $10^{-2}$ & $5\times 10^{-3}$ & $4$ & $0.8$ & $0.3$ & & $(20,10)$ & $0$ & $150$ & $20$  & $0.9$ & $0.03$ \\
\hline
    \ref{example6}(a) &  $10^{-2}$ & $5 \times 10^{-3}$ & $4$ & $0.6$ & $0.1$ & & $(10,10)$ & $0$ & $200$ & $20$  & $0.9$ & $0.01$ \\
    % \ref{example6}(b)  & $10^{-2}$ & $5 \times 10^{-3}$ & $4$ & $0.8$ & $0.1$ & & $(30,10)$ & $0$ & $60$ & $20$  & $0.9$ & $0.4$ \\
    \ref{example6}(b)  & $8 \times 10^{-3}$ & $5 \times 10^{-3}$ & $4$ & $0.8$ & $0.1$ & & $(20,10)$ & $0$ & $60$ & $20$  & $0.9$ & $0.05$ \\
\bottomrule
\end{tabular}
\end{threeparttable}
\end{table}

\begin{figure}[htbp]
	\centering\setlength{\tabcolsep}{0pt}
	\begin{tabular}{ccccccc}
		\includegraphics[width=1in]{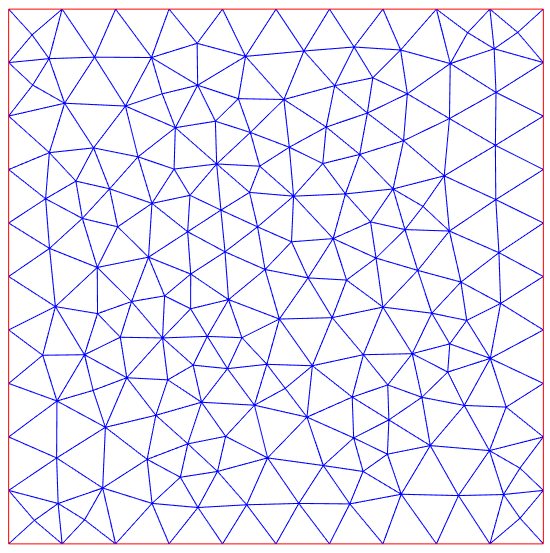}  
        &\includegraphics[width=1in]{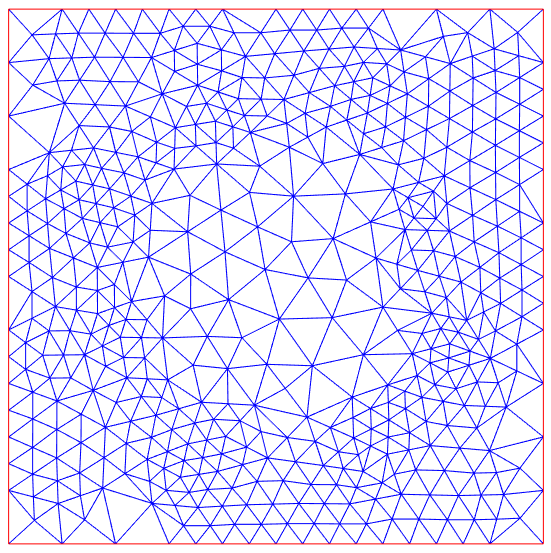}
		% &\includegraphics[width=1in]{SquareMax13TH.png}
		&\includegraphics[width=1in]{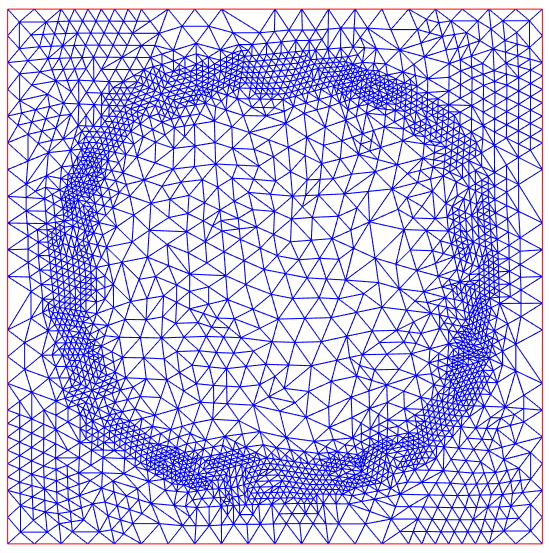}
		% &\includegraphics[width=1in]{SquareMax15TH.png} 
        &\includegraphics[width=1in]{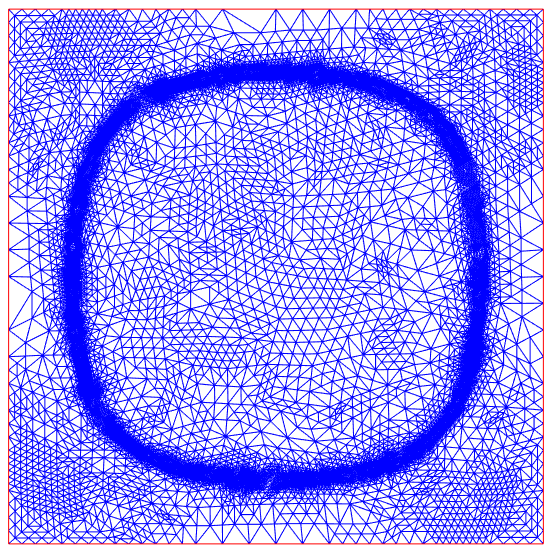}
        \\		
        \includegraphics[width=1in]{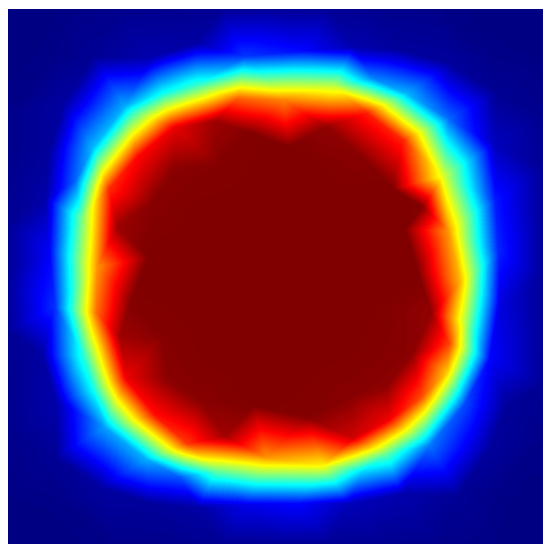}  
        &\includegraphics[width=1in]{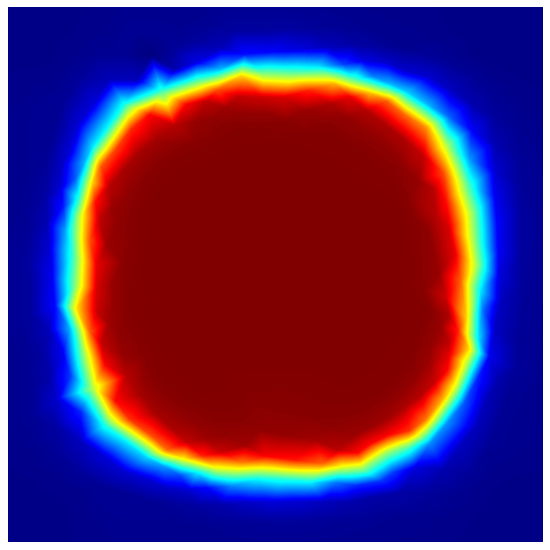}
        % &\includegraphics[width=1in]{SquareMax1FinalDesign3TH.png} 
        &\includegraphics[width=1in]{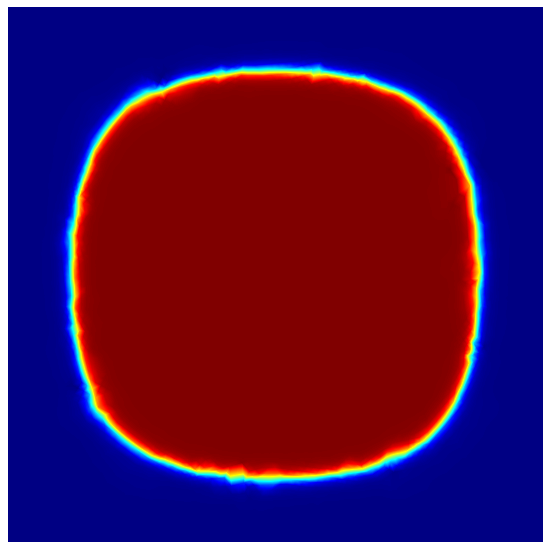}
        % &\includegraphics[width=1in]{SquareMax1FinalDesign5TH.png}
        &\includegraphics[width=1in]{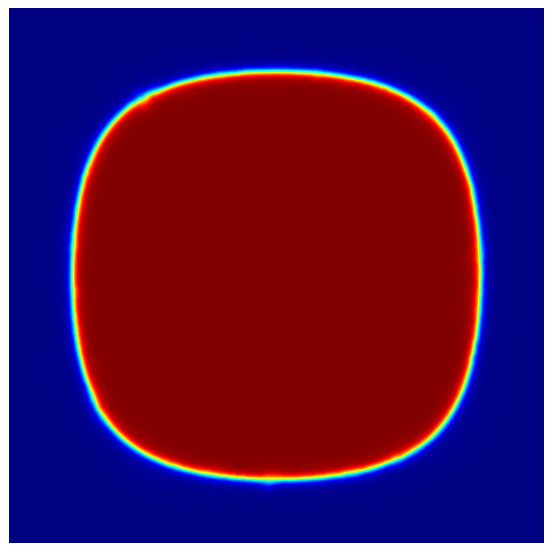}
        \\
        \includegraphics[width=1in]{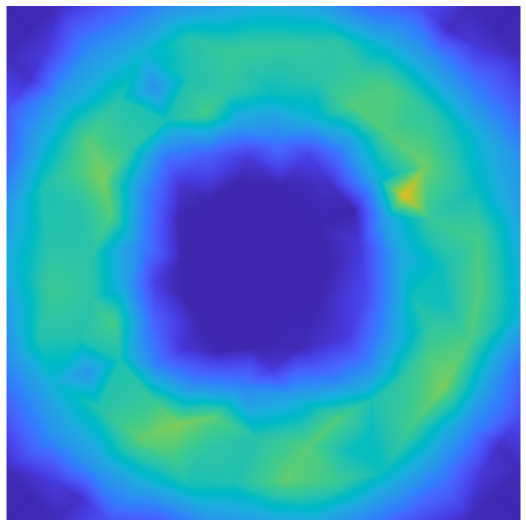}
        &\includegraphics[width=1in]{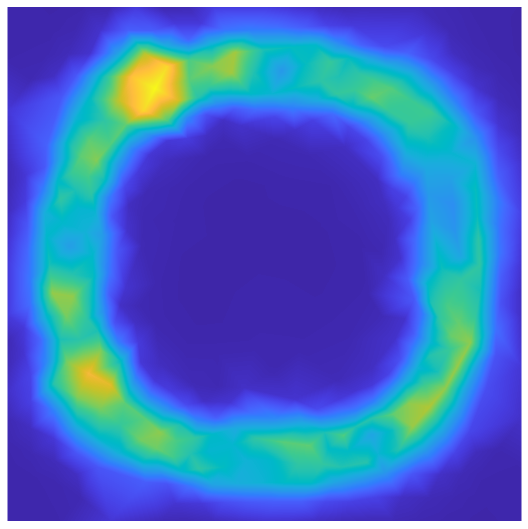}
        % &\includegraphics[width=1in]{SquareMax13th1EstimatorS.png}
        &\includegraphics[width=1in]{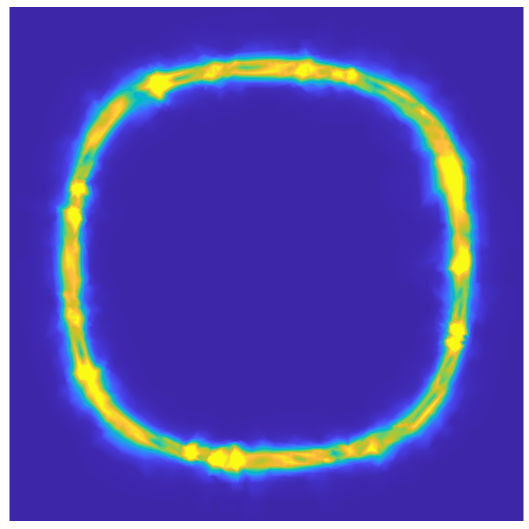}
        % &\includegraphics[width=1in]{SquareMax15th1EstimatorS.png}
        &\includegraphics[width=1in]{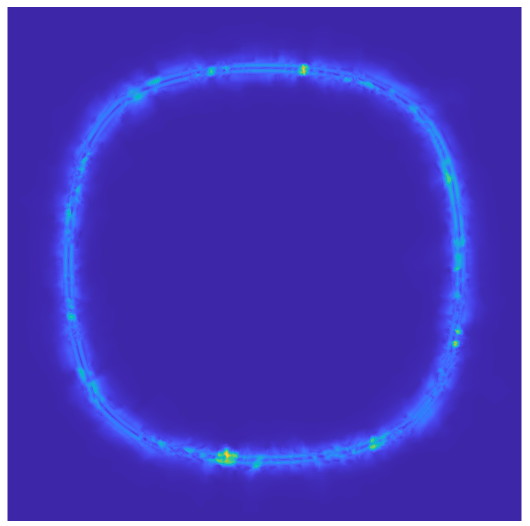}
        &\includegraphics[width=0.2in]{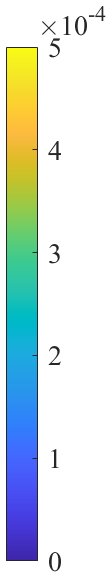}
        \\
        \includegraphics[width=1in]{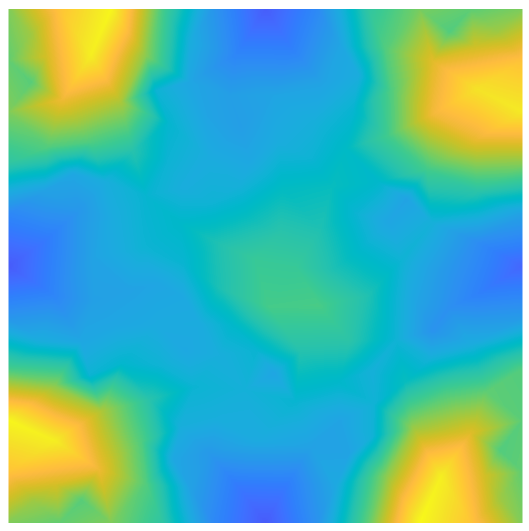}
        &\includegraphics[width=1in]{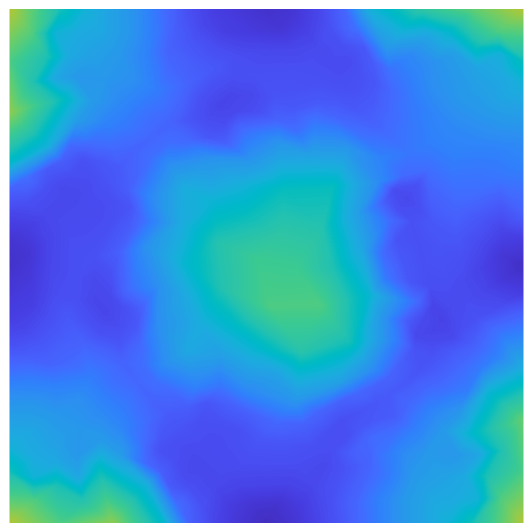}
        % &\includegraphics[width=1in]{SquareMax13th2EstimatorS.png}
        &\includegraphics[width=1in]{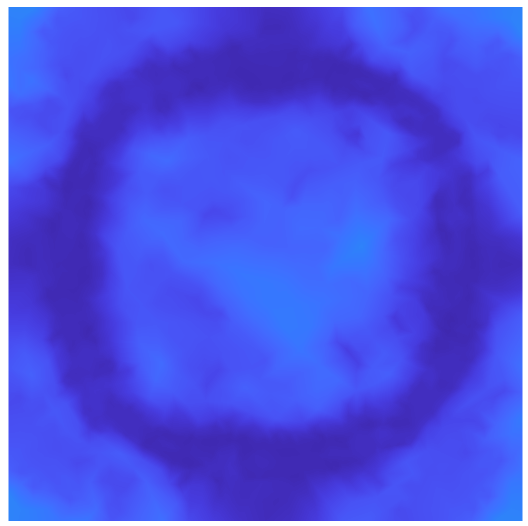}
        % &\includegraphics[width=1in]{SquareMax15th2EstimatorS.png}
        &\includegraphics[width=1in]{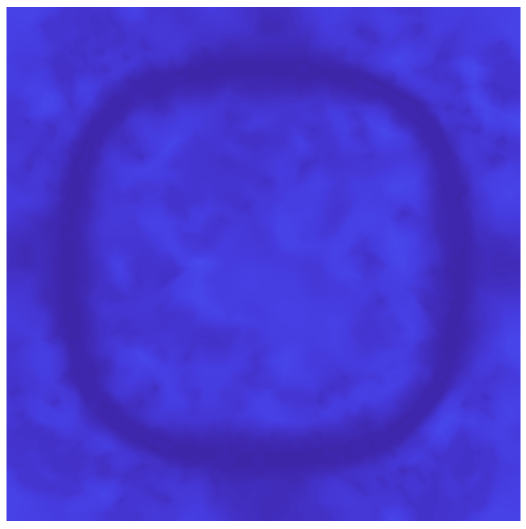}
        &\includegraphics[width=0.17in]{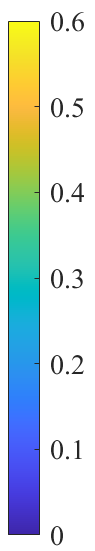}
	\end{tabular}
\caption{Evolution of adaptive mesh levels $k=0,1,3,5$ for Example \ref{example1}, with the number of vertices on each mesh being 185, 494, 2702 and 13331. The 2nd, 3rd and 4th rows show optimized designs $\phi_k^\ast$ with a numerical optimal shape  highlighted in red over each $\cT_k$, two error indicators $\eta_{k,0}$ and $\eta_{k,1}$ respectively.}\label{AdaptiveExp1Max1}
\end{figure}

\begin{figure}[htbp]
	\centering\setlength{\tabcolsep}{0pt}
	\begin{tabular}{ccccc}
		\includegraphics[width=1in]{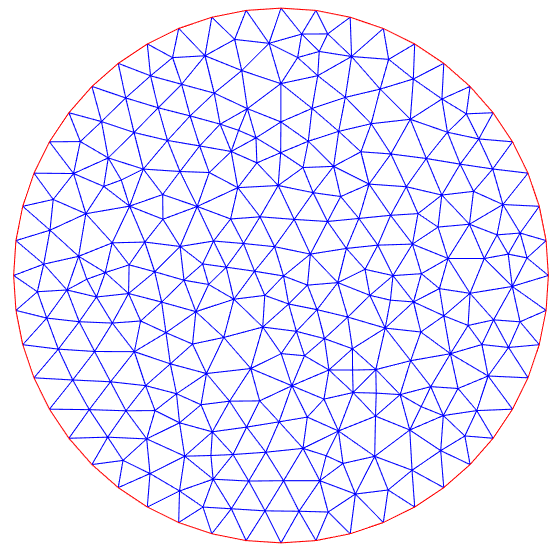}  
        &\includegraphics[width=1in]{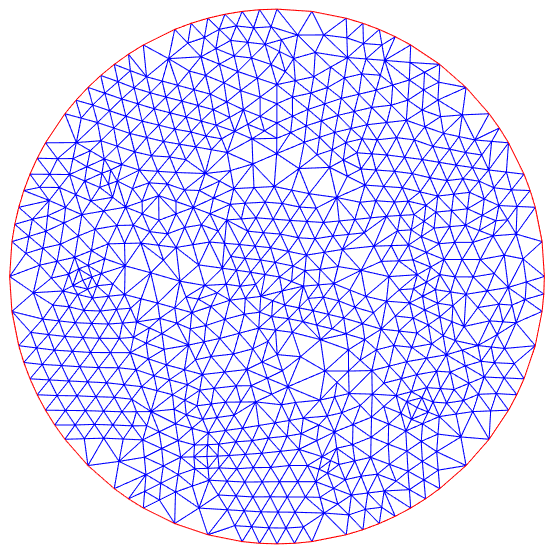}
		%&\includegraphics[width=1in]{CircleMin13TH.png}
		&\includegraphics[width=1in]{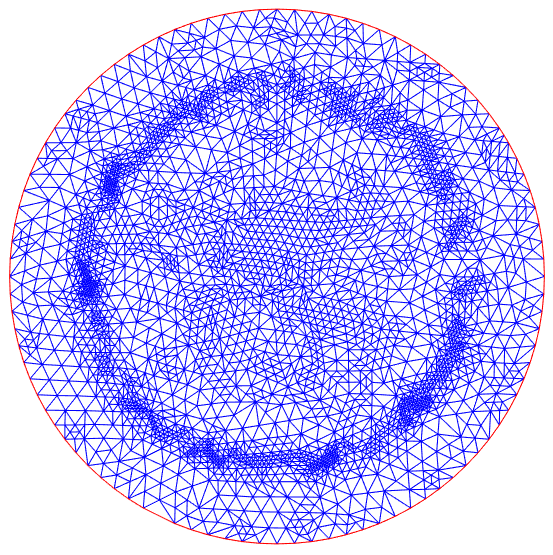}
		%&\includegraphics[width=1in]{CircleMin15TH.png} 
        &\includegraphics[width=1in]{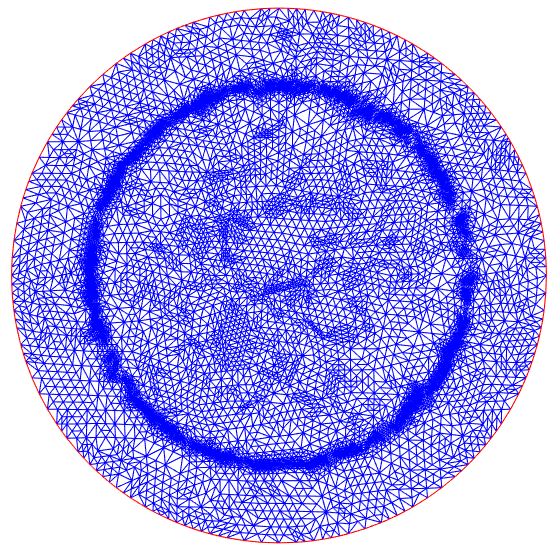}\\		%\includegraphics[width=1in]{CircleMin1Initialphi.png}
        \includegraphics[width=1in]{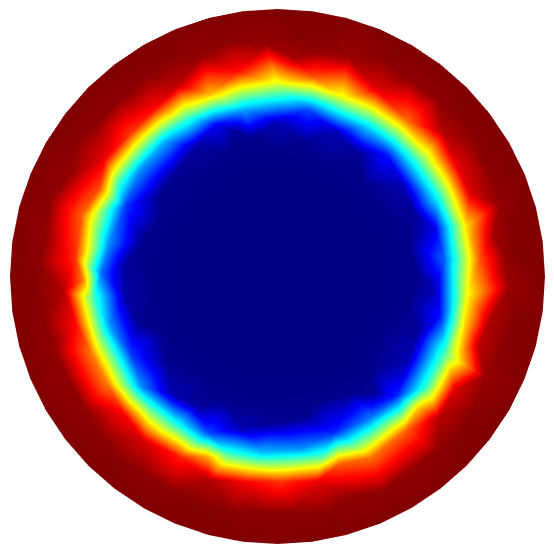}               
        &\includegraphics[width=1in]{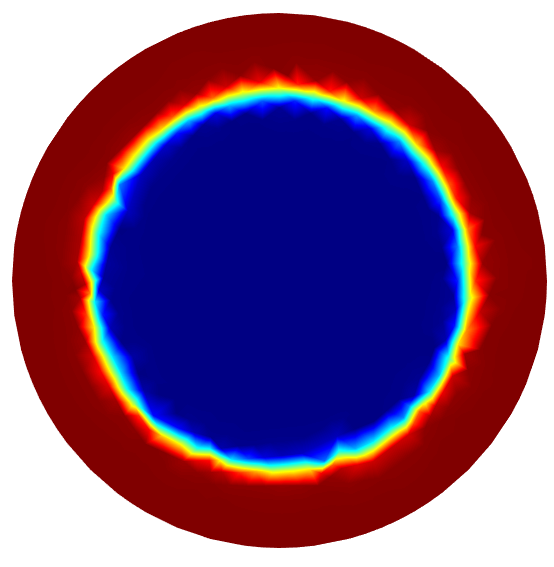}
        %&\includegraphics[width=1in]{CircleMin1FinalDesign3TH.png}  
        &\includegraphics[width=1in]{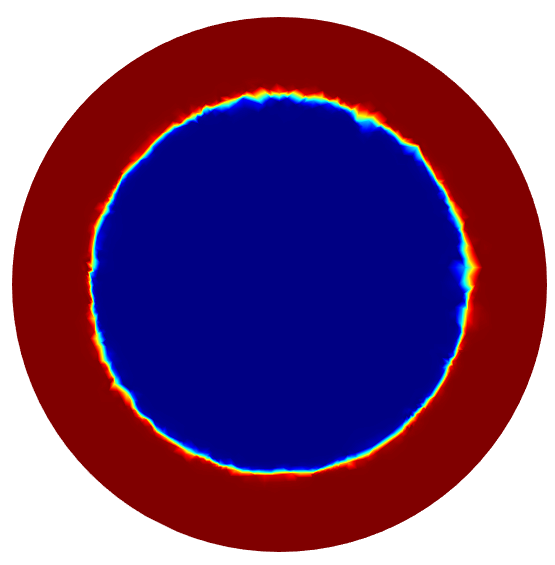}
    	%&\includegraphics[width=1in]{CircleMin1FinalDesign5TH.png}
        &\includegraphics[width=1in]{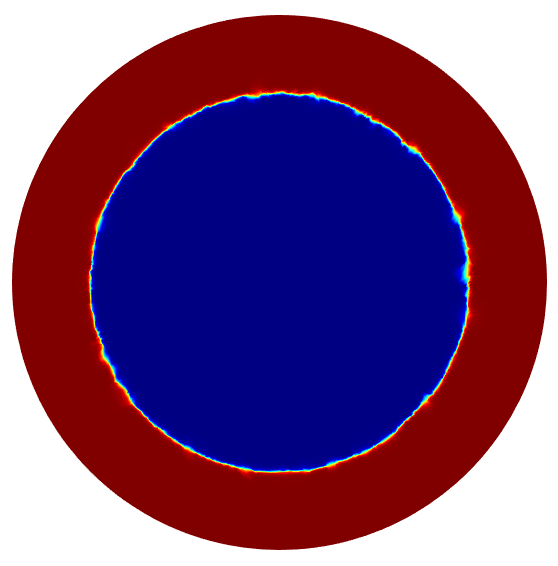}&
        \\
        \includegraphics[width=1in]{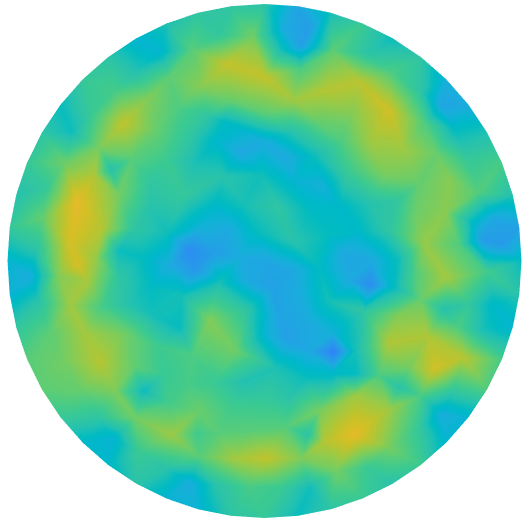}
        &\includegraphics[width=1in]{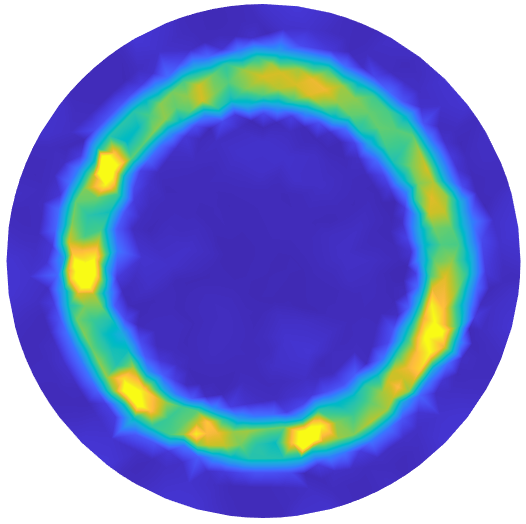}
        %&\includegraphics[width=1in]{CircleMin13th1EstimatorS.png}
        &\includegraphics[width=1in]{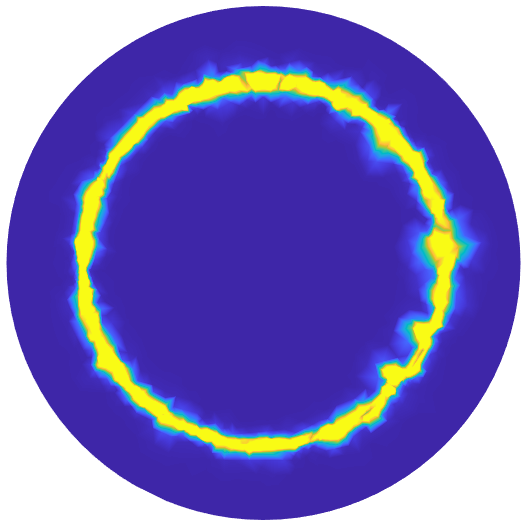}
        %&\includegraphics[width=1in]{CircleMin15th1EstimatorS.png}
        &\includegraphics[width=1in]{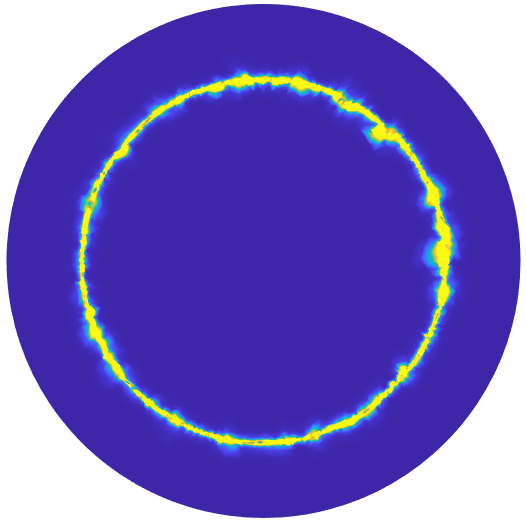}
        &\includegraphics[width=0.2in]{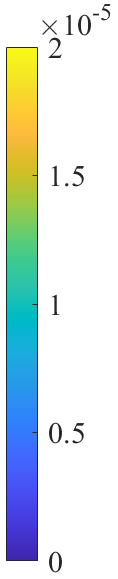}\\
        \includegraphics[width=1in]{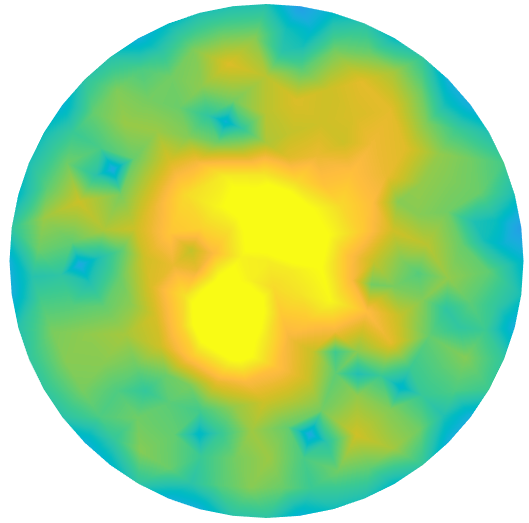}
        &\includegraphics[width=1in]{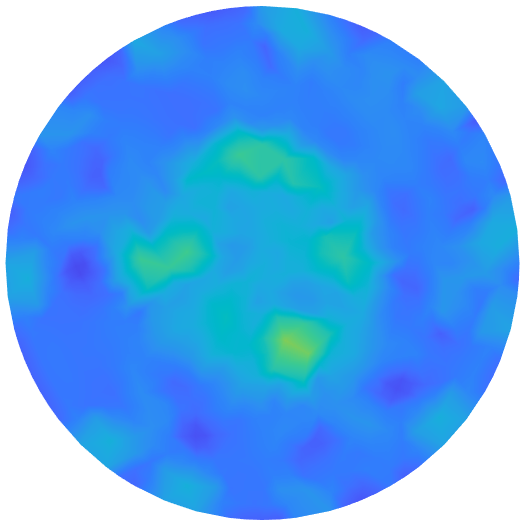}
        %&\includegraphics[width=1in]{CircleMin13th2EstimatorS.png}
        &\includegraphics[width=1in]{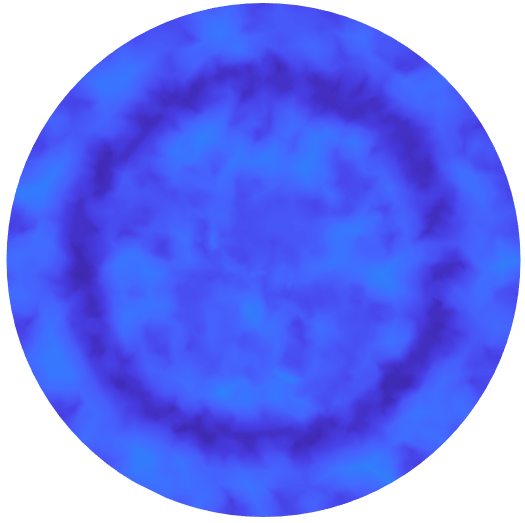}
        %&\includegraphics[width=1in]{CircleMin15th2EstimatorS.png}
        &\includegraphics[width=1in]{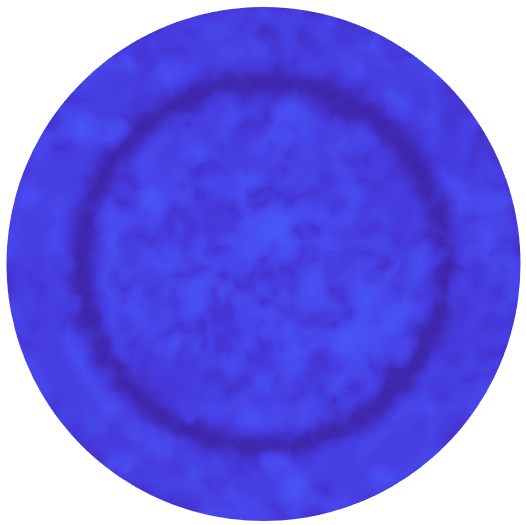}
        &\includegraphics[width=0.17in]{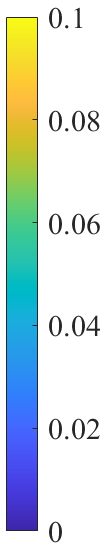}
	\end{tabular}
\caption{Evolution of adaptive mesh levels $k = 0,1,3,5$ for Example \ref{example2}, with the number of vertices on each mesh being 269, 817, 3095 and 10102. The 2nd, 3rd and 4th rows show optimized designs $\phi_k^\ast$ with a numerical optimal shape highlighted in red over each $\cT_k$ and two error indicators $\eta_{k,0}$, $\eta_{k,1}$ respectively.
}\label{AdaptiveExp2Min1}
\end{figure}

\begin{figure}[htbp]
\centering
\subfigure{\includegraphics[width=0.3\textwidth]{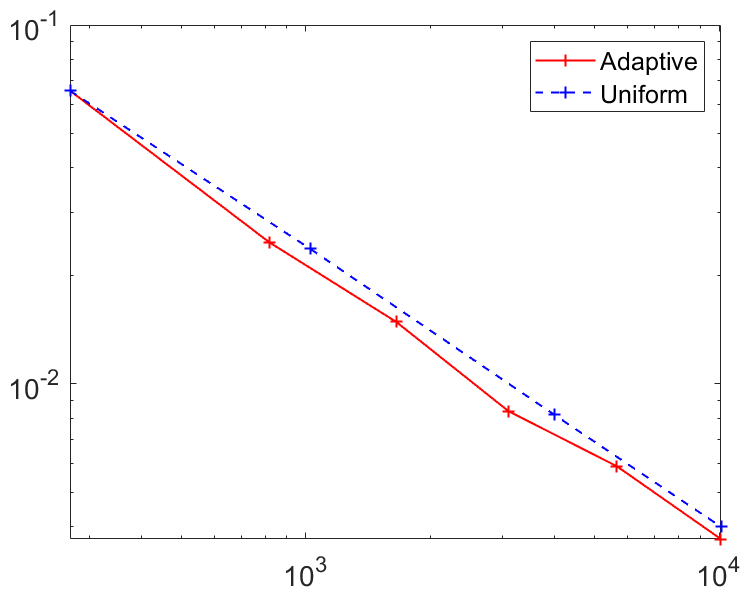}}
\subfigure{\includegraphics[width=0.3\textwidth]{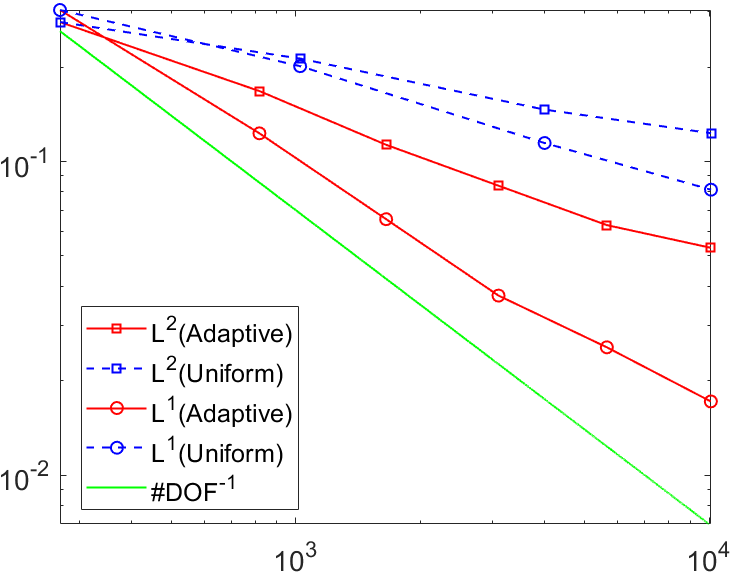}}
\subfigure{\includegraphics[width=0.29\textwidth]{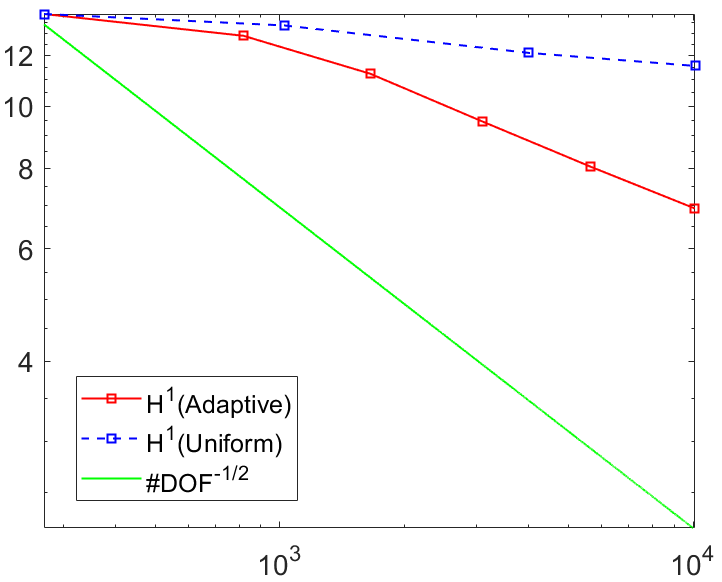}}
\caption{Errors in minimizing $\lambda_{1}$ (left) and approximating the associated phase-field function (middle and right) as functions of degrees of freedom for the adaptive refinement and the uniform refinement in Example \ref{example2}. $|\lambda_{k,1}^{\eps,\ast}-\lambda_{1}^{\mathrm{ref}}|$ (left), $L^1(D)$ and $L^2(D)$ errors (middle), and $H^1(D)$ error (right) of numerical phase-field functions.}\label{ObjExp2Min1Error}
\end{figure}

\begin{figure}[htbp]
	\centering\setlength{\tabcolsep}{0pt}
	\begin{tabular}{ccccc}
		\includegraphics[width=1in]{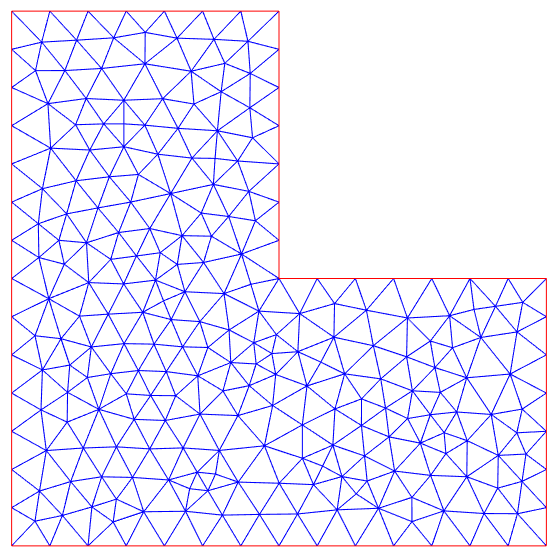}
		&\includegraphics[width=1in]{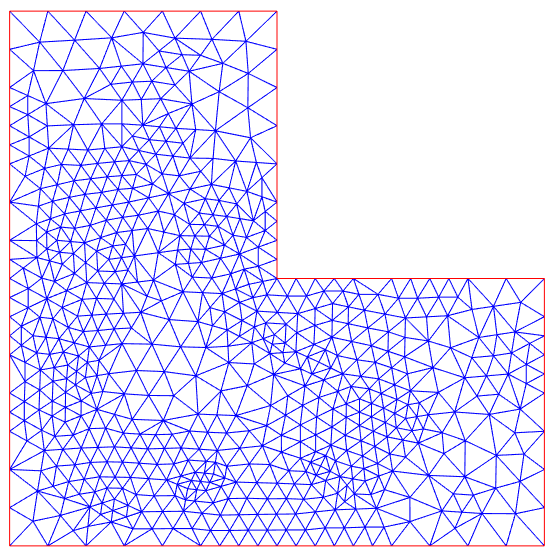}
		%&\includegraphics[width=1in]{LMin13TH.png}
		&\includegraphics[width=1in]{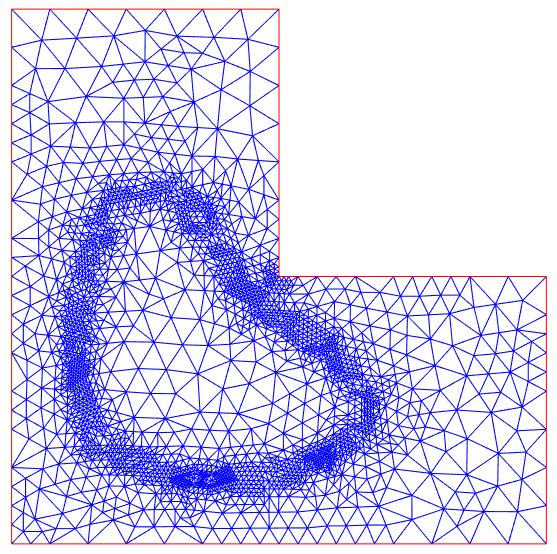}
		%&\includegraphics[width=1in]{LMin15TH.png}
		&\includegraphics[width=1in]{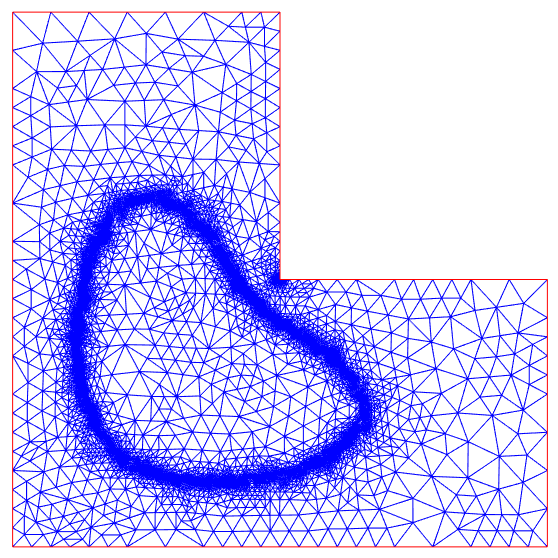}&\\
		\includegraphics[width=1in]{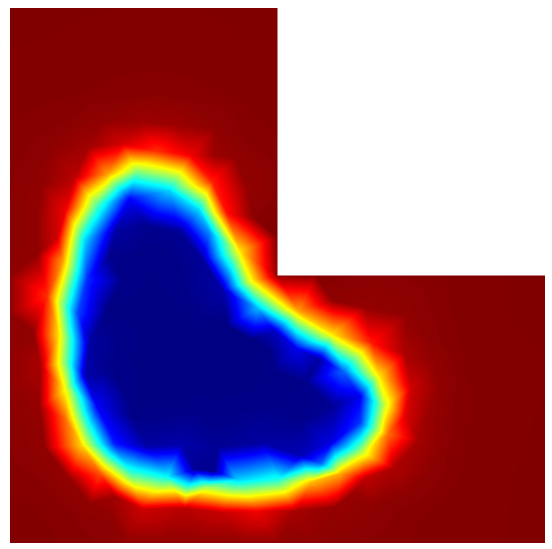}  &\includegraphics[width=1in]{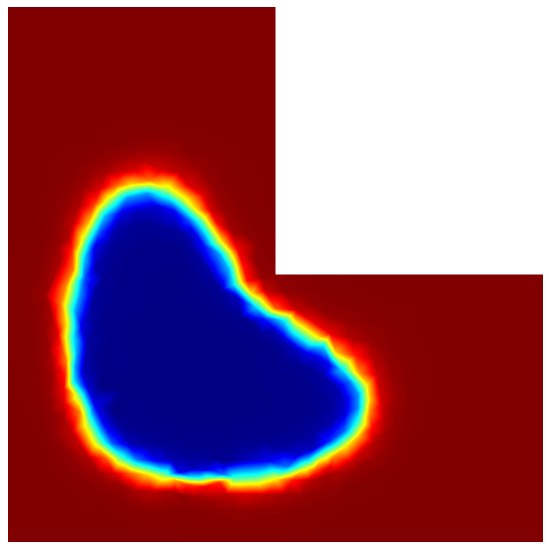}
		%&\includegraphics[width=1in]{LMin1FinalDesign3TH.png} 
        &\includegraphics[width=1in]{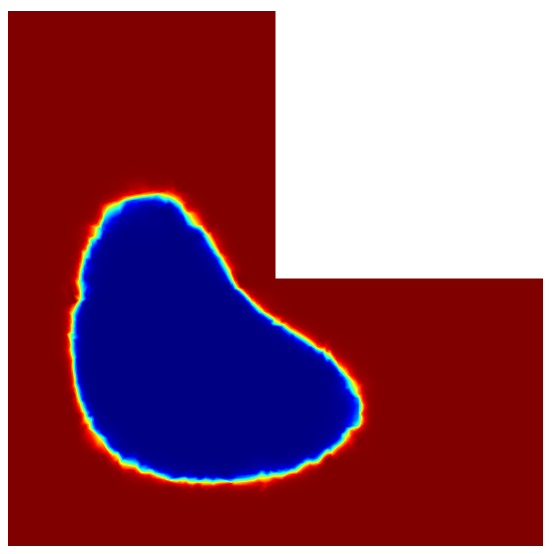}
		%&\includegraphics[width=1in]{LMin1FinalDesign5TH.png}
        &\includegraphics[width=1in]{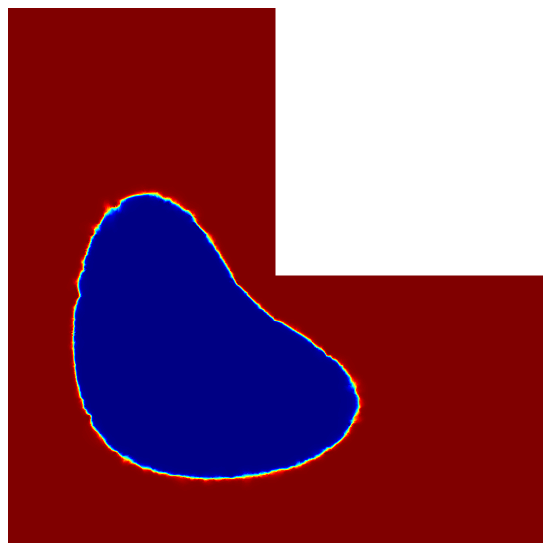}&\\
        \includegraphics[width=1in]{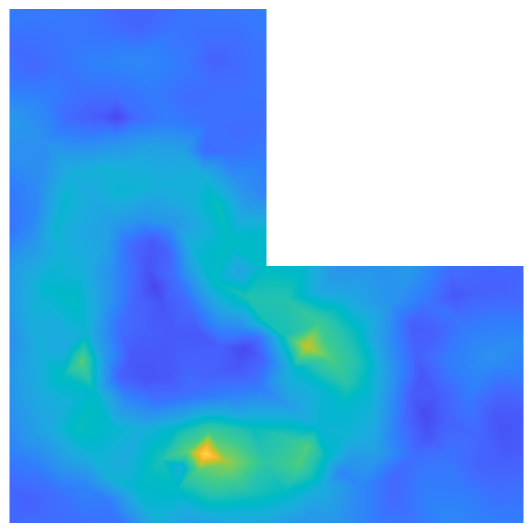}
        &\includegraphics[width=1in]{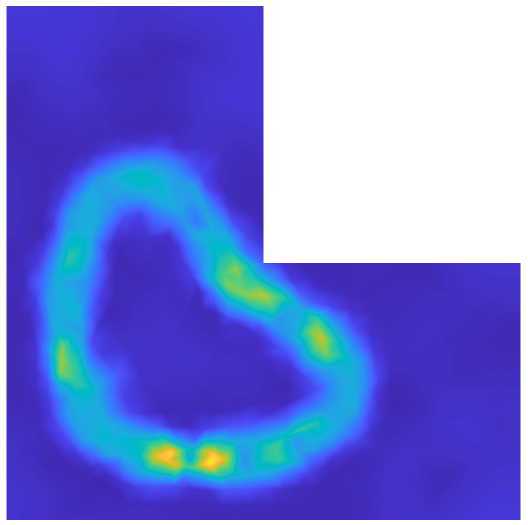}
        %&\includegraphics[width=1in]{LMin13th1EstimatorS.png}
        &\includegraphics[width=1in]{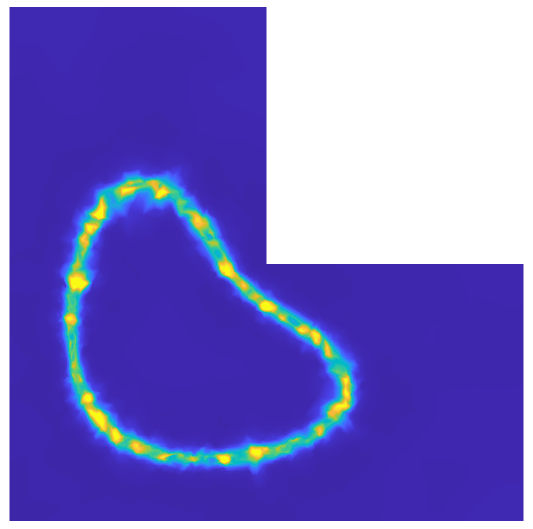}
        %&\includegraphics[width=1in]{LMin15th1EstimatorS.png}
        &\includegraphics[width=1in]{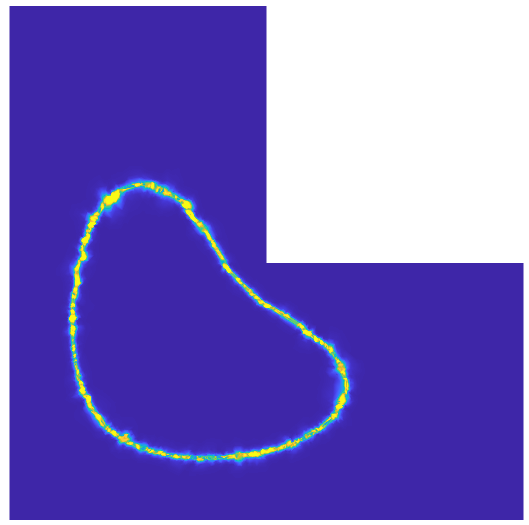}
        &\includegraphics[width=0.2in]{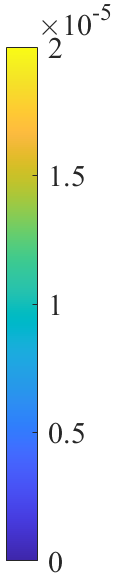}\\
        \includegraphics[width=1in]{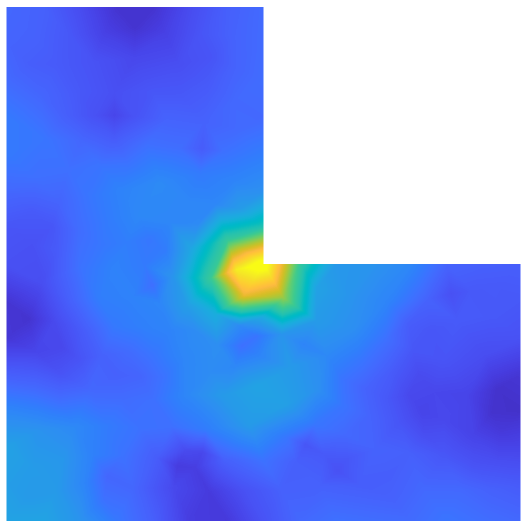}
        &\includegraphics[width=1in]{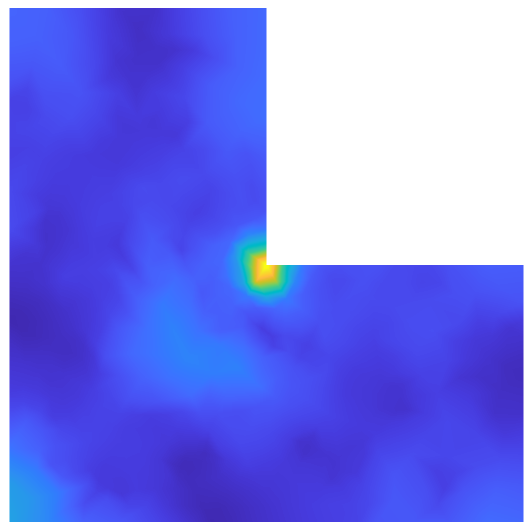}
        %&\includegraphics[width=1in]{LMin13th2EstimatorS.png}
        &\includegraphics[width=1in]{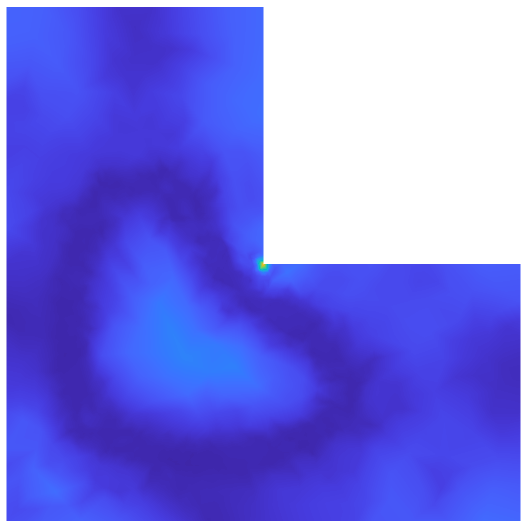}
        %&\includegraphics[width=1in]{LMin15th2EstimatorS.png}
        &\includegraphics[width=1in]{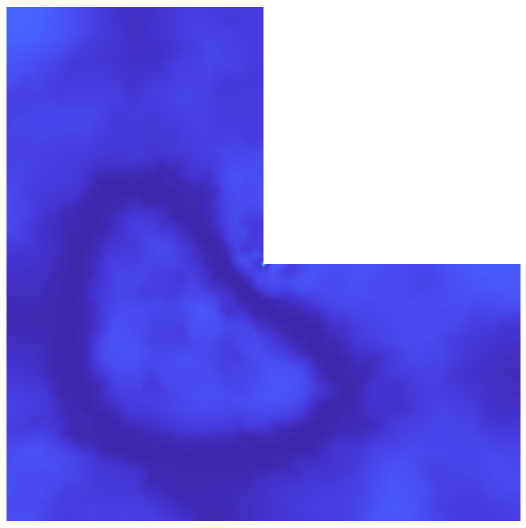}
        &\includegraphics[width=0.17in]{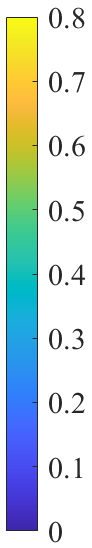}
	\end{tabular}
\caption{Evolution of adaptive mesh levels $k = 0,1,3,5$ for Example \ref{example4}, with the number of vertices on each mesh being 259, 663, 2462 and 9083. The 2nd, 3rd and 4th rows show optimized designs $\phi_k^\ast$ with a numerical optimal shape highlighted in red over each $\cT_k$ and two error indicators $\eta_{k,0}$, $\eta_{k,1}$ respectively.}\label{AdaptiveExp4Min1}
\end{figure}

\begin{figure}[htbp]
\centering
\subfigure{\includegraphics[width=0.32\textwidth]{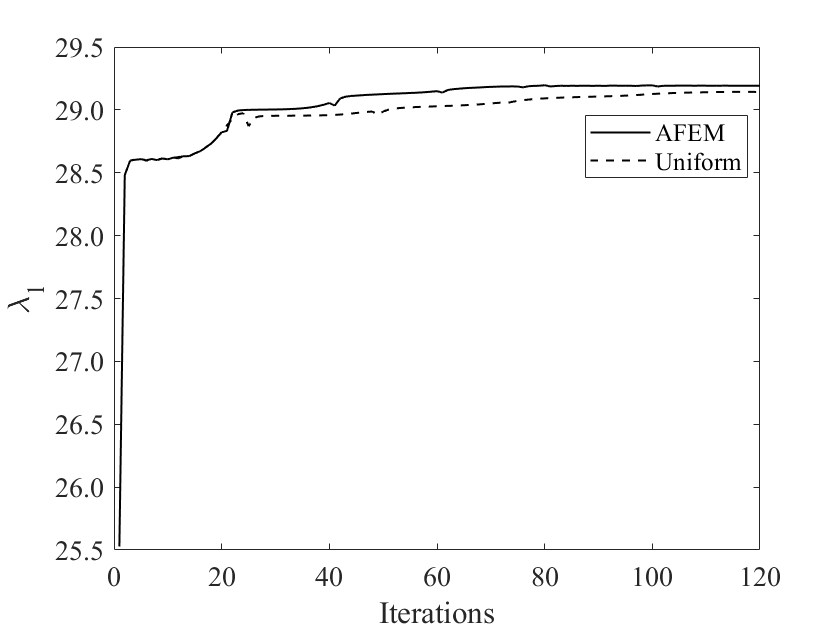}}
\subfigure{\includegraphics[width=0.32\textwidth]{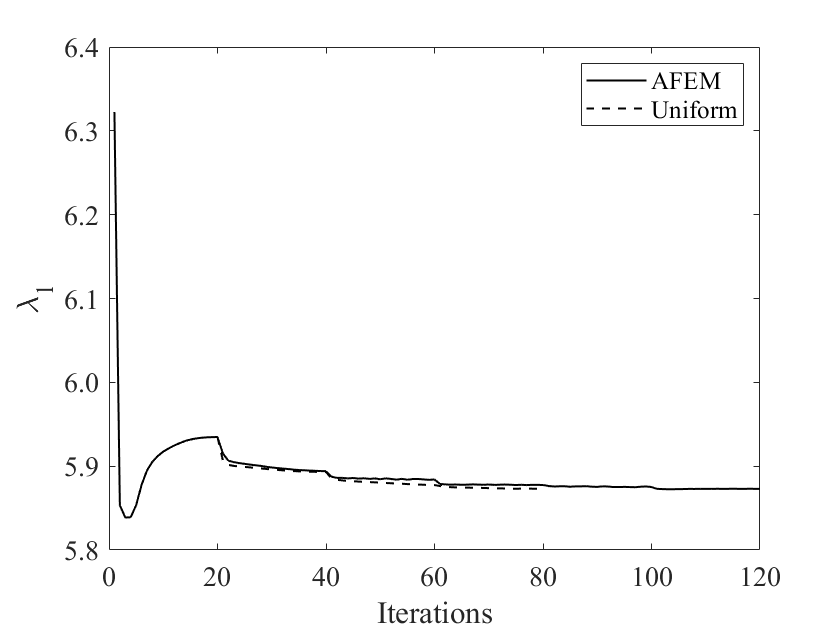}}
\subfigure{\includegraphics[width=0.32\textwidth]{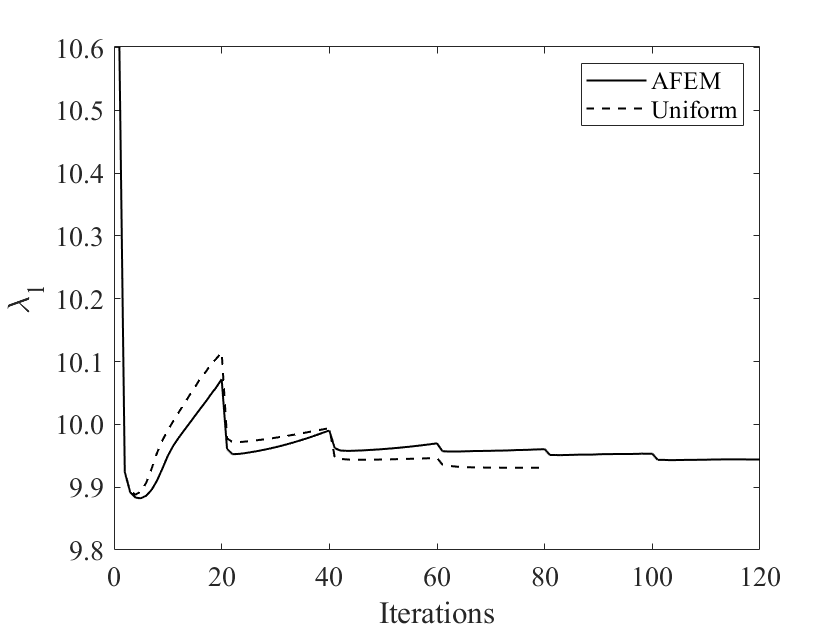}}
\caption{Convergence history of objectives as functions of the total number ($K \times N$) of outer-iterations performed in Algorithm \ref{Alga} for Example \ref{example1} (left), Example \ref{example2} (middle) and Example \ref{example4} (right). Example \ref{example1}: $K=6$, $N=20$ (adaptive), and $K=5$, $N=24$ (uniform). Example \ref{example2}: $K=6$, $N=20$ (adaptive), and $K=4$, $N=20$ (uniform). Example \ref{example4}: $K=6$, $N=20$ (adaptive), and $K=4$, $N=20$ (uniform).
% $K-1=5$ refinements, $N=20$ iterations per mesh and $K-1=3$ refinements, $N=20$ iterations per mesh by adaptive and uniform strategies respectively.
}\label{ObjExp1to3}
\end{figure}

\begin{figure}[htbp]
\centering
\subfigure{\includegraphics[width=0.32\textwidth]{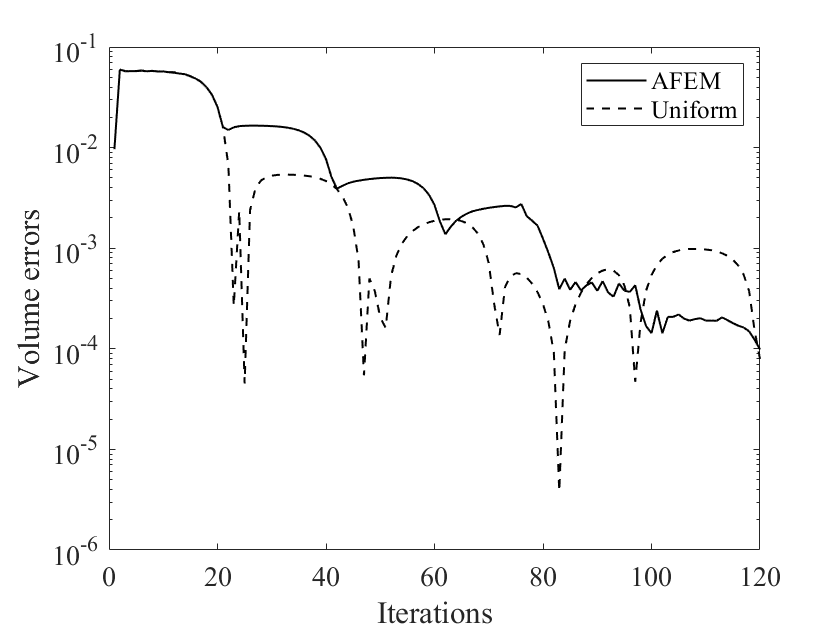}}
\subfigure{\includegraphics[width=0.32\textwidth]{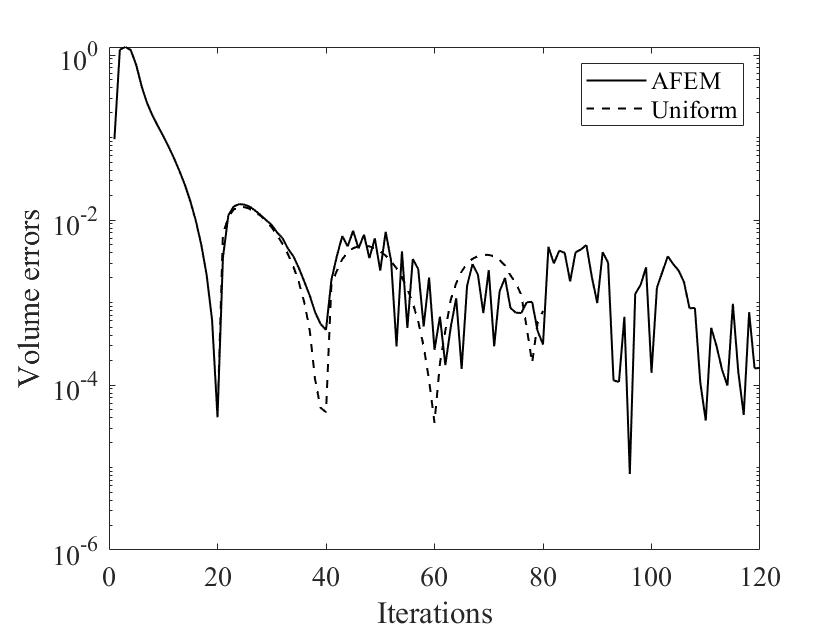}}
\subfigure{\includegraphics[width=0.32\textwidth]{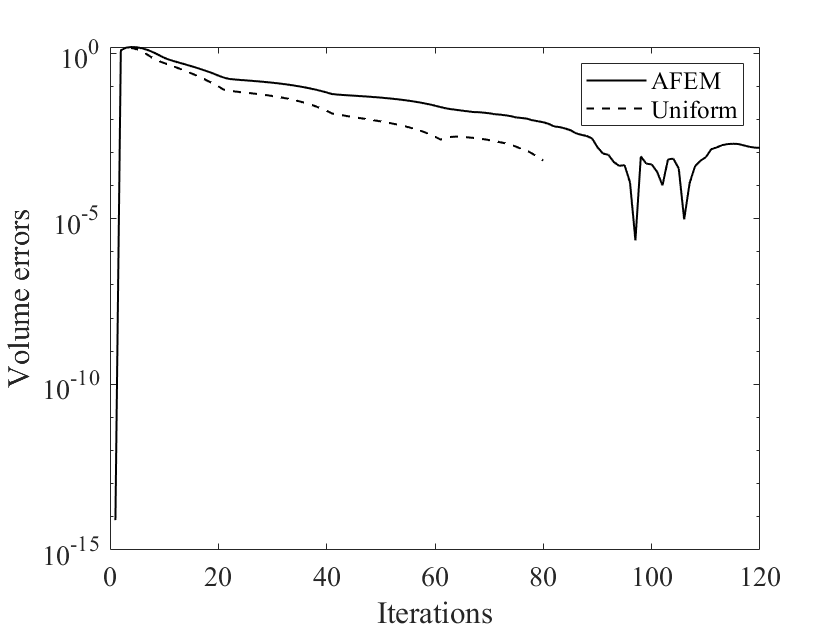}}
\caption{Errors in the volume constraint as functions of the total number ($K \times N$) of outer-iterations performed in Algorithm \ref{Alga} for Example \ref{example1} (left), Example \ref{example2} (middle) and Example \ref{example4} (right).
% $K-1=5$ refinements, $N=20$ iterations per mesh and $K-1=3$ refinements, $N=20$ iterations per mesh by adaptive and uniform strategies respectively.
}\label{Volerrors1to3}
\end{figure}

\begin{figure}[htbp]
\centering
\subfigure{\includegraphics[width=0.34\textwidth]{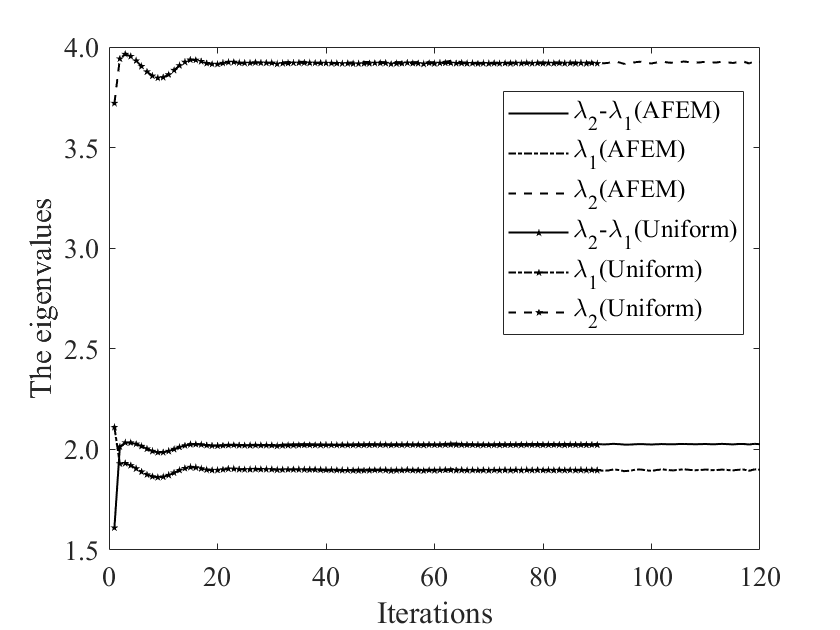}}
\subfigure{\includegraphics[width=0.34\textwidth]{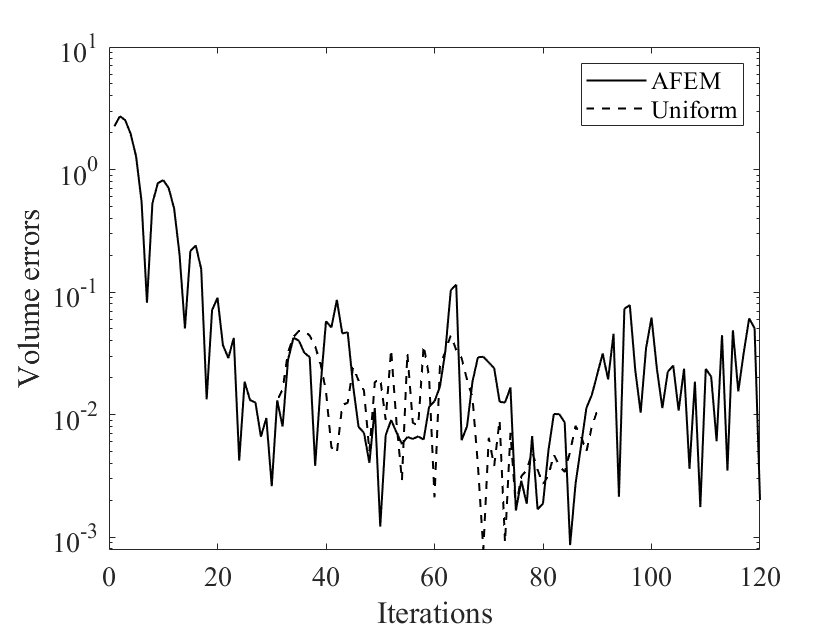}}
\caption{Convergence history of objectives (left) and volume constraint errors (right) as functions of the total number ($K\times N$) of outer-iterations performed in Algorithm \ref{Alga} for Example \ref{example7}. $K=4$ mesh levels, $N=30$ iterations per mesh and $K = 3$ mesh levels, $N=30$ iterations per mesh by adaptive and uniform strategies respectively.}
	\label{ObjExp7Max21}
\end{figure}

\begin{figure}[htbp]
	\centering
	\begin{tabular}{c}
        \includegraphics[width=1.7in]{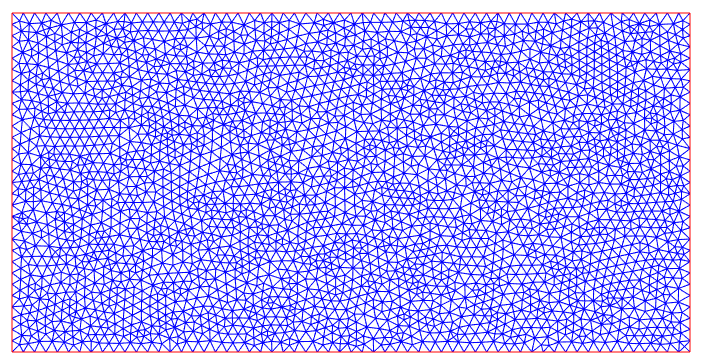}
        \includegraphics[width=1.7in]{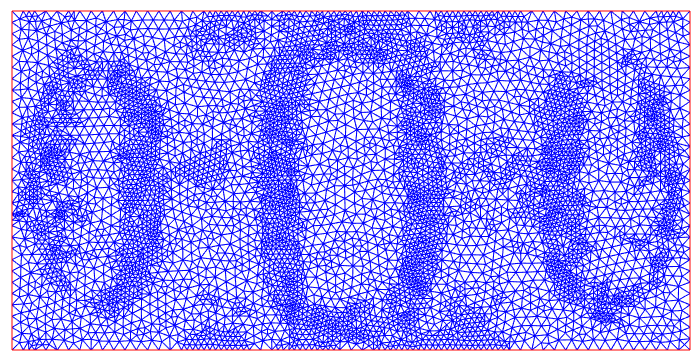}
		\includegraphics[width=1.7in]{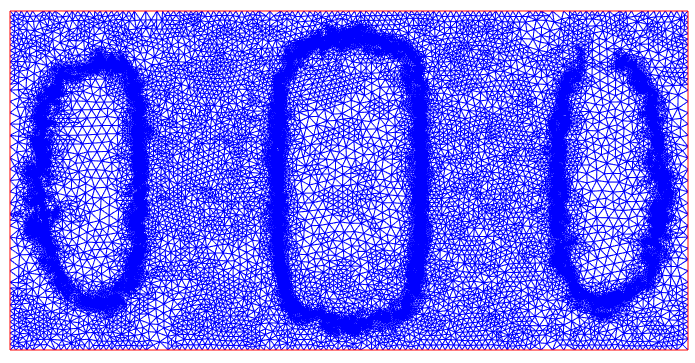}\\
        \includegraphics[width=1.7in]{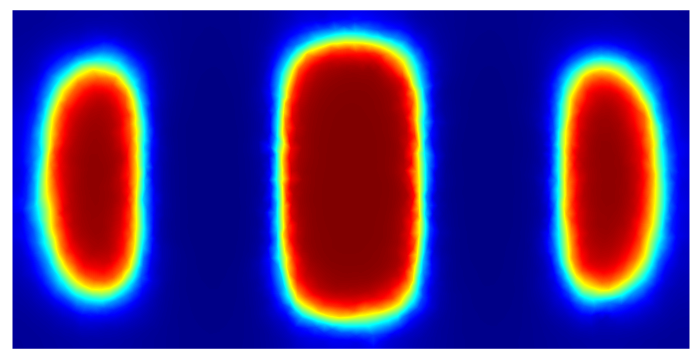}
        \includegraphics[width=1.7in]{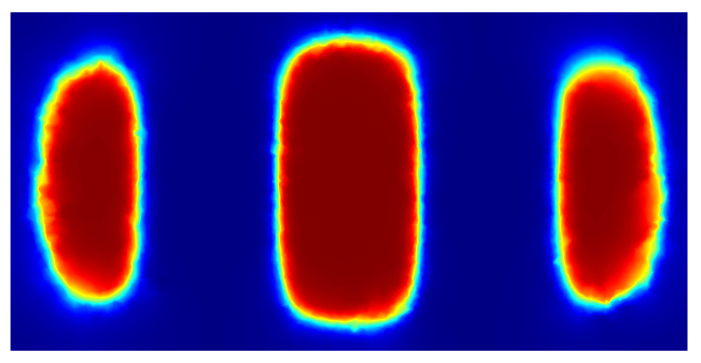}
        \includegraphics[width=1.7in]{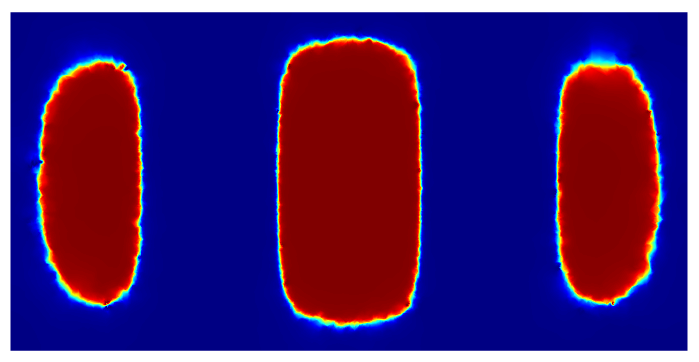}
  %       \\
		% \includegraphics[width=1.8in]{RectMax32FinalDesign1THColorBar.png}
		% \includegraphics[width=1.8in]{RectMax32FinalDesign2THColorBar.png}
		% \includegraphics[width=1.8in]{RectMax32FinalDesign3THColorBar.png}
		% \includegraphics[width=1.8in]{RectMax32FinalDesign4THColorBar.png}
	\end{tabular}
\caption{Evolution of adaptive mesh levels $k = 0,1,3$ for Example \ref{example3}, with the number of vertices on each mesh being 3118, 6254 and 23837. The 2nd row shows optimized designs $\phi_k^\ast$ with a numerical optimal shape highlighted in red over each $\cT_k$.
%The evolution of the adaptive mesh level from $k=0$ (initial) to 3 (top left to bottom right in the 1st and 2nd rows) to maximize $\lambda_{3}-\lambda_{2}$ in Example \ref{example3}(b), with the number of vertices of each mesh being 3118, 6254, 12282 and 23837. The 3rd and 4th rows (top left to bottom right) show optimized designs over the associated meshes.
}\label{AdaptiveExp3Max32}
\end{figure}

\begin{figure}[htbp]
\centering
\subfigure{\includegraphics[width=0.34\textwidth]{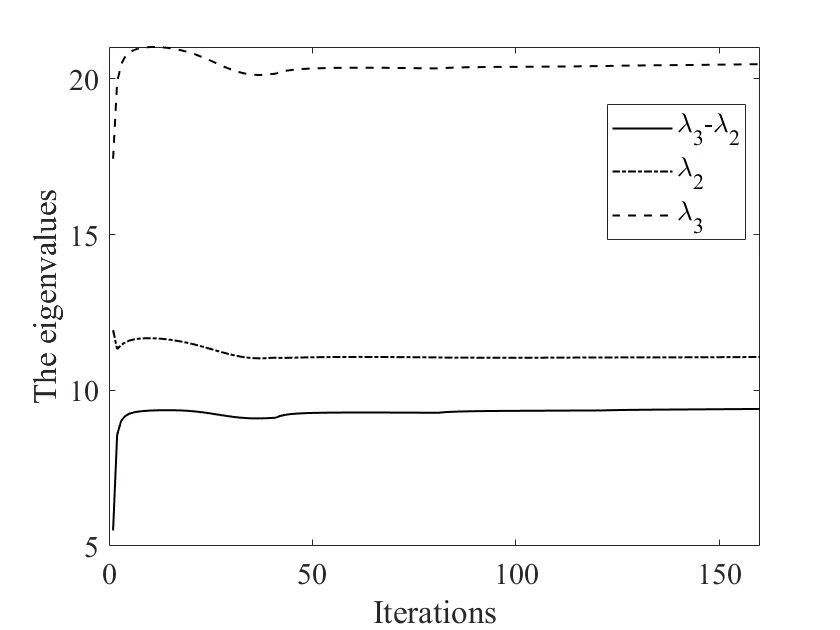}}
\subfigure{\includegraphics[width=0.34\textwidth]{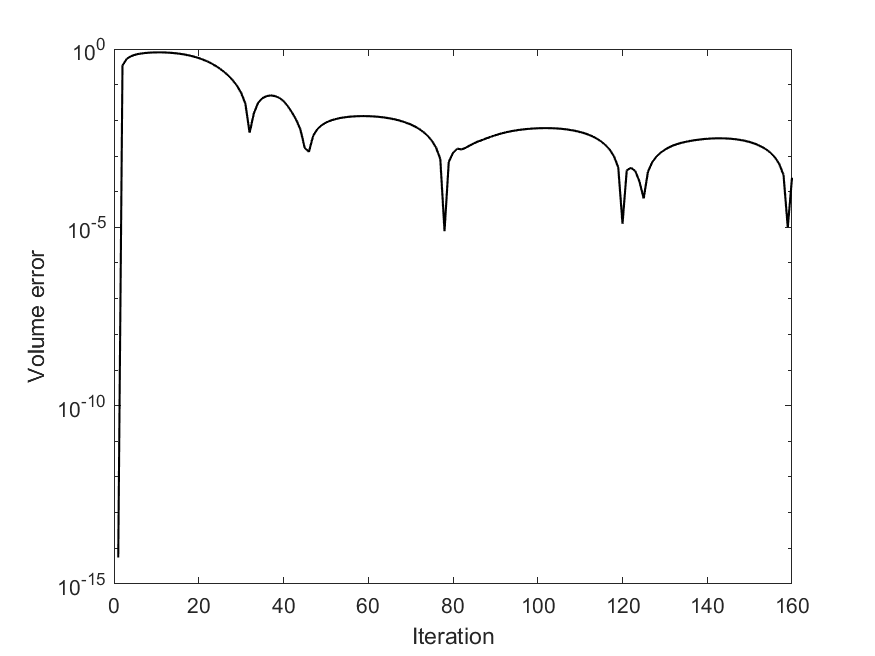}}
\caption{Convergence history of the objective (left) and the volume constraint error (right) as functions of the total number ($K \times N$) of  outer-iterations performed in Algorithm \ref{Alga} for Example \ref{example3}. $K=4$ mesh levels and $N=40$ iterations per mesh in the adaptive process.}\label{ObjExp3Max32}
\end{figure}

\begin{figure}[htbp]
	\centering
	\begin{tabular}{cc} 
        % \includegraphics[width=1in]{SquareMin1OptimizedDesign6TH.png}
        % \includegraphics[width=1in]{SquareMin1UFinalDesign5TH.png}
        % \hspace{2mm}
        \includegraphics[width=1in]{SquareMax1FinalDesign6TH.png}
	  \includegraphics[width=1in]{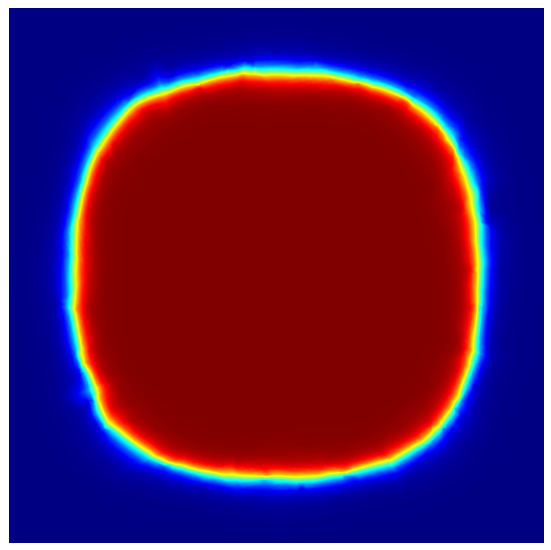}
        &%\hspace{2mm}
        \includegraphics[width=1in]{CircleMin1FinalDesign6TH.png}
        \includegraphics[width=1in]{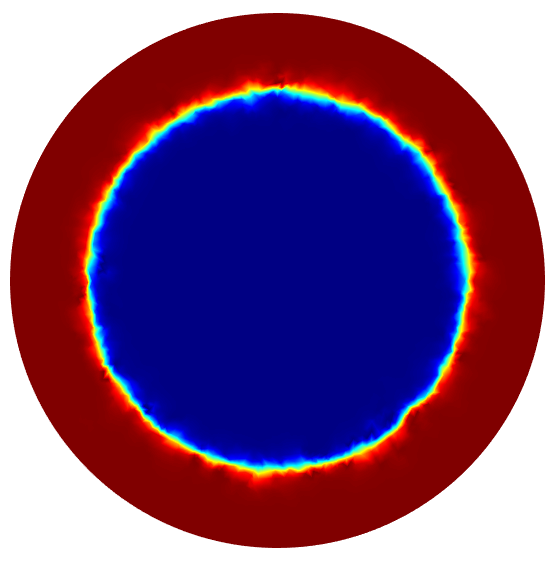}
\\
        % (a)  $\min\lambda_{1}$ in Example \ref{example1}(a)
        % \hspace{8mm}
         $\max\lambda_{1}$ in Example \ref{example1} & %\hspace{12mm}
         $\min\lambda_{1}$ in Example \ref{example2}
\\
        \includegraphics[width=1in]{LMin1FinalDesign6TH.png}
        \includegraphics[width=1in]{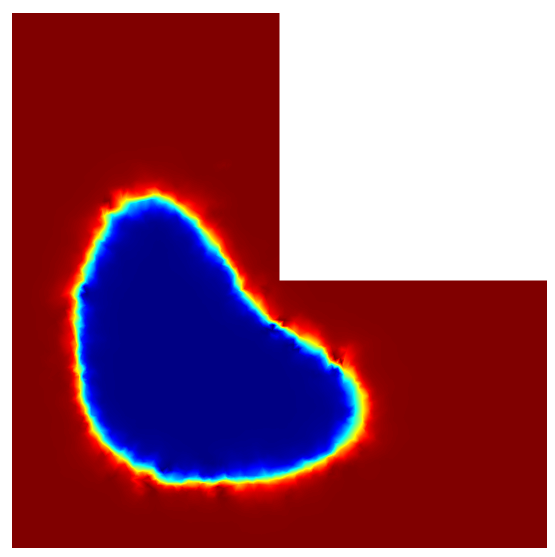}
         &%\hspace{2mm}
         \includegraphics[width=1in]{IrrPMax21FinalDesign4TH.png}
        \includegraphics[width=1in]{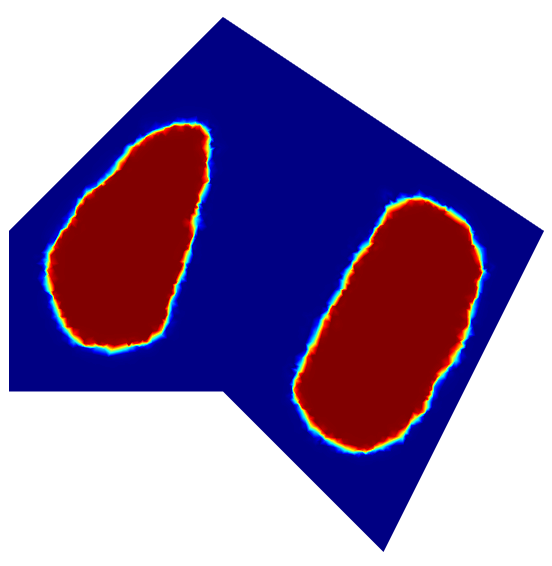}
  \\
         $\min\lambda_{1}$ in Example \ref{example4}
        &%\hspace{5mm}
        % (e) $\max\lambda_{1}$ in Example \ref{example4}(b)
        % \includegraphics[width=1in]{IrrPMin1FinalDesign4TH.png}		\includegraphics[width=1in]{IrrPMin1UFinalDesign3TH.png}
        %  \hspace{2.5mm}
        % (f) $\min\lambda_{1}$ in Example \ref{example7}(a)
        % \hspace{5mm}
         $\max\lambda_{2}-\lambda_{1}$ in Example \ref{example7}
	\end{tabular}
\caption{Comparison of optimal designs $\phi_{K-1}^\ast$ by adaptive (left within each pair) and uniform (right within each pair) refinements with numerical optimal shapes highlighted in red for Examples \ref{example1}-\ref{example7}. %Within each pair, the left design and the right design are produced by adaptive and uniform refinements respectively.
}\label{fig:designCompare}
\end{figure}

In Figures \ref{AdaptiveExp1Max1}-\ref{fig:designCompare}, we present the results for Examples \ref{example1}-\ref{example3}, including the evolution of meshes, numerical optimal designs and error estimators during the adaptive refinement process. We also compare errors of phase-field functions, objective function values and volume constraints by adaptive and uniform refinement strategies. It is clear that as the adaptive algorithm progresses, more local refinements are performed to identify the diffuse interfaces in the interior and capture singularities around corners. In particular, when optimizing a single eigenvalue, the discrete minimizers $\phi_0^\ast$ over the initial meshes ($k=0$) provide rough sketches of the desired designs and then adaptive mesh refinements yield more fine detail of the diffuse interfaces as the accuracy of discrete eigenfunctions improves. We note that this observation holds true as well in maximizing the band gap between two consecutive eigenvalues as illustrated in Figure \ref{AdaptiveExp7Max21} (Example \ref{example7}) and Figure \ref{AdaptiveExp3Max32} (Example \ref{example3}). %The true design in Example \ref{example2} is a circle $\{x|x_{1}^2+x_{2}=R^2\}$, where the value of $R>0$ is chosen to satisfy the volume constraint with $C=0.5$. 
Since the optimal shape exists for Example \ref{example2}, Figure \ref{ObjExp2Min1Error} presents the absolute errors of the numerical first eigenvalues as well as the errors of the numerical phase-field functions in $L^1(D)$, $L^2(D)$ and $H^1(D)$ norms for both adaptive and uniform refinements, plotted as functions of the degrees of freedom. The reference phase-field function is an approximation of the Heaviside function defined on a high-resolution mesh with more than 50,000 degrees of freedom and given by $\psi = \sqrt{x_{1}^2+x_{2}^2}-1/\sqrt{2}$ %Denote by $H(\cdot)$ the Heaviside function. The true phase field function for reference is set as a smoothing $H(\psi)$  
via \cite{QianHuZhu:2022}
\[
H(\psi) \approx\left\{\begin{array}{cl}
\frac{1}{2}\left(1+\frac{\psi}{\tau_{1}}+\frac{1}{\pi} \sin \left(\frac{\pi \psi}{\tau_{1}}\right)\right) & \text { if }|\psi| \leq \tau_{1}, \\
1 & \text { if } \psi>\tau_{1}, \\
0 & \text { if } \psi<-\tau_{1},
\end{array}\right.
\]
where $\tau_{1}=0.02$ is prescribed as twice the mesh size. The corresponding first eigenvalue of \eqref{eigen_vp_phasefield_disc} using the aforementioned $H(\psi)$ over the same high-resolution mesh, denoted by $\lambda_1^{\mathrm{ref}}$, is considered as the reference minimum. %In order to compare with the true solution, the eigenvalue problem \eqref{eigen_vp_phasefield_disc} is computed over the same fine mesh is taken as the reference optimum. The error curves between adaptive solution and true solution of $\lambda_{1}$ and optimal phase field function for Example \ref{example2} are presented in Figure \ref{ObjExp2Min1Error}. 
The convergence history in Figure \ref{ObjExp2Min1Error} demonstrates that the AFEM achieves more accurate approximations of the reference phase-field function compared to the uniform refinement method, despite both methods utilizing the same number of degrees of freedom. These findings underscore the effectiveness and efficiency of the proposed adaptive algorithm.

Another noteworthy observation is the decrease of each estimator $\eta_{k,i}$ ($0\leq i\leq 2$) towards zero, as displayed in Figures \ref{AdaptiveExp1Max1}, \ref{AdaptiveExp2Min1}, \ref{AdaptiveExp4Min1} and \ref{AdaptiveExp7Max21} for Examples \ref{example1}-\ref{example7}. This trend corresponds to the conclusion of vanishing limits in Theorem \ref{thm:estimator->zero}. Furthermore, $\eta_{k,0}$ detects jumps in the discrete phase-field functions at internal interfaces, while $\eta_{k,1}$ and $\eta_{k,2}$ capture singularities in the associated eigenfunction around re-entrant corners.

Figures \ref{ObjExp1to3}, \ref{Volerrors1to3}, \ref{ObjExp7Max21} and \ref{ObjExp3Max32} compare the convergence of associated objective function values and volume constraint errors, defined as $|G(\phi)|$ at the beginning of Section 3.1, by the adaptive refinement and the uniform refinement in Examples \ref{example1}-\ref{example7}. Here $K$ and $N$ denote total mesh levels and the number of iterations per mesh in Algorithm \ref{Alga} respectively. The levels of uniform mesh refinements are appropriately selected such that the numbers of vertices on the final meshes are comparable to those by adaptive refinement processes. Notably, no significant difference in performance of the two refinement strategies is observed in terms of the objective functions and the volume constraint in $\U$ is well preserved during the computations by both strategies. We also note that the objective function value initially decreases or increases very rapidly and then fluctuates steadily. This observation is also evident in maximizing the band gap $\lambda_3-\lambda_2$, as noted in Figure \ref{ObjExp3Max32} for Example \ref{example3}. In addition, the final objective function values by both refinement strategies are nearly identical in Examples \ref{example1}-\ref{example7}, consistent with the similar optimized designs depicted in Figure \ref{fig:designCompare}.

Table \ref{tab:compare_results} presents more quantitative results for both refinement strategies, including the value of each objective function in Examples \ref{example1}-\ref{example7}
(achieved at the final mesh generated by adaptive and uniform refinements), the number of vertices on the final mesh and the computing time (in second). As previously mentioned, the optimal objective function values obtained through adaptive mesh refinements are very close to those by uniform mesh refinements. However, compared to the uniform refinement, the adaptive algorithm saves $12\%$-$43\%$ of computing time to attain a similar level of objective function values. Figures \ref{ObjExp1to3} (middle and right) and \ref{ObjExp7Max21} (left) in Examples \ref{example2}-\ref{example7} indicate that, despite requiring fewer iterations in the augmented Lagrangian method, the standard uniform mesh refinement incurs more computing time than the adaptive mesh refinement. Additionally, Figure \ref{fig:designCompare} displays the optimized designs generated by the two strategies in Examples \ref{example1}-\ref{example7}, revealing that the diffuse interfaces produced by the adaptive algorithm are smoother and exhibit fewer small oscillations compared to those generated by uniform mesh refinements. This highlights the advantage of the adaptive strategy in better identifying the diffuse interface and the singularity of the associated eigenfunction.

% \begin{figure}
% 	\centering
% 	\begin{tabular}{c}
% 		\includegraphics[width=1in]{CircleMin1FinalDesign6TH.png}
% 		\includegraphics[width=1in]{CircleMin1UFinalDesign4TH.png}
% 	\end{tabular}
% 	\caption{Optimized designs by adaptive (left) and uniform (right) refinement for Example \ref{example2}.}\label{AdaptiveExp2Min_FinalDesigns}
% \end{figure}

%%%%remove
% \begin{figure}[htbp]
% \centering
% \subfigure{\includegraphics[width=0.34\textwidth]{CircleMin1eigsV.png}}
% \subfigure{\includegraphics[width=0.34\textwidth]{CircleMin1VolErrV.png}}
% \caption{Convergence history of the objective function value (left) and the error of the volume constraint (right) as functions of the total number ($K \times N$) of outer-iterations performed in Algorithm \ref{Alga} to minimize $\lambda_{1}$ in Example \ref{example2}. $K-1=5$ refinements, $N=20$ iterations per mesh and $K-1=3$ refinements, $N=20$ iterations per mesh by adaptive and uniform strategies respectively.}\label{ObjExp2Min1}
% \end{figure}

% \begin{table}[htbp]
% 	\centering
% 	\caption{Comparison between adaptive and uniform refinements in Example \ref{example2}.}
% 	\begin{tabular}{|c|cc|c|}
% 		\hline
% 		\multirow{2}{*}{Method} & \multicolumn{2}{c|}{Final mesh}                 & \multirow{2}{*}{CPU-time(s)} \\ \cline{2-3}
% 		& \multicolumn{1}{c|}{Vertices} & Obiective &                              \\ \hline
% 		AFEM                    & \multicolumn{1}{c|}{10102}    & 5.8729          & 2041.61                      \\ \hline
% 		Uniform                 & \multicolumn{1}{c|}{10111}    & 5.8732          & 2796.11                      \\ \hline
% 	\end{tabular}\label{TimeExp1Min1}
% \end{table}
%%%%%%%%%%%%%%%%%%%%%%%%%%%%%%%%%%%%%%%%%%%%%%%%%%%%%%%%%%%%%%%%%%%%%%%%%%%%%
%Set $\beta=100$ and $C=0.75$ for minimizing $\lambda_{1}$ on the L-shaped domain. The initial phase function is 0.75. Fig. \ref{AdaptiveExp4Min1} shows that the adaptive algorithm converges to an optimized design after five adaptive mesh refinements. The refinement of the mesh is not only on the interface, but also at the corner points of the L-shaped area.
%Fig. \ref{ObjExp4Min1} demonstrates the both convergence histories and volume error. The uniform result is shown in Fig. \ref{fig:designCompare}(d), and it can be seen that the diffusion interface of adaptive refinement is clearer than the uniform. From Table \ref{tab:compare_results}, the objective value of uniform refinement is identical with the adaptive one. However, the CPU time of adaptive algorithm is reduced by 42.74\%.

%%%remove
% \begin{figure}[htbp]
% \centering
% \subfigure{\includegraphics[width=0.34\textwidth]{LMin1eigsV.png}}
% \subfigure{\includegraphics[width=0.34\textwidth]{LMin1VolErrV.png}}
% \caption{Convergence history of the objective function value (left) and the error of the volume constraint (right) as functions of the total number ($K \times N$) of outer-iterations performed in Algorithm \ref{Alga} to minimize $\lambda_{1}$ in Example \ref{example4}(a). $K-1=5$ refinements, $N=20$ iterations per mesh and $K-1=3$ refinements, $N=20$ iterations per mesh by adaptive and uniform strategies respectively.}\label{ObjExp4Min1}
% \end{figure}

\begin{table}[hbtp]
\centering
\begin{threeparttable}
\caption{Comparison of quantitative results by adaptive and uniform refinement strategies for Examples \ref{example1}-\ref{example7}: the number of vertices on the final mesh $\cT_{K-1}$, the associated objective function value and the total computing time. Saving refers to the saving in the
overall computing time. \label{tab:compare_results}}
\begin{tabular}{|cl|ccc|ccc|c|}
    \hline
    \multicolumn{2}{|c|}{\multirow{2}{*}{Example}}
    & \multicolumn{3}{c}{adaptive}&\multicolumn{3}{|c|}{uniform} & \multirow{2}{*}{saving}\\
    \cmidrule(lr){3-5} \cmidrule(lr){6-8}&
     & vertices & objective & time (sec) & vertices & objective & time (sec) & \\
    \hline
    % \ref{example1}(a) & min $\lambda_{1}$ & 12909    & 20.4782   & 1819.80 & 12917    & 20.5363   & 1964.09 & 7.35\% \\
    \ref{example1} & max $\lambda_{1}$& 13331    & 29.1916   & 1944.73 & 13347    & 29.1420   & 2230.13 & 12.80\% \\ \hline
     \ref{example2} & min $\lambda_{1}$ & 10102    & 5.8729    & 2041.61 & 10111    & 5.8732    & 2796.11 & 26.98\%\\\hline
    \ref{example4} & min $\lambda_{1}$ & 9083     & 9.9310    & 1174.42 & 9096     & 9.9310    & 2051.03 & 42.74\%\\
    % \ref{example4}(b) & max $\lambda_{1}$ & 7543     & 10.3026   & 754.92  & 7563     & 10.2839   & 1022.06 & 26.14\% \\
    \hline
    % \ref{example7}(a) & min $\lambda_{1}$ & 5661     & 1.7557    & 592.45  & 5646     & 1.7598    &  700.33  & 15.40\%\\
     \ref{example7} & max $\lambda_{2}-\lambda_{1}$ & 7777     & 2.0262    & 1285.73 & 7740     & 2.0222    & 1791.91 & 28.25\%\\
    \bottomrule
    \end{tabular}
\end{threeparttable}
\end{table}

%%%%%%%%%%%%%%%%%%%%%%%%%%%%%%%%%%%%%%%%%%%%%%%%%%%%%%%%%%%%%%%%%%%%%%%%%%%%%%%%%%%%%%%%%%%%%%

\begin{figure}[htbp]
    \centering\centering\setlength{\tabcolsep}{0pt}
		\begin{tabular}{c}
			\includegraphics[width=1.65in]{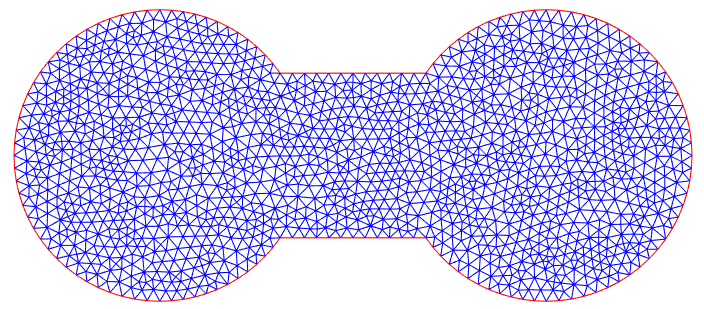}
			\includegraphics[width=1.65in]{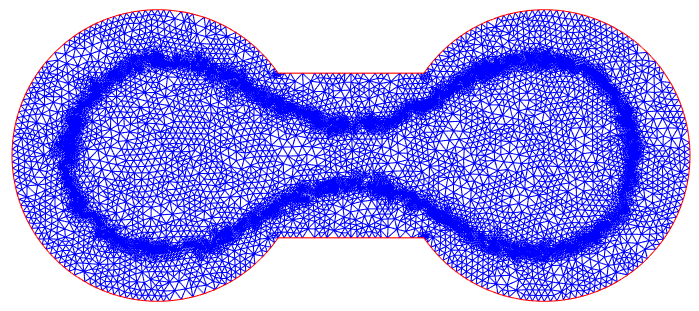}		\includegraphics[width=1.65in]{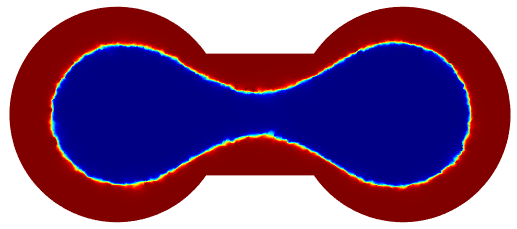}
   \\
            \includegraphics[width=1.65in]{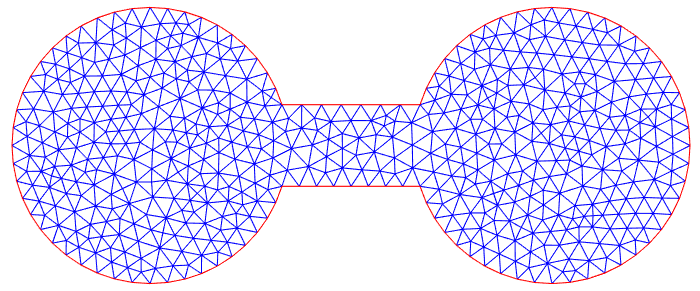}
            \includegraphics[width=1.65in]{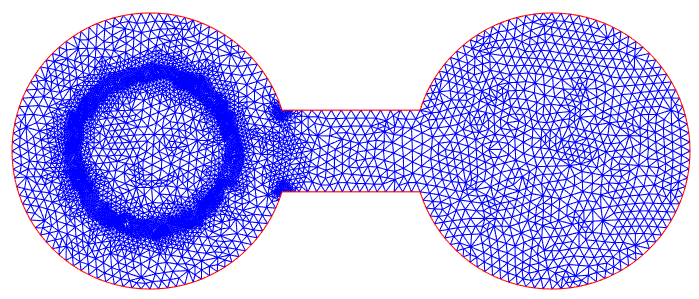}
	      \includegraphics[width=1.65in]{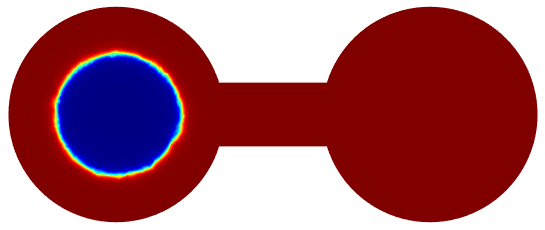}
	\end{tabular}
\caption{Adaptive mesh levels $k=0$, $k=3$ and the numerical optimal shape highlighted in red over $\cT_3$ to minimize $\lambda_{1}$ in Example \ref{example5}(a) (top) and Example \ref{example5}(b) (bottom) with the number of vertices on each mesh being 1546, 11184 and 578, 5062 respectively.}\label{AdaptiveExp5Min1}
\end{figure}

\begin{figure}[htbp]
    \centering\centering\setlength{\tabcolsep}{0pt}
    \begin{tabular}{cccccc}
  %       \includegraphics[width=0.88in]{AnnulusMax1a1Initialmesh.png}
		% % &\includegraphics[width=0.88in]{AnnulusMax1a12TH.png}
		% % &\includegraphics[width=0.88in]{AnnulusMax1a13TH.png}
		% &\includegraphics[width=0.88in]{AnnulusMax1a14TH.png}
  %       &\includegraphics[width=0.88in]{AnnulusMax1a1FinalDesign4TH.png}\\
		\includegraphics[width=0.88in]{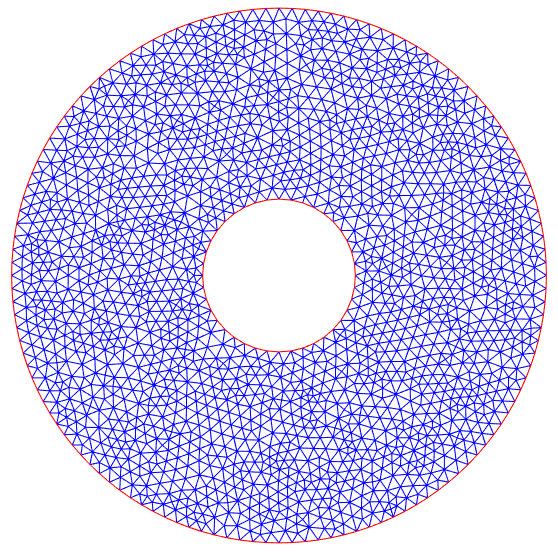}
		% &\includegraphics[width=0.88in]{AnnulusMin1a102TH.png}
		% &\includegraphics[width=0.88in]{AnnulusMin1a103TH.png}
		&\includegraphics[width=0.88in]{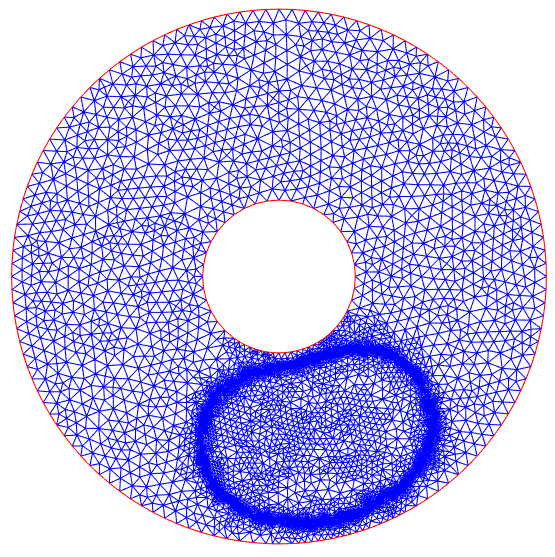}
        &\includegraphics[width=0.88in]{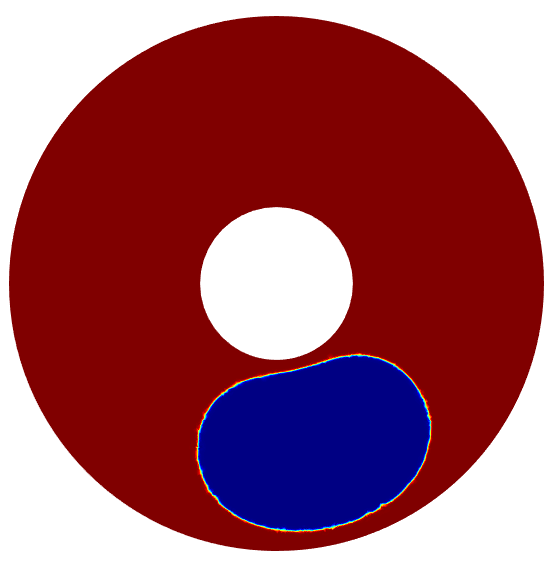}
        &\includegraphics[width=0.88in]{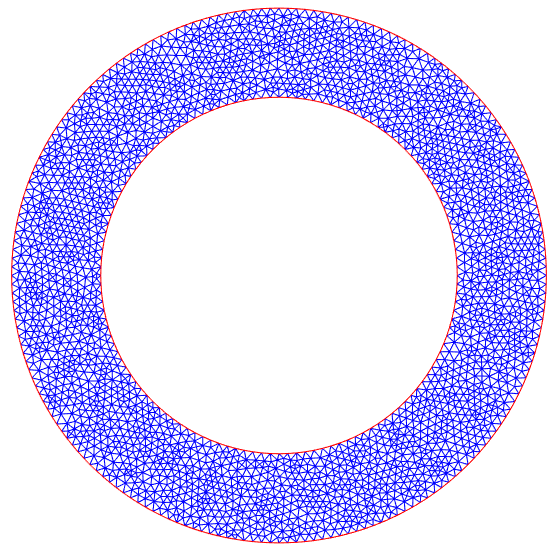}
        &\includegraphics[width=0.88in]{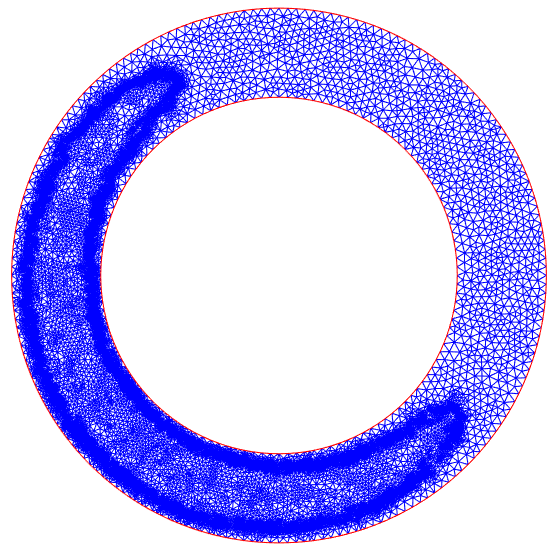}
        &\includegraphics[width=0.88in]{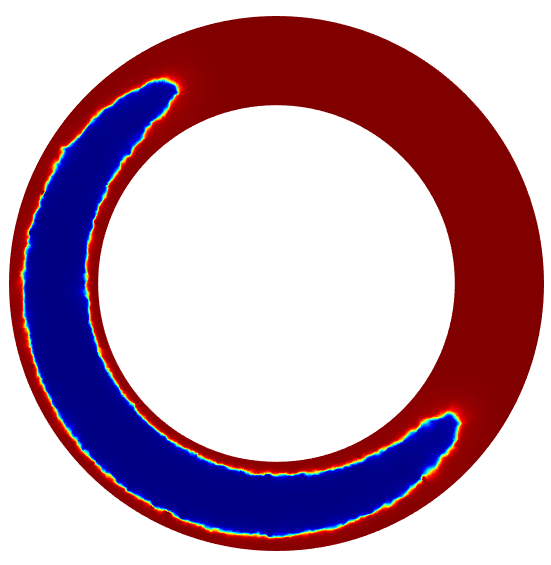}\\
        (a) & (b) & (c) & (d) & (e) & (f)        
    \end{tabular}
\caption{Adaptive mesh levels $k=0$, $k=3$ and the numerical optimal shape highlighted in red over $\cT_3$ to maximize $\lambda_1$ ((a)-(c)) in Example \ref{example6}(a) with the number of vertices being 2779 and 15116, and to minimize $\lambda_1$ ((d)-(f)) in Example \ref{example6}(b) with the number of vertices being 1929 and 6143.} %Meshes and final design by adaptive refinements in Example \ref{example6} ($\eta=3.5$). The first row ($\alpha=1$ for maximizing $\lambda_{1}$): the number of vertices of each mesh is 2779, 7001, 9946 and 15116, the optimized design on the final mesh. The second row ($\alpha=10$ for minimizing $\lambda_{1}$): the number of vertices of each mesh is 1929, 2627, 3816 and 6143, the optimized design on the final mesh.}
 \label{AdaptiveExp6Eta3.5}
\end{figure}

In Examples \ref{example5}-\ref{example6}, we investigate the intriguing symmetry preservation/breaking phenomenon in a dumbbell and an annulus, respectively \cite{Chanillo:2000}. As stated earlier, the evolution of meshes generated by the adaptive algorithm in Figures \ref{AdaptiveExp5Min1} and \ref{AdaptiveExp6Eta3.5} demonstrates that the refinements primarily concentrate in regions where the diffuse interfaces are located.

Figure \ref{AdaptiveExp5Min1} (top) illustrates that in a symmetric dumbbell $D_\mathcal{H}$ with $\mathcal{H}=0.6$ and $C=0.5$, the optimal shape (highlighted in red) maintains symmetry with respect to the $x_2$-axis. In contrast, this symmetry is vanished in the complement (highlighted in blue) of the optimal shape in Figure \ref{AdaptiveExp5Min1} (bottom), where the dumbbell has a handle width of $\mathcal{H}=0.3$ and a larger prescribed volume $\Omega$ ($C=0.85$). This corresponds to the assertion in Theorem 7 of \cite{Chanillo:2000} that symmetry breaking occurs when the handle width $2\mathcal{H}$ of $D_\mathcal{H}$ is below a critical value and the prescribed volume exceeds a certain threshold. Moreover, as noted in \cite{Chanillo:2000}, reflecting the blue region across the $x_2$-axis indicates a lack of uniqueness for the optimal shape in this case.

In Figure \ref{AdaptiveExp6Eta3.5}, we observe a lack of rotational symmetry in the numerical optimal shape (highlighted in red) of $D$ for Example \ref{example6}(a), with $\alpha=10$ and $C=0.85$ in \eqref{eigen-opt_phasefield}-\eqref{eigen_vp_phasefield}. Additionally, this absence of symmetry is also evident in Example \ref{example6}(b), as asserted in Theorem 6 of \cite{Chanillo:2000}. These findings emphasize the deviation from expected symmetry in optimal shapes under the specified conditions.

%For $\mathcal{H}=0.6$, set $\beta=80$, $C=0.5$, $\zeta=0.03$ and $\gamma=10^{-3}$. The initial phase function is uniform with 0.8. Let $\theta_{1}=0.7$ and $\theta_{2}=0.1$.
%Fig. \ref{AdaptiveExp5Min1} shows adaptive meshes and optimized designs with symmetry preserving property. %The convergence histories of objective function and volume error are displayed in Fig.\ref{ObjExp5Min1}.
%For $\mathcal{H}=0.3$, set $\beta=150$, $C=0.85$, $\zeta=0.03$ and $\gamma=5\times10^{-3}$. The initial phase function is 0.5. Let $\theta_{1}=0.8$ and $\theta_{2}=0.3$.
%Fig. \ref{AdaptiveExp5Min2} shows the adaptive meshes and symmetry-breaking optimized designs when the volume of $\Omega$ and the handle of dumbbell reduce. As the computation progresses, the refined mesh catches the interface precisely and becomes more and more apparent.

%%%remove
% \begin{figure}[htbp]
%     \centering\centering\setlength{\tabcolsep}{0pt}
%     \begin{tabular}{ccccc}
% 		\includegraphics[width=0.88in]{AnnulusMin1a11TH.png}
% 		&\includegraphics[width=0.88in]{AnnulusMin1a12TH.png}
% 		&\includegraphics[width=0.88in]{AnnulusMin1a13TH.png}
% 		&\includegraphics[width=0.88in]{AnnulusMin1a14TH.png}
% &\includegraphics[width=0.88in]{AnnulusMin1a1FinalDesign4TH.png}
%     \end{tabular}
%     \caption{Evolution of adaptive mesh levels from $k=0$ (initial) to 3 and the optimized design $\phi_3^\ast$ with red corresponding to a numerical optimal shape over $\cT_3$ to minimize $\lambda_1$ in Example \ref{example6}(c), with the number of vertices of each mesh being 2775, 5381, 10141 and 19088 respectively.} %Meshes and final design for minimizing $\lambda_{1}$ ($\alpha=1$ and $\eta=1.5$) by adaptive refinements in Example \ref{example6}. The number of vertices of each mesh is 2775, 5381, 10141 and 19088.}
%     \label{AdaptiveExp6Min1eta1.5}
% \end{figure}
		
		%In this section, we discuss the convergence of Algorithm %\ref{alg_afem_eigenvalue} without any stopping criterion. The %analysis is divided into three parts. We explain the reason %to use the estimators in the module ESTIMATE in subsection %\ref{subsect:estimator}. In subsection \ref{subsect:limit}, %the limiting behavior of Algorithm \ref{alg_afem_eigenvalue} %characterized by an auxiliary limiting problem is carefully %examined. Finally, the convergence result of Algorithm %\ref{alg_afem_eigenvalue} is presented in subsection %\ref{subsect:conv}.
		
\section{Origin of error estimator}\label{sec:estimator}
In this section, we explain the reason to use the estimators in the module \textsf{ESTIMATE}. The motivation of our error estimators $\eta_{k,0}(\phi_k^\ast)$ and $\eta_{k,j}(\phi_{k}^{\ast},\lambda_{k,i_j}^{\eps,\ast},w_{k,i_j}^{\ast})$ $(1\leq j\leq l)$ comes from the following residuals associated with \eqref{optsys_phase-field}
\begin{equation}\label{residual_0}
    \langle \mathcal{R}_0(\phi_k^\ast), \phi \rangle : = \gamma\mathcal{F}_\eps'(\phi^\ast_k)\phi + \displaystyle\alpha\sum_{j=1}^l\frac{\p \Psi}{\p \lambda_{i_j}^{\eps}}(\lambda^{\eps,\ast}_{k,i_1},\lambda^{\eps,\ast}_{k,i_2},\cdots,\lambda^{\eps,\ast}_{k,i_l})\int_{\Omega}(w_{k,i_j}^{\ast})^2 \phi   \dx   \quad \forall\phi \in \U,
\end{equation}
\begin{equation}\label{residual_j}
    \langle \mathcal{R}_j(\phi_{k}^{\ast},\lambda_{k,i_j}^{\eps,\ast},w_{k,i_j}^{\ast}), v \rangle : = \int_D \bold{\nabla} w_{k,i_j}^{\ast}  \cdot \bold{\nabla} v \dx + \alpha\int_D \phi_k^\ast w_{k,i_j}^{\ast} v \dx - \lambda_{k,i_j}^{\eps,\ast} \int_D \phi_k^\ast w_{k,i_j}^{\ast} v \dx \quad \forall v \in H^1_0(D).
\end{equation}
Our argument relies on a quasi-interpolation operator
$\Pi_{k}: L^1(D)\to V_k^1$ \cite{JinXu:2020}, where $V_k^1$ is the $H^1$-conforming linear finite element space, i.e., $V_k$ with $r=1$. Let $\mathcal{V}_k$ be the set of all vertices on $\cT_k$ and $\{\varphi_{x}\}_{\bold{x}\in\mathcal{V}_k}$ be the nodal basis functions in $V^1_k$. For each $\bold{x}\in\mathcal{V}_k$, $\omega_x$ stands for the support of $\varphi_x$, i.e., the union of all elements in $\cT_k$ with $\bold{x}$ being a vertex, and $\cT_k(\omega_x)$ is the triangulation of $\omega_x$	with respect to $\cT_k$. The interpolation operator $\Pi_k:L^1(D)\to V_k^1$ is given by
\begin{equation}\label{def:cn_int}
    \Pi_k v := \sum_{\bold{x}\in\mathcal{N}_k} \frac{1}{|\omega_x|}\int_{\omega_x}v \dx \varphi_x.
\end{equation}
    It is easy to see that $0\leq \Pi_k v \leq 1$ if $0\leq v \leq1$ a.e. in $D$. Moreover, the quadrature rule implies
	\[
		\begin{aligned}
			\int_D \Pi_k v \dx & = \sum_{\bold{x}\in\mathcal{N}_k} \frac{1}{|\omega_x|}\int_{\omega_x} v \dx \int_{\omega_x} \varphi_x \dx = \sum_{\bold{x}\in\mathcal{N}_k} \frac{1}{|\omega_x|}\int_{\omega_x} v \dx \sum_{T\in \cT_k(\omega_x)}\int_{T} \varphi_x \dx \\
			&= \sum_{\bold{x}\in\mathcal{N}_k} \frac{1}{|\omega_x|}\int_{\omega_x} v \dx \sum_{T\in \cT_k(\omega_x)} \frac{|T|}{d+1}= \frac{1}{d+1} \sum_{\bold{x}\in\mathcal{N}_k} \sum_{T\in \cT_k(\omega_x)} \int_T v \dx = \int_D v \dx.
		\end{aligned}
	\]
So $\int_{D}\Pi_k v \dx = V$ if $\int_D v \dx = V$. Consequently, $\Pi_k \phi \in \U_k$ for $\phi\in \U$. Noting $\phi_k^\ast$ is a minimizer to \eqref{eigen-opt_phasefield_disc}-\eqref{eigen_vp_phasefield_disc} and assuming that each $\lambda^{\eps,\ast}_{k,i_j}$ $(1\leq j\leq l)$ is simple, we may assert as in \cite{GarckeHuttlKahleKnopfLaux:2023, GarckeHuttlKnopf:2022} that $\phi_k^\ast$ and associated eigenpairs $\left(\lambda^{\eps,\ast}_{k,i_j}, w_{k,i_j}^\ast\right)$ $(1\leq j\leq l)$ satisfy
	\begin{equation}\label{optsys_phase-field_disc}
			\left\{
			\begin{array}{ll}				\gamma\mathcal{F}_\eps'(\phi^\ast_k)(\phi_k-\phi_k^\ast) + \displaystyle\alpha\sum_{j=1}^l\frac{\p \Psi(\lambda^{\eps,\ast}_{k,i_1},\lambda^{\eps,\ast}_{k,i_2},\cdots,\lambda^{\eps,\ast}_{k,i_l})}{\p \lambda_{i_j}^{\eps}}\int_{D}(w_{k,i_j}^{\ast})^2 (\phi_k -\phi_k^\ast)  \dx  \geq 0 \quad \forall\phi_k \in \U_k,    \\[1ex]
				\displaystyle\int_{D} \bold{\nabla} w_{k,i_j}^\ast \cdot \bold{\nabla} v_k \dx +  \alpha\int_D \phi^\ast_k w^\ast_{k,i_j} v_k \dx = \lambda^{\eps,\ast}_{k,i_j}\int_{D} w_{k,i_j}^\ast v_k \dx \quad  \forall v_k \in  V_k^0,~\forall j \in \{1,2,\cdots, l\}.
			\end{array}
			\right.
		\end{equation}
		
		Now inserting $\phi_k=\Pi_k \phi\in \U_k$ for any $\phi\in\U$ in the first variational inequality of \eqref{optsys_phase-field_disc} and using \eqref{residual_0} result in
		\begin{align}
			\langle \mathcal{R}_0(\phi_k^\ast), \phi - \phi_k^\ast \rangle & =	\gamma\mathcal{F}_\eps'(\phi^\ast_k)(\phi-\Pi_k\phi) + \displaystyle\alpha\sum_{j=1}^l\frac{\p \Psi(\lambda^{\eps,\ast}_{k,i_1},\lambda^{\eps,\ast}_{k,i_2},\cdots,\lambda^{\eps,\ast}_{k,i_l})}{\p \lambda_{i_j}^{\eps}}\int_{D}(w_{k,i_j}^{\ast})^2 (\phi - \Pi_k\phi)  \dx \nonumber \\
			& \quad + \gamma\mathcal{F}_\eps'(\phi^\ast_k)(\Pi_k\phi-\phi_k^\ast) + \displaystyle\alpha\sum_{j=1}^l\frac{\p \Psi(\lambda^{\eps,\ast}_{k,i_1},\lambda^{\eps,\ast}_{k,i_2},\cdots,\lambda^{\eps,\ast}_{k,i_l})}{\p \lambda_{i_j}^{\eps}}\int_{D}(w_{k,i_j}^{\ast})^2 (\Pi_k\phi -\phi_k^\ast)  \dx \nonumber\\
			&\geq  \gamma\mathcal{F}_\eps'(\phi^\ast_k)(\phi-\Pi_k\phi) + \displaystyle\alpha\sum_{j=1}^l\frac{\p \Psi(\lambda^{\eps,\ast}_{k,i_1},\lambda^{\eps,\ast}_{k,i_2},\cdots,\lambda^{\eps,\ast}_{k,i_l})}{\p \lambda_{i_j}^{\eps}}\int_{D}(w_{k,i_j}^{\ast})^2 (\phi - \Pi_k\phi)  \dx. \label{eqn:est_derive01}
		\end{align}
		By the elementwise integration by parts, Cauchy-Schwarz inequality, error estimates for $\Pi_k$ \cite[Lemma 5.3]{JinXu:2020}, we proceed with the right hand side in the inequality of \eqref{eqn:est_derive01}
		\begin{align}
			\left|\langle \mathcal{R}_0(\phi_k^\ast), \phi - \Pi_k\phi_k \rangle\right| & = \left|\sum_{T\in\cT_k} \int_T R_0(\phi_k^\ast) (\phi-\Pi_k\phi) \dx + \sum_{F\in\mathcal{F}_k}\int_F J_0(\phi_k^\ast) (\phi-\Pi_k\phi) \dd s \right| \nonumber\\
			& \leq \sum_{T\in\cT_k} \left( \|R_0(\phi_k^\ast)\|_{L^2(T)}\|\phi-\Pi_k\phi\|_{L^2(T)}  + \sum_{F\subset\p T} \|J_0(\phi_k^\ast)\|_{L^2(F)}\|\phi-\Pi_k\phi\|_{L^2(F)}\right) \nonumber\\
			& \leq c \sum_{T\in\cT_{k}} \left( h_T \|R_0(\phi_k^\ast)\|_{L^2(T)} + h_T^{1/2}\sum_{T\subset\p T}\|J_0(\phi_k^\ast)\|_{L^2(F)}\right)\|\bold{\nabla}\phi\|_{L^2(\omega_k(T))}\nonumber\\
			& \leq c \eta_{k,0}(\phi_k^\ast)\|\phi\|_{H^1(D)}\quad \forall \phi \in \U. \label{eqn:est_derive02}
		\end{align}
		From \eqref{eqn:est_derive01}-\eqref{eqn:est_derive02}, it follows that the residual \eqref{residual_0}, associated with the variational inequality in optimality system \eqref{optsys_phase-field}, is related to  $\eta_{k,0}(\phi_k^\ast)$ although the computable quantity does not provide a reliable bound. This motivates the involvement of $\eta_{k,0}$ in the module \textsf{ESTIMATE}.
		
Noting $\Pi_k$ does not preserve the zero boundary condition in $H_0^1(D)$, next we resort to the Scott-Zhang interpolation operator $I_k: H_0^1(D)\to V_k^0$ \cite{ScottZhang:1990}.  By \eqref{residual_j}, the operator $I_k$ and the second equation of \eqref{optsys_phase-field_disc}  with $v_k = I_k v \in V_k^0$ for any $v\in H_0^1(D)$, similar arguments yield for $1\leq j \leq l$
		\begin{align}
			&\quad \left| \langle  \mathcal{R}_j(\phi_{k}^{\ast},\lambda_{k,i_j}^{\eps,\ast},w_{k,i_j}^{\ast}), v \rangle   \right|  = \left| \textcolor{red}{\langle}\mathcal{R}_j(\phi_{k}^{\ast},\lambda_{k,i_j}^{\eps,\ast},w_{k,i_j}^{\ast}), v - I_k v \rangle   \right| \quad \forall v \in H_0^1(D) \nonumber\\
			& =\left| \sum_{T\in \cT_k}\int_T R_j(\phi_{k,i_j}^{\eps,\ast}, \lambda_{k,i_j}^{\eps,\ast}, w_{k,i_j}^{\ast}) ( v - I_k v) \dx + \sum_{F\in\mathcal{F}_k}\int_F J_j(w_{k,i_j}^{\ast}) ( v - I_k v) \dd s \right| \nonumber\\
			&\leq c\eta_{k,j}(\phi_{k,i_j}^{\eps,\ast}, \lambda_{k,i_j}^{\eps,\ast}, w_{k,i_j}^{\ast}) \|v\|_{H^1(D)}. \label{eqn:est_derive03}
		\end{align}
		As before, each $\eta_{k,j}$ provides an upper bound of the corresponding residual in \eqref{residual_j}, so it naturally comes into use in the module \textsf{ESTIMATE} as well.
		
		It is worth mentioning that even though the above argument fails to induce the reliability of each estimator $\eta_{k,j}$ $(0 \leq j \leq l)$
		in terms of the error bound, \eqref{eqn:est_derive01}-\eqref{eqn:est_derive03} forms the basis of our convergence analysis. In particular, in the next section we are able to show $\{\phi_{k}^\ast,(\lambda_{k,i_j}^{\eps,\ast},w_{k,i_j}^\ast)\}_{k\geq 0}$
		has a strongly convergent subsequence. Then our desired result in Theorem \ref{thm:alg_conv} follows from \eqref{eqn:est_derive01}-\eqref{eqn:est_derive03} once after we prove the zero limit of estimators associated with the convergent subsequence in the right hand sides and the convergence of all residuals in the left hand sides in Section \ref{sec:conv}.
		
		\section{Limiting behavior}\label{sec:limit}
		
		In this section, we establish an auxiliary convergence result for Algorithm \ref{alg_afem_eigenvalue}. Since $\{\U_k\}_{k\geq0}$ and $\{V_k^0\}_{k\geq 0}$ given by Algorithm \ref{alg_afem_eigenvalue} are both nested, it is possible to define two closures in the $H^1(D)$-norm and the $H_0^1(D)$-norm respectively
		\[
		\U_\infty:= \overline{\bigcup_{k\geq 0} \U_k},\quad V_\infty^0:= \overline{\bigcup_{k\geq 0} V^0_k}.
		\]
		Now we consider the following limiting problem
		\begin{equation}\label{eigen-opt_phasefield_limit}
			\begin{aligned}
				\inf_{\phi_\infty\in \U_\infty}\left\{\J_{\infty}^\eps(\phi_\infty):=\Psi(\lambda^{\eps}_{\infty,i_1},\lambda^{\eps}_{\infty,i_2},\cdots,\lambda^{\eps}_{\infty,i_l})+\gamma \mathcal{F}_\eps(\phi_\infty)\right\}~
				\text{subject to}~ \{\lambda^{\eps}_{\infty,i_j}\}_{j=1}^l~ \text{given by}
			\end{aligned}
		\end{equation}
		\begin{equation}\label{eigen_vp_phasefield_limit}
			\begin{aligned}
				\int_{D} \bold{\nabla} w_\infty & \cdot \bold{\nabla} v_\infty \dx + \alpha \int_D \phi_\infty w_\infty v_\infty \dx = \lambda_{\infty}^\eps(\phi_\infty)\int_{D} w_\infty v_\infty \dx\quad  \forall v_\infty \in V^0_\infty \\  &\text{for}~(\lambda^\eps_{\infty}, w_\infty)\in \mathbb{R} \times V_\infty^0 \setminus 0~\text{with} ~\|w_\infty\|_{L^2(D)} = 1.
			\end{aligned}
		\end{equation}
		The goal of this subsection is to prove that the adaptively generated sequence of discrete minimizers has a subsequence convergent to a minimizer of \eqref{eigen-opt_phasefield_limit}-\eqref{eigen_vp_phasefield_limit} in $H^1(D)$. To achieve this, we need some preliminary results, which may be viewed as the counterparts of ingredients in the analysis for the continuous case \cite{GarckeHuttlKahleKnopfLaux:2023, GarckeHuttlKnopf:2022}. First it is not difficult to find from the definition that $V_\infty^0$ is a closed subspace of $H_0^1(D)$. Due to the compact embedding of $H_0^1(D)$ into $L^2(D)$, $V_\infty^0\hookrightarrow L^2(D)$ is also compact. Then by the spectrum theorem for compact self-adjoint operators, we immediately know
		that for any $\phi_\infty\in \U_\infty$, state equation \eqref{eigen_vp_phasefield_limit} has an infinite  nondecreasing sequence $\{\lambda_{\infty,i}^\eps(\phi_\infty)\}_{i\geq1}$ of positive eigenvalues, each of which is counted with its geometric multiplicity, tending to $\infty$ and has a sequence  $\{w_{\infty,i}\}_{i\geq1}\subset V_\infty^0$ of corresponding eigenfunctions forming an orthonormal basis of $L^2(D)$.
		Moreover, there holds the following Courant-Fischer characterization.
		
		\begin{lemma}\label{lem:Courant-Fisher_limit}
			With any $\phi_\infty\in \U_\infty$ given, define a bilinear form and the Rayleigh quotient respectively as
			\[
			a(u,v):= \int_{D}\bold{\nabla}u\cdot\bold{\nabla}v\dx + \alpha\int_{D}\phi_\infty u v\dx\quad \forall u,v\in V_\infty^0,\quad
			\mathcal{R}(v):=\frac{a(v,v)}{\|v\|^2_{L^2(D)}}\quad \forall  v\in V_{\infty}^0\setminus{\{0\}}.
			\]
			%\[
			%    %\mathcal{R}(v):=\frac{\int_{D}|\bold{\nabla}v|^2\dx+\int_{D}\phi_\infty|v|^2\dx}{\|v\|^2_{L^2(D)}}\quad %\forall  v\in V_{\infty}^0\setminus{\{0\}}.
			%\]
			For $j\in \mathbb{N}$, the corresponding $j$-th eigenvalue $\lambda_{\infty,j}^\eps$ in \eqref{eigen_vp_phasefield_limit} is characterized by the maximum-minimum principle
			\begin{equation}\label{max-min_limit}
				\lambda_{\infty,j}^\eps =  \max_{\substack{S_{j-1}\subset V_\infty^0\\ \mathrm{dim}S_{j-1}=j-1}}\min_{v \in S_{j-1}^{\perp}\cap V_\infty^0} \mathcal{R}(v) = \min_{\substack{v\in V_{\infty}^0, a(v,w_{\infty,i}) = 0 \\ i=1,2,\cdots,j-1}}\mathcal{R}(v),
			\end{equation}
			where $S_{j-1}^\perp$ is the orthogonal complement of $S_j\subset L^2(D)$ with respect to the associated scalar product and the maximum is attained at $S_{j-1}^\ast:=\mathrm{span}\{w_{\infty,1}, w_{\infty,2},\cdots,w_{\infty,j-1}\}$, or equivalently, by the minimum-maximum principle
			\begin{equation}\label{min-max_limit}
				\lambda_{\infty,j}^\eps = \min_{\substack{S_j\subset V_\infty^0 \\ \mathrm{dim}S_j=j}} \max_{v\in S_j} \mathcal{R}(v) = \max_{v\in \textrm{span}\{w_{\infty,1},w_{\infty,2},\cdots,w_{\infty,j}\}}\mathcal{R}(v).
			\end{equation}
			
		\end{lemma}
		
		\begin{proof}
			We focus on the second conclusion as the first one can be proved as in Theorem 3.2(b) of \cite{GarckeHuttlKnopf:2022}. As $S_j$ is a $j$ dimensional subspace, it clearly has non-empty intersection with $S=\overline{\mathrm{span}}\{w_{\infty,j}, w_{\infty,j+1},\cdots\}$. Thus any nonzero $\widetilde{v}\in S_j\cap S$ may be written as $\widetilde{v}=\sum_{i=j}^{\infty}\widetilde{\alpha}_{\infty,i} w_{\infty,i}$ with $\widetilde{\alpha}_{\infty,i} = \int_D \widetilde{v} w_{\infty, i} \dx$. Noting each $\lambda_{\infty,i}^\eps$ is positive and nondecreasing with respect to $i$, we deduce
			\begin{equation}\label{lem:Courant-Fisher_limit_pf1}
				\sup_{v\in S_{j}}\mathcal{R}(v)\geq \sup_{\widetilde{v}\in S_{j}\cap S}\mathcal{R}(\widetilde{v})=\sup_{\widetilde{v}\in S_{j}\cap S}\frac{\sum_{i=j}^{\infty}\widetilde{\alpha}_{\infty,i}^2\lambda_{\infty,i}^\eps}{\sum_{i=j}^{\infty}\widetilde{\alpha}_{\infty,i}^2}\geq \lambda_{\infty,j}^\eps.
				%\quad \forall \widetilde{v} \in (S_j\cap %S)\setminus\{0\}.
				%\quad \geq \lambda_{\infty,j}^\eps.
			\end{equation}
			%On the other hand, from \eqref{eigen_vp_phasefield_limit} we know $\left\{\dfrac{w_{\infty,i}}{\sqrt{\lambda^\eps_{\infty,i}}}\right\}_{i\geq 1}$ is an orthonormal basis of $V_\infty^0$ equipped with the inner product $a(\cdot,\cdot)$. In view of this, we assume that the $j$ dimensional subspace $S_j\subset V_{\infty}^0$ is spanned by certain finite eigenfunctions $\{w_{\infty,l_1}, w_{\infty,l_2}, \cdots, w_{\infty,l_j}\}$. To be exact, any nonzero $v\in S_j$  may be represented as
			%\[
			%v=\sum_{k=1}^j \frac{1}{\lambda_{\infty,l_k}^\eps} a(v,w_{\infty,l_k}) w_{\infty,l_k}  = \sum_{k=1}^j \int_{D} v w_{\infty,l_k} \dx w_{\infty,l_k} = \sum_{k=1}^j \alpha_{\infty,k} w_{\infty,l_k}.
			%\]
			%This further implies that
			%\begin{equation}\label{lem:Courant-Fisher_limit_pf2}
			%	\mathcal{R}(v) = \frac{\sum_{k=1}^{j}\alpha_{\infty,k}^2\lambda_{\infty,l_k}^\eps}{\sum_{k=1}^{j}\alpha_{\infty,k}^2}\leq \lambda_{\infty,l_j}^\eps\quad \forall v \in S_j\setminus\{0\}.
			%\end{equation}
			%Since $\mathcal{R}(v)$ is bounded from above over $S_j$,  there exists a maximizing sequence $\left\{v_i\right\}_{i\geq 1}\subset S_{j}\setminus \{0\}$ such that $\mathcal{R}(v_i)\to \sup_{v\in S_{j}}\mathcal{R}(v)$ as $i\to\infty$ and consequently $\{\left\|v_i/\|v_i\|_{L^2(D)}\right\|_{H^1(D)}\}_{i\geq1}$ is bounded. Noting the finite dimensionality of $S_j$ and the continuity of the $H^1_0(D)$-norm, we can use a compactness argument to obtain a maximizer to $\mathcal{R}(v)$ over $S_j$. 
            Moreover, as $S_j$ is arbitrary, it follows from \eqref{lem:Courant-Fisher_limit_pf1} that
			\begin{equation}\label{lem:Courant-Fisher_limit_pf3}
				\inf_{\substack{S_j\subset V_\infty^0 \\ \mathrm{dim}S_j=j}} \sup_{v\in S_j} \mathcal{R}(v) \geq \lambda_{\infty,j}^\eps.
			\end{equation}
On the other hand, taking $S_j=S^\ast_j:=\mathrm{span}\{w_{\infty,1},w_{\infty,w_2},\cdots,w_{\infty,j}\}$, we obtain $$\mathcal{R}(v) = \frac{\sum_{i=1}^{j}\alpha_{\infty,i}^2\lambda_{\infty,i}^\eps}{\sum_{i=1}^{j}\alpha_{\infty,i}^2}\leq \mathcal{R}(w_{\infty,j}) =\lambda_{\infty,j}^\eps$$ for any nonzero $v\in S^\ast_j$, which, together with \eqref{lem:Courant-Fisher_limit_pf3}, yields the first equation in \eqref{min-max_limit} and indicates that the minimum is attained at $S^\ast_j$, i.e., the second equation in \eqref{min-max_limit}.
            %Now taking $S_j:=S^\ast_j:=\mathrm{span}\{w_{\infty,1},w_{\infty,w_2},\cdots,w_{\infty,j}\}$ in \eqref{lem:Courant-Fisher_limit_pf2}, we obtain $\mathcal{R}(v)\leq \lambda_{\infty,j}^\eps$ for any nonzero $v\in S^\ast_j$, which, together with \eqref{lem:Courant-Fisher_limit_pf3}, yields the first equation in \eqref{min-max_limit} and indicates that the minimum is attained at $S^\ast_j$, i.e., the second equation in \eqref{min-max_limit}.
		\end{proof}
		
		Next we turn to an estimate of eigenvalues of \eqref{eigen_vp_phasefield_disc}.

		\begin{lemma}\label{lem:eigenvalue_disc_estimate}
			Let $\left\{\U_k \times V_k^0\right\}_{k\geq0}$ be as given by Algorithm \ref{alg_afem_eigenvalue}, $(\widetilde{\lambda}_{j},\widetilde{w}_j)\in\mathbb{R}\times H_0^1(\Omega)$ be the $j$-th eigenpair of the Laplacian operator satisfying
			\begin{equation}\label{eigen_vp_Laplace}
				\int_{D} \bold{\nabla} \widetilde{w}_j\cdot \bold{\nabla} v \dx = \widetilde{\lambda}_{j} \int_D \tilde{w}_j v \dx \quad \forall v \in H_0^1(D).
			\end{equation}
			and $\left\{(\widetilde{\lambda}_{k,j},\widetilde{w}_{k,j})\right\}_{k \geq K} \subset \mathbb{R} \times \bigcup_{k \geq K} V_k^0$, with nonnegative integer $K$ being the minimum index such that $j\leq \mathrm{dim}V_K^0$,
			be the sequence of the $j$-th eigenpairs of discrete problems
			\begin{equation}\label{eigen_vp_Laplace_disc}
				\int_{D} \bold{\nabla} \widetilde{w}_{k,j}\cdot \bold{\nabla} v_k \dx = \widetilde{\lambda}_{k,j} \int_D \widetilde{w}_{k,j} v_k \dx \quad \forall v_k \in V_k^0.
			\end{equation}
			Then for any sequence $\{\phi_k\}_{k\geq0}\subset\bigcup_{k\geq0}\U_k$ and the corresponding sequence $\{\lambda_{k,j}^{\eps}\}_{k\geq0}$ of the $j$-th eigenvalues, each given by \eqref{eigen_vp_phasefield_disc} over $V_k^0$, there holds
			\begin{equation}\label{eigenvalue_disc_estimate}
				\widetilde{\lambda}_{j}\leq \lambda^\eps_{k,j} \leq \widetilde{\lambda}_{K,j} + \alpha\quad \forall k\geq K.
			\end{equation}
			% $\widetilde{\lambda}_{0,j}$ is the %$j$-eigenvalue over
		\end{lemma}
		
		\begin{proof}
			The proof follows that of Lemma 3.7 in \cite{GarckeHuttlKahleKnopfLaux:2023}. As $\phi_k\in [0,1]$, it is straightforward to see
			\[
			\frac{\|\bold{\nabla}v_k\|^2_{L^2(D)}}{\|v_k\|^2_{L^2(D)}}\leq \frac{\int_{D}|\bold{\nabla}v_k|^2\dx+\alpha\int_{D}\phi_k|v_k|^2\dx}{\|v_k\|^2_{L^2(D)}} \leq \frac{\|\bold{\nabla}v_k\|^2_{L^2(D)}}{\|v_k\|^2_{L^2(D)}} +
			\alpha \quad \forall v_k\in V_k^0\setminus\{0\}.
			\]
			So the minimum-maximum principles for \eqref{eigen_vp_Laplace_disc} and \eqref{eigen_vp_phasefield_disc} over $V_k^0$ imply
			\begin{equation}\label{lem:eigenvalue_disc_estimate_pf01}
				\widetilde{\lambda}_{k,j} \leq \lambda_{k,j}^\eps \leq \widetilde{\lambda}_{k,j} + \alpha \quad \forall k\geq K.
			\end{equation}
			As $\{V_k^0\}_{k\geq0}$ are nested closed subspaces of $H_0^1(\Omega)$, by the minimum-maximum principles for \eqref{eigen_vp_Laplace} and \eqref{eigen_vp_Laplace_disc} again we know $\{\widetilde{\lambda}_{k,j}\}_{k\geq K}$ is a decreasing sequence bounded below by $\widetilde{\lambda}_j$, from which and \eqref{lem:eigenvalue_disc_estimate_pf01} we infer \eqref{eigenvalue_disc_estimate}.
		\end{proof}
		
		Now we are ready to show the continuity of discrete eigenpairs $\{(\lambda_{k,j}^{\eps}, w_{k,j})\}_{k\geq 0} \subset \mathbb{R}_{>0} \times \bigcup_{k\geq0}V^0_{k}  $ with respect to convergent discrete phase-fields $\{\phi_k\}_{k\geq 0} \subset \bigcup_{k\geq0} \U_k$ as $k\to\infty$.
		As in \cite{GarckeHuttlKnopf:2022}, the analysis is first performed for the first eigenvalue of \eqref{eigen_vp_phasefield_limit} and then extended to other ones by induction. But our argument is modified in the current situation since we deal with a sequence of discrete eigenpairs
		instead of continuous ones in \cite{GarckeHuttlKnopf:2022}.
		
		\begin{lemma}\label{lem:eigen-1_weak_cont}
			If $\{\phi_k\}_{k\geq0}\subset\bigcup_{k\geq0} \U_k$ converges to some $\phi_\infty\in \U_\infty$ weakly in $H^1(D)$, then for the sequence   $\left\{\left(\lambda^\eps_{k,1}(\phi_k),w_{k,1}(\phi_k)\right)\right\}_{k\geq 0}\subset \mathbb{R}_{>0}\times \bigcup_{k\geq 0} V_k^0$ of the first eigenpairs to \eqref{eigen_vp_phasefield_disc} over $\{\cT_k\}_{k\geq0}$, $\left\{\lambda^\eps_{k,1}(\phi_k)\right\}_{k\geq0}$  converges to the first eigenvalue $\lambda_{\infty,1}^\eps(\phi_\infty)$ of \eqref{eigen_vp_phasefield_limit}, i.e.,
			\begin{equation}\label{eigenvalue-1_weak_cont}
				\lambda_{k,1}^\eps \to \lambda_{\infty,1}^\eps \quad \text{as}~k \to\infty,
			\end{equation}
			and there exists a corresponding eigenfunction $\overline{w}_{\infty,1}(\phi_\infty)\in V_\infty^0$ associated with $\lambda_{\infty,1}^\eps$ such that up to a subsequence,
			\begin{equation}\label{eigenfunction-1_weak_cont}
				w_{k,1}\to \overline{w}_{\infty,1}\quad \text{strongly in}~H^1_0(D)\quad \text{as}~ k \to \infty.%,\quad w_{k,1}\to w_{\infty,1}\quad %\text{strongly in}~L^2(\Omega).
			\end{equation}
			%as $k\to\infty$.
			%$ %solving  respectively and the first eigenpair %$(\lambda_{\infty,1}^\ast(\phi_\infty),w_{\infty,1}(\phi_\infty))\in\mathbb{R}\times %V_\infty^0$ to \eqref{eigen_vp_phasefield_limit}.    %there hold %eigenpair $(\lambda^\eps_{\infty,1}, %w_{\infty,1}) %\in %\mathbb{R}\times V_\infty^0$ such %that as $k\to %\infty$
			
		\end{lemma}
		
		\begin{proof}
			In view of Lemma \ref{lem:eigenvalue_disc_estimate}, $\left\{\lambda_{k,1}^\eps\right\}_{k\geq0}$ is a bounded sequence. From \eqref{eigen_vp_phasefield_disc}, we further know that $\|\bold{\nabla}w_{k,1}\|_{L^2(D)}\leq \lambda_{k,1}^\eps$, which implies $\left\{\|w_{k,1}\|_{H^1(D)}\right\}_{k\geq0}$ is bounded. Since $V_\infty$ is a closed subspace of $H_0^1(D)$, by the reflexivity of $H^1_0(D)$ and Sobolev compact embedding theorem there exists a subsequence $\left\{w_{k_n,1}\right\}_{n \geq 0}$ and some $\overline{w}_{\infty,1}\in V_\infty^0$ such that as $n\to\infty$,
			\begin{equation}\label{lem:eigen-1_weak_cont_pf1}
				w_{k_n,1} \rightharpoonup \overline{w}_{\infty,1} \quad \text{weakly in}~H^1_0(D),\quad w_{k_n,1} \to \overline{w}_{\infty,1} \quad \text{strongly in}~L^2(D).
			\end{equation}
			By Sobolev compact embedding theorem again, the assumption of the $H^1(D)$ weak convergence on $\{\phi_k\}_{k\geq0}$ implies the strong convergence in $L^2(D)$, which further allows us to extract a subsequence $\{\phi_{k_n}\}_{n \geq 0}$ with
			\begin{equation}\label{lem:eigen-1_weak_cont_pf2}
				\phi_{k_n} \to \phi_\infty \quad \text{a.e.~in}~D \quad \text{as}~n \to \infty.
			\end{equation}
			The rest of the proof consists of three steps. We establish the convergence of $\left\{\lambda^\eps_{k,1}\right\}_{k\geq0}$ in the first and the second steps while the third step discusses the convergence of corresponding eigenfunctions.
			
			\noindent \textit{Step 1.}
			By Lebesgue's dominated convergence theorem, the $L^2(D)$ strong convergence in \eqref{lem:eigen-1_weak_cont_pf1} and the pointwise convergence in \eqref{lem:eigen-1_weak_cont_pf2}, we may deduce that as $n \to \infty$,
			\begin{align}
				\left|\int_D (\phi_{k_n} w_{k_n,1}^2 - \phi_\infty \overline{w}_{\infty,1}^2)  \dx \right| & \leq  \left|\int_D \phi_{k_n} ( w_{k_n,1}^2 - \overline{w}_{\infty,1}^2 ) \dx \right| + \left| \int_D ( \phi_{k_n}  - \phi_\infty ) \overline{w}_{\infty,1}^2 \dx \right| \nonumber \\
				& \leq \int_D | w_{k_n,1}^2 - \overline{w}_{\infty,1}^2 | \dx + \int_D | ( \phi_{k_n}  - \phi_\infty ) \overline{w}_{\infty,1}^2 | \dx \to 0.\label{lem:eigen-1_weak_cont_pf3}
			\end{align}
			This, together with the $H^1_0(D)$ weak convergence in \eqref{lem:eigen-1_weak_cont_pf1} and \eqref{eigen_vp_phasefield_disc} again, implies
			\begin{align*}
				\liminf_{n\to\infty} \lambda_{k_n,1}^\eps & = \liminf_{n\to\infty} \left(\int_D|\bold{\nabla}w_{k_n,1}|^2\dx+ \alpha\int_D \phi_{k_n} w_{k_n,1}^2\dx \right) \\
				&= \liminf_{n\to\infty}\int_D|\bold{\nabla}w_{k_n,1}|^2\dx + \alpha\lim_{n\to\infty}\int_D \phi_{k_n} w_{k_n,1}^2\dx \geq
				\int_D |\bold{\nabla}\overline{w}_{\infty,1}|^2 \dx + \alpha\int_D \phi_{\infty} \overline{w}_{\infty,1}^2\dx.
			\end{align*}
			On the other hand, $\|w_{k_n,1}\|_{L^2(\Omega)}=1$ and the strong convergence in \eqref{lem:eigen-1_weak_cont_pf1} ensure $\|\overline{w}_{\infty,1}\|_{L^2(\Omega)}=1$. So from the maximum-minimum principle for $\lambda^\eps_{\infty,1}$ in Lemma \ref{lem:Courant-Fisher_limit}, which reduces to $\lambda_{\infty,1}^\eps=\min _{v\in V_\infty^0} \mathcal{R}(v)$, it follows that
			\begin{equation}\label{lem:eigen-1_weak_cont_pf4}
				\liminf_{n\to\infty}\lambda_{k_n,1}^\eps \geq \lambda_{\infty,1}^\eps.
			\end{equation}
			
			\noindent{\textit{Step 2.}} According to the definition, any $v_\infty \in V_\infty^0$ admits a strongly convergent sequence $\{v_k\}_{k\geq0}\subset\bigcup_{k\geq0}V_k^0$ with respect to the $H^1_0(D)$-topology. Then by the maximum-minimum principle for $\lambda_{k_n,1}^\eps$,
			\[
			\lambda_{k_n,1}^\eps \leq \frac{\int_D |\bold{\nabla}v_{k_n}|^2 \dx + \alpha\int_D \phi_{k_n}v_{k_n}^2  \dx}{\int_D |v_{k_n}|^2 \dx}.
			\]
			Passing to the limit $n\to\infty$, arguing as in Step 1 and using the maximum-minimum principle for $\lambda_{\infty,1}^\eps$ again, we obtain
			\begin{equation}\label{lem:eigen-1_weak_cont_pf5}
				\limsup_{n\to\infty} \lambda_{k_n,1}^\eps \leq \lambda_{\infty,1}^\eps.
			\end{equation}
			Combining \eqref{lem:eigen-1_weak_cont_pf4} with \eqref{lem:eigen-1_weak_cont_pf5} yields $\displaystyle\lim_{n\to\infty}\lambda_{k_n,1}^\eps \to \lambda_{\infty,1}^\eps$, from which we deduce  \eqref{eigenvalue-1_weak_cont} for the whole sequence by the standard subsequence contradiction argument.
			
			\noindent\textit{Step 3.} We shall prove that the weak limit $\overline{w}_{\infty,1}\in V_\infty^0$ in \eqref{lem:eigen-1_weak_cont_pf1} solves \eqref{eigen_vp_phasefield_limit} as an eigenfunction associated with $\lambda_{\infty,1}^\eps$. As $\{V_k^0\}_{k\geq 0}$ is nested, then for $k_n\geq K$ with some $K\in \mathbb{N}_0$ fixed
			\[
			\int_D \bold{\nabla} w_{k_n,1} \cdot \bold{\nabla} v_{K} \dx + \alpha \int_D \phi_{k_n} w_{k_n,1} v_K \dx = \lambda^\eps_{k_n,1} \int_D w_{k_n,1} v_K \dx \quad \forall v_K \in V_{K}^0.
			\]
			Letting $n$ tend to $\infty$, we immediately know
			\[
			\int_D \bold{\nabla} w_{k_n,1} \cdot \bold{\nabla} v_{K} \dx \to \int_D \bold{\nabla} \overline{w}_{\infty,1} \cdot \bold{\nabla} v_{K} \dx, \quad
			\lambda^\eps_{k_n,1} \int_D w_{k_n,1} v_K \dx \to \lambda^\eps_{\infty,1} \int_D \overline{w}_{\infty,1} v_K \dx.
			\]
			Further \eqref{lem:eigen-1_weak_cont_pf1}, \textcolor{blue}{\eqref{lem:eigen-1_weak_cont_pf2}} and Lebesgue's dominated convergence theorem imply that as $n\to\infty$,
			\[
			\begin{aligned}
				\left| \int_D ( \phi_{k_n} w_{k_n,1} - \phi_{\infty} \overline{w}_{\infty,1} ) v_K \dx \right|  & = \left| \int_D \phi_{k_n} ( w_{k_n,1} - \overline{w}_{\infty,1} )  v_K \dx
				+ \int_D (\phi_{k_n} - \phi_{\infty}) \overline{w}_{\infty,1}   v_K \dx \right| \\
				& \leq \|w_{k_n,1} - \overline{w}_{\infty,1}\|_{L^2(\Omega)}\|v_K\|_{L^2(\Omega)} + \int_D |(\phi_{k_n} - \phi_{\infty}) \overline{w}_{\infty,1}   v_K| \dx \to 0.
			\end{aligned}
			\]
			Therefore, it follows that
			\[
			\int_D \bold{\nabla} \overline{w}_{\infty,1} \cdot \bold{\nabla} v_{K} \dx + \alpha \int_D \phi_{\infty} \overline{w}_{\infty,1} v_K \dx = \lambda^\eps_{\infty,1} \int_D \overline{w}_{\infty,1} v_K \dx \quad \forall v_K \in V_{K}^0.
			\]
			By the density of $\bigcup_{k\geq 0}V_k^0$ in $V_\infty^0$, we arrive at
			\begin{equation}\label{lem:eigen-1_weak_cont_pf6}
				\int_D \bold{\nabla} \overline{w}_{\infty,1} \cdot \bold{\nabla} v_{\infty} \dx + \alpha \int_D \phi_{\infty} \overline{w}_{\infty,1} v_\infty \dx = \lambda^\eps_{\infty,1} \int_D \overline{w}_{\infty,1} v_\infty \dx \quad \forall v_\infty \in V_{\infty}^0.
			\end{equation}
			Now taking $v_\infty = \overline{w}_{\infty,1}$ in \eqref{lem:eigen-1_weak_cont_pf6} and recalling  $\|\overline{w}_{\infty,1}\|_{L^2(\Omega)}=1$, we have
			\[
			\|\bold{\nabla}\overline{w}_{\infty,1}\|^2_{L^2(\Omega)} + \alpha \int_{D} \phi_\infty \overline{w}_{\infty,1}^2 \dx = \lambda_{\infty,1}^\eps.
			\]
			From \eqref{eigen_vp_phasefield_disc} again and \eqref{eigenvalue-1_weak_cont}, it follows that
			$\|\bold{\nabla}w_{k_n,1}\|^2_{L^2(\Omega)} + \alpha \int_{D} \phi_{k_n} w_{k_n,1}^2 \dx = \lambda_{k_n,1}^\eps \to \lambda_{\infty,1}^\eps$ as $n\to\infty$. Taking account of \eqref{lem:eigen-1_weak_cont_pf3}, we obtain $\|\bold{\nabla}w_{k_n,1}\|^2_{L^2(\Omega)}\to \|\bold{\nabla}\overline{w}_{\infty,1}\|^2_{L^2(\Omega)}$ as $n\to\infty$, which, along with the $H^1_0(D)$ weak convergence in \eqref{lem:eigen-1_weak_cont_pf1} again, leads to \eqref{eigenfunction-1_weak_cont}.
		\end{proof}
		
		\begin{lemma}\label{lem:eigen_weak_cont}
			Under the same assumption of Lemma \ref{lem:eigen-1_weak_cont}, for the sequence    $\left\{\left(\lambda^\eps_{k,j}(\phi_k),w_{k,j}(\phi_k)\right)\right\}_{k\geq 0}\subset \mathbb{R}_{>0}\times \bigcup_{k\geq 0} V_k^0$ of the $j$-th eigenpairs ($j\in \mathbb{N}$) to \eqref{eigen_vp_phasefield_disc} over $\{\cT_k\}_{k\geq0}$, $\left\{\lambda^\eps_{k,j}(\phi_k)\right\}_{k\geq0}$  converges to the $j$-th eigenvalue $\lambda_{\infty,j}^\eps(\phi_\infty)$ of \eqref{eigen_vp_phasefield_limit}, i.e.,
			\begin{equation}\label{eigenvalue_weak_cont}
				\lambda_{k,j}^\eps \to \lambda_{\infty,j}^\eps \quad \text{as}~k\to\infty,
			\end{equation}
			and there exists a corresponding eigenfunction $\overline{w}_{\infty,j}(\phi_\infty)\in V_\infty^0$ associated with $\lambda_{\infty,j}^\eps$ such that up to a subsequence,
			\begin{equation}\label{eigenfunction_weak_cont}
				w_{k,j}\to \overline{w}_{\infty,j}\quad \text{strongly in}~H^1_0(D)\quad  \text{as}~k\to\infty.%,\quad w_{k,1}\to w_{\infty,1}\quad %\text{strongly in}~L^2(\Omega).
			\end{equation}
		\end{lemma}
		
		\begin{proof}
			By induction, \eqref{eigenvalue_weak_cont} and \eqref{eigenfunction_weak_cont} for $j=1$ are available in Lemma \ref{lem:eigen-1_weak_cont}. Assuming the conclusion is true for $j - 1\in \mathbb{N}$, we shall verify the two desired results for $j$ in the rest of the proof.
			
			From Lemma \ref{lem:eigenvalue_disc_estimate}, \eqref{eigen_vp_phasefield_disc}, Sobolev compact embedding theorem and the weak convergence of $\{\phi_k\}_{k\geq0}$, we may extract subsequences $\left\{w_{k_n,j}\right\}_{n \geq 0}$ and $\left\{\phi_{k_n}\right\}_{n \geq 0}$ as in the proof of Lemma \ref{lem:eigen-1_weak_cont} such that as $n\to\infty$,
			\begin{equation}\label{lem:eigen_weak_cont_pf1}
				w_{k_n,j} \rightharpoonup \overline{w}_{\infty,j} \quad \text{weakly in}~H^1_0(D),\quad w_{k_n,j} \to \overline{w}_{\infty,j} \quad \text{strongly in}~L^2(D),
				\quad \phi_{k_n} \to \phi_\infty \quad \text{a.e.~in}~D
			\end{equation}
			with some $L^2(D)$-normalized $\overline{w}_{\infty,j}\in V_\infty^0$. Since $\lambda^\eps_{k_n,j}$ is the $j$-th eigenvalue of \eqref{eigen_vp_phasefield_disc} over $\cT_{k_n}$, we argue as in Step 1 of the proof of Lemma \ref{lem:eigen-1_weak_cont} to obtain
			\begin{equation}\label{lem:eigen_weak_cont_pf2}
				\liminf_{n\to\infty}\lambda_{k_n,j}^\eps \geq  \int_D |\bold{\nabla}\overline{w}_{\infty,j}|^2 \dx + \alpha \int_D \phi_{\infty} \overline{w}_{\infty,j}^2\dx.
			\end{equation}
			By the induction hypothesis, for each $i\in \{1,2,\cdots,j-1\}$ the sequence $\{w_{k,i}\}_{k\geq0}$ of the $i$-th discrete eigenfunctions given by  \eqref{eigen_vp_phasefield_disc} over $\{\cT_{k}\}_{k\geq0}$ contains a subsequence $\{w_{k_n,i}\}_{n \geq 0}$ converging strongly in $H_0^1(D)$  to some $\overline{w}_{\infty,i}\in V_\infty^0$ associated with the $i$-th eigenvalue $\lambda_{\infty,i}^\eps$ of \eqref{eigen_vp_phasefield_limit}, namely
			\begin{equation}\label{lem:eigen_weak_cont_pf3}
				w_{k_n,i}\to \overline{w}_{\infty,i}\quad \text{strongly in}~H_0^1(D)\quad \text{as}~n \to \infty.
			\end{equation}
			As each sequence $\{w_{k,i}\}_{i=1}^{\mathrm{dim}V_{k}^0}$ is $L^2(D)$-orthonormal, it follows from \eqref{lem:eigen_weak_cont_pf3} that $\{\overline{w}_{\infty,i}\}_{i=1}^{j-1}$ is also $L^2(D)$-orthonormal and consequently linearly independent. With $j^\ast\leq j-1$ being the maximal index such that $\lambda_{\infty, j^\ast}^\eps < \lambda_{\infty,j}^\eps$ and the finite $L^2(D)$-orthonormal basis $\{w_{\infty,i}\}_{i=1}^{j^\ast}$ given by \eqref{eigen_vp_phasefield_limit} invoked, the associated eigenvalues of which are identical with those of $\{\overline{w}_{\infty,i}\}_{i=1}^{j^\ast}$, we can deduce that
			\begin{equation}\label{lem:eigen_weak_cont_pf4} \mathrm{span}\{\overline{w}_{\infty,1},\overline{w}_{\infty,2},\cdots,\overline{w}_{\infty,j^\ast}\} = \mathrm{span}\{w_{\infty,1},w_{\infty,2},\cdots,w_{\infty,j^\ast}\}.
			\end{equation}
			On the other hand, the $L^2(D)$-orthogonality of $\{w_{k,i}\}_{i=1}^{\mathrm{dim}V_{k}^0}$ and  \eqref{lem:eigen_weak_cont_pf1} imply that
			\begin{equation}\label{lem:eigen_weak_cont_pf5}
				\int_{D} \overline{w}_{\infty,j}\overline{w}_{\infty,i} \dx=0\quad i\in \{1,2,\cdots, j-1\}.
			\end{equation}
			Now combining \eqref{lem:eigen_weak_cont_pf4} and \eqref{lem:eigen_weak_cont_pf5} leads to
			\[
			\overline{w}_{\infty,j} = \sum_{i=1}^\infty \overline{\alpha}_i w_{\infty,i} = \sum_{i=j^\ast+1}^\infty \overline{\alpha}_i w_{\infty,i} \quad \overline{\alpha}_i=\int_D \overline{w}_{\infty,j}w_{\infty,i} \dx.
			\]
			Noting $\|\overline{w}_{\infty,j}\|_{L^2(D)}=1$ and the series in the above representation converges in $H_0^1(D)$, inserting it into the right hand side of \eqref{lem:eigen_weak_cont_pf2} and using \eqref{eigen_vp_phasefield_limit}, we proceed
			\begin{equation}\label{lem:eigen_weak_cont_liminf}
				\liminf_{n\to\infty}\lambda_{k_n,j}^\eps \geq  \sum_{i=j^\ast+1}^\infty \lambda_{\infty,i}^\eps \overline{\alpha}_i^2\geq \lambda_{\infty,j}^\eps.
			\end{equation}
			On the other hand, for any nonzero $\overline{v}_\infty\in \widetilde{S}_{j-1}^\perp\cap V_\infty^0:=(\mathrm{span}\{\overline{w}_{\infty,1},\overline{w}_{\infty,2},\cdots, \overline{w}_{\infty,j-1}\})^\perp\cap V_\infty^0$, the construction of $V_\infty^0$ allows the existence of a sequence $\{v_k\}_{k\geq0}\subset \bigcup_{k\geq 0}V_k^0$ with the property
			\begin{equation}\label{lem:eigen_weak_cont_pf6}
				v_k\to \overline{v}_\infty \quad \text{strongly in}~ H_0^1(D) \quad  \text{as}~ k \to \infty.
			\end{equation}
			Defining $\overline{v}_k:=v_k - \sum_{i=1}^{j-1} \alpha_{k,i} w_{k,i}$ with $\alpha_{k,i} = \int_{D} v_k w_{k,i} \dx$, we immediately know $\overline{v}_k\in (S^\ast_{k,j-1})^\perp\cap V_k^0:=(\mathrm{span}\{w_{k,1},w_{k,2},\cdots,
			w_{k,j-1}\})^\perp \cap V_k^0$. More importantly, in view of the maximum-minimum principle for the $j$-th eigenvalue associated with \eqref{eigen_vp_phasefield_disc} over $\cT_k$, it holds
			\begin{equation}\label{lem:eigen_weak_cont_pf7}
				\lambda_{k,j}^\eps = \min_{v \in (S_{k,j-1}^\ast)^\perp\cap V_k^0 } \frac{\int_{D} |\bold{\nabla} v |^2 \dx + \alpha\int_{D}\phi_k v^2\dx}{\|v\|^2_{L^2(D)}}\leq \frac{\int_{D} |\bold{\nabla} \overline{v}_k |^2 \dx + \alpha\int_{D}\phi_k \overline{v}_k^2\dx}{\|\overline{v}_k\|^2_{L^2(D)}}.
			\end{equation}
			Next we consider the limit of the right hand side of \eqref{lem:eigen_weak_cont_pf7} along the subsequence $\left\{\overline{v}_{k_n}\right\}_{n \geq 0}$. The definition of $\overline{v}_{k_n}$ leads to
			\[
			\begin{aligned}
				&\quad \int_{D} |\bold{\nabla} \overline{v}_{k_n}|^2 \dx + \alpha \int_{D} \phi_{k_n} \overline{v}_{k_n}^2  \dx  =
				\int_{D} |\bold{\nabla} v_{k_n}|^2 \dx + \alpha \int_{D} \phi_{k_n} v_{k_n}^2 \dx \\
				& - 2 \sum_{i=1}^{j-1} \alpha_{k_n,i} \left( \int_D \bold{\nabla} w_{k_n,i} \cdot \bold{\nabla} v_{k_n} \dx  +  \alpha \int_D \phi_{k_n}  w_{k_n,i} v_{k_n}  \dx \right) \\
				& + \int_D \sum_{i=1}^{j-1} \alpha_{k_n,i} \bold{\nabla} w_{k_n,i} \cdot \sum_{i=1}^{j-1} \alpha_{k_n,i} \bold{\nabla} w_{k_n,i} \dx + \alpha \int_{D} \phi_{k_n} \left( \sum_{i=1}^{j-1} \alpha_{k_n,i}  w_{k_n,i} \right)^2 \dx.
			\end{aligned}
			\]
			Using the strong convergence in \eqref{lem:eigen_weak_cont_pf6}, the pointwise convergence in \eqref{lem:eigen_weak_cont_pf1}, Lebegue's dominated convergence theorem and the argument for \eqref{lem:eigen-1_weak_cont_pf3}, we obtain
			\[
			\int_{D} |\bold{\nabla} v_{k_n}|^2 \dx \to \int_{D} |\bold{\nabla} \overline{v}_\infty|^2 \dx, \quad
			\int_{D} \phi_{k_n} v_{k_n}^2 \dx \to \int_{D} \phi_{\infty} \overline{v}_{\infty}^2 \dx \quad \text{as}~n\to\infty.
			\]
			As $\left\{w_{k_n,i}\right\}_{i=1}^{j-1}$ is $L^2(D)$-orthonormal, it follows from \eqref{eigen_vp_phasefield_disc} over $\cT_{k_n}$ that
			\[
			\sum_{i=1}^{j-1} \alpha_{k_n,i} \left( \int_D \bold{\nabla} w_{k_n,i} \cdot \bold{\nabla} v_{k_n} \dx +  \alpha \int_D \phi_{k_n}  w_{k_n,i} v_{k_n}  \dx \right)  = \sum_{i=1}^{j-1} \left( \int_D v_{k_n} w_{k_n,i} \dx \right)^2\lambda_{k_n,i}^\eps,
			\]
			\[
			\int_D \left|\sum_{i=1}^{j-1} \alpha_{k_n,i} \bold{\nabla} w_{k_n,i} \right|^2 \dx + \alpha \int_{D} \phi_{k_n} \left( \sum_{i=1}^{j-1} \alpha_{k_n,i}  w_{k_n,i} \right)^2 \dx = \sum_{i=1}^{j-1} \left( \int_D v_{k_n} w_{k_n,i} \dx \right)^2\lambda_{k_n,i}^\eps.
			\]
			Invoking the induction hypothesis that $\displaystyle\lim_{k\to\infty}\lambda_{k,i}^\eps =\lambda_{\infty,i}^\eps$ for $i=1,2,\cdots,j-1$,   \eqref{lem:eigen_weak_cont_pf3} again as well as  \eqref{lem:eigen_weak_cont_pf6} and noting $\overline{v}_\infty\in \widetilde{S}_{j-1}^\perp$, we have
			\[
			\sum_{i=1}^{j-1} \left( \int_D v_{k_n} w_{k_n,i} \dx \right)^2\lambda_{k_n,i}^\eps \to \sum_{i=1}^{j-1} \left( \int_D \overline{v}_{\infty} \overline{w}_{\infty,i} \dx \right)^2\lambda_{\infty,i}^\eps = 0 \quad \text{as} ~ n\to \infty.
			\]
			As a result of the above arguments, it follows
			\begin{equation}\label{lem:eigen_weak_cont_pf8}
				\int_{D} |\bold{\nabla} \overline{v}_{k_n}|^2 \dx + \alpha \int_{D} \phi_{k_n} \overline{v}_{k_n}^2  \dx \to
				\int_{D} |\bold{\nabla} \overline{v}_\infty|^2 \dx + \alpha \int_{D} \phi_{\infty} \overline{v}_{\infty}^2 \dx \quad \text{as} ~ n\to \infty.
			\end{equation}
			Likewise, \eqref{lem:eigen_weak_cont_pf6} and the induction hypothesis yield
			\begin{align}
				\|\overline{v}_{k_n}\|_{L^2(D)}^2 & = \|v_{k_n}\|_{L^2(D)}^2 - 2 \sum_{i=1}^{j-1}\alpha_{k_n,i} \int_{D} v_{k_n} w_{k_n,i} \dx  + \int_{D} \left( \sum_{i=1}^{j-1} \alpha_{k_n,i}  w_{k_n,i} \right)^2 \dx \nonumber \\
				& = \|v_{k_n}\|_{L^2(D)}^2 -  \sum_{i=1}^{j-1}\left( \int_{D} v_{k_n} w_{k_n,i} \dx\right)^2 \to  \|\overline{v}_{\infty}\|_{L^2(D)}^2 \quad \text{as} ~ n\to \infty. \label{lem:eigen_weak_cont_pf9}
			\end{align}
			Now passing to the limit along the subsequence index  $n\to\infty$ in \eqref{lem:eigen_weak_cont_pf7}, we find from \eqref{lem:eigen_weak_cont_pf8} and \eqref{lem:eigen_weak_cont_pf9} that for any nonzero $ \overline{v}_\infty\in \widetilde{S}_{j-1}^\perp\cap V_\infty^0$,
			\[
			\limsup_{n\to\infty}\lambda_{k_n,j}^\eps \leq \limsup_{n\to\infty} \frac{\int_{D} |\bold{\nabla} \overline{v}_{k_n}|^2 \dx + \alpha \int_{D} \phi_{k_n} \overline{v}_{k_n}^2  \dx }{\|\overline{v}_{k_n}\|_{L^2(D)}^2} = \frac{\int_{D} |\bold{\nabla} \overline{v}_\infty|^2 \dx + \alpha\int_{D} \phi_{\infty} \overline{v}_{\infty}^2 \dx}{\|\overline{v}_{\infty}\|_{L^2(D)}^2},
			\]
			which implies that $\displaystyle\limsup_{n\to\infty}\lambda_{k_n,j}^\eps \leq \inf_{\overline{v}_\infty\in \widetilde{S}_{j-1}^\perp\cap V_\infty^0}\mathcal{R}(\overline{v}_\infty)$. The infimum in the right hand side allows us to choose a minimizing sequence $\left\{v_i\right\}_{i\geq 1}\subset (\widetilde{S}_{j-1}^\perp\cap V_\infty^0)\setminus \{0\}$ such that $\displaystyle \mathcal{R}(v_i)\to \inf_{\overline{v}_\infty\in \widetilde{S}_{j-1}^\perp\cap V_\infty^0}\mathcal{R}(\overline{v}_\infty)$ as $i\to\infty$. Since $\left\{\left\|v_i/\|v_i\|_{L^2(D)}\right\|_{H^1(D)}\right\}_{i\geq1}$ is bounded, the compactness argument ensures the infimum is attained. Therefore, we infer from the maximum-minimum principle \eqref{max-min_limit} in Lemma \ref{lem:Courant-Fisher_limit} that
			\begin{equation}\label{lem:eigen_weak_cont_limsup}
				\limsup_{n\to\infty}\lambda_{k_n,j}^\eps \leq \min_{\overline{v}_\infty\in \widetilde{S}_{j-1}^\perp\cap V_\infty^0}\mathcal{R}(\overline{v}_\infty)
				\leq \lambda_{\infty,j}^\eps.
			\end{equation}
			As a result of \eqref{lem:eigen_weak_cont_liminf} and  \eqref{lem:eigen_weak_cont_limsup}, $\displaystyle\lim_{n\to\infty}\lambda_{k_n,j}^\eps = \lambda_{\infty,j}^\eps$. Considering the limit is independent of the subsequence, the conclusion is true for the whole sequence. Finally, the existence of an eigenfunction of \eqref{eigen_vp_phasefield_limit} associated with $\lambda_{\infty,j}^\eps$ and the strong convergence with respect to the $H_0^1(D)$-topology in \eqref{eigenfunction_weak_cont} can be justified analogously as in Lemma \ref{lem:eigen-1_weak_cont} (see Step 3 of the proof).
		\end{proof}
		
		\begin{theorem}\label{thm:alg_conv_limit}
			The sequence $\{\phi_k^\ast\}_{k\geq0}$ of discrete minimizers to \eqref{eigen-opt_phasefield_disc}-\eqref{eigen_vp_phasefield_disc} generated by Algorithm \ref{alg_afem_eigenvalue} contains a subsequence converging to a minimizer $\phi_\infty^\ast$ to \eqref{eigen-opt_phasefield_limit}-\eqref{eigen_vp_phasefield_limit} with respect to the $H^1(D)$-topology. Moreover, each of the corresponding sequences $\left\{\left(\lambda_{k,i_j}^{\eps,\ast}, w^\ast_{k,i_j}\right)\right\}_{k\geq0}$ ($j=1,2,\cdots,l$) of discrete eigenpairs has a subsequence converging to the $i_j$-th  eigenpair $\left(\lambda_{\infty,i_j}^{\eps,\ast},w^\ast_{\infty,i_j}\right)$ of \eqref{eigen_vp_phasefield_limit} associated with $\phi_\infty^\ast$ in the sense that
			\begin{equation}\label{alg_conv_limit_eigen}
				\lambda_{k_n, i_j}^{\eps,\ast}\to  \lambda_{\infty,i_j}^{\eps,\ast},\quad w_{k_n,i_j}^{\ast} \to w_{\infty,i_j}^\ast\quad  \text{strongly in}~H^1_0(D) \quad \text{as}~n\to\infty.
			\end{equation}
			
		\end{theorem}
		
\begin{proof}
    The definition clearly implies that $\U_\infty$ is closed in $H^1(D)$. For any $\phi_\infty$, $\psi_\infty\in \U_\infty$, there exist two sequences $\{\phi_k\}_{k\geq0}$ and $\{\psi_k\}_{k\geq0}$ in $\bigcup_{k\geq0}\U_{k}$ such that $\phi_k\to\phi_\infty$ and $\psi_k\to\psi_\infty$ strongly in $H^1(D)$ as $k\to\infty$. Thanks to the convexity of each $\U_k$, $\{t\phi_k+(1-t)\psi_k\}_{k\geq0}\subset \bigcup_{k\geq0}\U_{k}$ for any $t\in (0,1)$. Since $\U_\infty$ is the $H^1(D)$ closure of $\bigcup_{k\geq0}\U_{k}$, $t\phi_k+(1-t)\psi_k\to t\phi_\infty+(1-t)\psi_\infty \in \U_\infty $ as $k\to\infty$ for any $t\in (0,1)$. So $\U_\infty$ is convex. Further, the strong convergence of $\{\phi_k\}_{k\geq0}$ in $H^1(D)$ implies that  $\phi_k\to \phi$ a.e. in $D$ and $\int_D \phi_k \dx \to \int_D \phi_\infty \dx $ as $k\to\infty$ along a subsequence. Noting the constraints $0\leq \phi_k\leq 1$ a.e. in $D$ and $\int_D \phi_k\dx = V$, we immediately know that $0\leq \phi_\infty\leq 1$ a.e. in $D$ and $\int_D \phi_\infty\dx = V$.  In summary, $\U_\infty$ is a closed convex subset of $\U$.
			
			It is easy to see that $\frac{V}{|D|}$ belongs to $\U_k$ for all $k\in \mathbb{N}_0$. Lemma \ref{lem:eigenvalue_disc_estimate} implies that sequences $\left\{\lambda_{k,i_j}^\eps(\phi^\ast_k)\right\}_{k\geq0}$ and $\left\{\lambda_{k,i_j}^\eps(\frac{V}{|D|})\right\}_{k\geq0}$ ($j=1,2,\cdots,l$) of eigenvalues are all positive and  uniformly bounded from above and below. Recalling $\Psi$ is a continuous function on $(\R_{>0})^l$ and taking $\phi_k = \frac{V}{|D|}$ in  $\J_k^\eps$ , without loss of generality
			we may deduce that there exist two positive constants $c_1$ and $c_2$ such that 
			\[
			-c_1 + \gamma \mathcal{F}_\eps(\phi_k^\ast) \leq \J_{k}^\eps(\phi_k^\ast) \leq \J_k^\eps\left(\frac{V}{|D|}\right)\leq c_2 + \frac{\gamma}{\eps}\int_{D} f\left(\frac{V}{|D|}\right)\dx.
			\]
			Hence, $\left\{|\phi_k^\ast|_{H^1(D)}\right\}_{k\geq0}$ is bounded. In addition, by H\"{o}lder inequality, it holds $\|\phi_k^\ast\|_{L^2(D)}\leq \sqrt{|D|}$ for each $k\in \mathbb{N}_0$. So $\left\{\phi_k^\ast\right\}_{k\geq0}$ is a bounded sequence in $H^1(D)$. As $\U_{\infty}$ is closed and convex, the reflexivity of $H^1(D)$ and Sobolev compact embedding theorem guarantee a subsequence $\left\{\phi_{k_n}^\ast\right\}_{n \geq 0}$ and some $\overline{\phi}_{\infty}\in\U_\infty$ such that as $n\to\infty$,
			\begin{equation}\label{thm:alg_conv_limit_pf01}
				\phi^\ast_{k_n}\rightharpoonup \overline{\phi}_\infty\quad \text{weakly in}~H^1(D), \quad \phi^\ast_{k_n}\to\overline{\phi}_\infty\quad \text{a.e. in}~D.
			\end{equation}
			Then Lemma \ref{lem:eigen_weak_cont} allows us to extract, from corresponding eigenpairs $\left\{\left(\lambda_{k,i_j}^{\eps,\ast},w_{k,i_j}^{\ast}\right)\right\}_{k\geq0}$ ($1\leq j\leq l$) given by Algorithm \ref{alg_afem_eigenvalue}, a subsequence $\left\{\left(\lambda_{k_n,i_j}^{\eps,\ast},w_{k_n,i_j}^{\ast}\right)\right\}_{n \geq 0}$, which strongly converges to the $i_j$-th ($1\leq j\leq l$) eigenpair $\left(\lambda_{\infty,i_j}^{\eps,\ast},w^\ast_{\infty,i_j}\right)\in \mathbb{R}_{>0} \times V_\infty^0 $ of \eqref{eigen_vp_phasefield_limit} with $\phi_\infty = \overline{\phi}_\infty$, i.e.,
			\begin{equation}\label{thm:alg_conv_limit_pf02}
				\lambda_{k_n,i_j}^{\eps,\ast}\to\lambda_{\infty,i_j}^{\eps,\ast},\quad
				w_{k_n,i_j}^{\ast}\to w_{\infty,i_j}^{\ast}\quad\text{strongly in}~H_0^1(D) \quad \text{as}~ n\to\infty.
			\end{equation}
			Then from the continuity of $\Psi$, the convergence of $\left\{\lambda_{k_n,i_j}^{\eps,\ast}\right\}_{n \geq 0}$ ($1\leq j\leq l$) in \eqref{thm:alg_conv_limit_pf02}, the pointwise convergence of $\left\{\phi_{k_n}^\ast\right\}_{n \geq 0}$ in \eqref{thm:alg_conv_limit_pf01} and Lebegue's dominated convergence theorem, we infer
			\begin{equation}\label{thm:alg_conv_limit_pf03}
				\Psi\left(\lambda_{k_n,i_1}^{\eps,\ast},\lambda_{k_n,i_2}^{\eps,\ast}\cdots,\lambda_{k_n,i_l}^{\eps,\ast}\right)\to \Psi\left(\lambda_{\infty,i_1}^{\eps,\ast},\lambda_{\infty,i_2}^{\eps,\ast},\cdots,\lambda_{\infty,i_l}^{\eps,\ast}\right),\quad \int_D f(\phi_{k_n}^\ast) \dx \to \int_D f(\overline{\phi}_\infty) \dx \quad \text{as}~ n\to\infty.
			\end{equation}
			On the other hand, any $\phi_\infty\in \U_\infty$ entails the strong $H^1(D)$ convergence of a sequence $\{\phi_{k}\}_{k\geq0}\subset\bigcup_{k\geq0}\U_k$. By Lemma \ref{lem:eigen_weak_cont} again and the standard subsequence contradiction argument, we deduce as above
			\begin{equation}\label{thm:alg_conv_limit_pf04}
                \Psi\left(\lambda_{k,i_1}^\eps,\lambda_{k,i_2}^\eps,\cdots,\lambda_{k,i_l}^\eps\right)\to \Psi\left(\lambda_{\infty,i_1}^\eps,\lambda_{\infty,i_2}^\eps,\cdots,\lambda_{\infty,i_l}^\eps\right),\quad
				\int_D f(\phi_{k}) \dx \to \int_D f(\phi_\infty) \dx \quad \text{as}~n\to\infty
			\end{equation}
			with $\left\{\left(\lambda_{\infty,i_j}^\eps,\overline{w}_{\infty,i_j}\right)\right\}_{j=1}^l$ solving \eqref{eigen_vp_phasefield_limit}. Now using the weak convergence in \eqref{thm:alg_conv_limit_pf01}, the weak lower semi-continuity of $H^1(D)$ semi-norm,  \eqref{thm:alg_conv_limit_pf03}, the minimizing property of $\{\J^\eps_k(\phi_k^\ast)\}$ over $\U_k$ and \eqref{thm:alg_conv_limit_pf04}, we obtain
			\[
			\begin{aligned}
				\J_\infty^\eps(\overline{\phi}_\infty) \leq \liminf_{n\to\infty}\J_{k_n}^\eps(\phi_{k_n}^\ast) \leq \limsup_{n\to\infty}\J_{k_n}^\eps(\phi_{k_n}^\ast)
				\leq \limsup_{k\to\infty}\J_{k}^\eps(\phi_{k}^\ast)
				\leq \limsup_{k\to\infty}\J_{k}^\eps(\phi_{k})=\J_\infty(\phi_\infty)\quad \forall \phi_\infty\in \U_\infty.
			\end{aligned}
			\]
			The existence of a minimizer to \eqref{eigen-opt_phasefield_limit}-\eqref{eigen_vp_phasefield_limit} is given by $\phi^\ast_\infty:=\overline{\phi}_\infty\in \U_\infty$ and the second assertion of the theorem holds due to \eqref{thm:alg_conv_limit_pf02}. Inserting $\phi_\infty=\phi_\infty^\ast$ in the last inequality, we further deduce
			\[
			\lim_{n\to\infty}\J_{k_n}^\eps(\phi_{k_n}^\ast)  = \J_\infty^\eps(\phi^\ast_\infty),
			\]
			which, along with \eqref{thm:alg_conv_limit_pf03} again, yields $\|\bold{\nabla}\phi_{k_n}^\ast\|_{L^2(D)}^2 \to \|\bold{\nabla}\phi_{\infty}^\ast\|_{L^2(D)}^2 $ as $n\to\infty$. A combination of this and the weak convergence in \eqref{thm:alg_conv_limit_pf01} concludes the proof.
		\end{proof}
		
		\begin{remark}\label{rem:alg_conv_limit}
			From \eqref{thm:alg_conv_limit_pf01} in the above proof, we can also assert that $\phi_{k_n}^\ast\to \phi_\infty^\ast$ almost everywhere in $D$ as $n\to\infty$.
		\end{remark}
		
		\section{Convergence}\label{sec:conv}
		%\subsection{Main result}\label{subsect:conv}
		
		This section is devoted to the convergence of Algorithm \ref{alg_afem_eigenvalue} without any stopping criterion.  As mentioned in the introduction, it includes two results involving the convergent subsequence
		\[
		\left\{\phi_{k_n}^\ast, \left(\lambda_{k_n,i_1}^{\eps,\ast},w_{k_n,i_1}\right),\cdots,\left(\lambda_{k_n,i_l}^{\eps,\ast},w_{k_n,i_l}\right)\right\}_{n \geq 0}
		\]
		from Theorem \ref{thm:alg_conv_limit}. We shall first prove that the $l+1$ relevant subsequences $\{\eta_{k_{n},j}\}_{n\geq0}$ ($0\leq j\leq l$) of  estimators given by Algorithm \ref{alg_afem_eigenvalue} all converge to zero as $n\to\infty$. In this step, the approach without recourse to reliability and efficiency of the estimator \cite{GanterPraetorius:2022} is adapted to our case. Then the limit
		\[
		\left\{\phi_\infty^\ast, \left(\lambda_{\infty,i_1}^{\eps,\ast},w_{\infty,i_1}\right), \cdots, \left(\lambda_{\infty,i_l}^{\eps,\ast},w_{\infty,i_l}\right)\right\}\in \U_\infty\times (\mathbb{R}_{>0}\times V_\infty^0)^l
		\]
		given by Theorem \ref{thm:alg_conv_limit} will be shown to solve the optimality system \eqref{optsys_phase-field}. The modules \textsf{ESTIMATE}, \textsf{MARK} and \textsf{REFINE} in Algorithm \ref{alg_afem_eigenvalue}, which are not utilized in proving Theorem \ref{thm:alg_conv_limit}, play an essential role at this stage. We note that these three modules all entail the mesh-size function $h_k$ defined in \eqref{meshsize_def}. The latter is linked to a decomposition of the sequence $\{\cT_k\}_{k\geq0}$ of meshes generated by Algorithm \ref{alg_afem_eigenvalue} as in \cite{Siebert:2011}
		\[
		\mathcal{T}_{k}^{+}:=\bigcap_{m\geq k}\mathcal{T}_{m},\quad
		\mathcal{T}_{k}^{0}:=\mathcal{T}_{k}\setminus\mathcal{T}_{k}^{+},\quad
		D_{k}^{+}:= \bigcup_{T\in\mathcal{T}^{+}_{k}}T,\quad
		D_{k}^{0}:=\bigcup_{T\in\mathcal{T}^{0}_{k}}T.
		\]
		The set $\mathcal{T}_{k}^{+}$ consists of all elements not refined after the $k$-th iteration while
		all elements in $\mathcal{T}_{k}^{0}$ are refined at least once after the $k$-th iteration. It is not difficult to observe that $\mathcal{T}_m^+\subset\mathcal{T}_n^+$ for $m<n$ and $\mathcal{M}_k\subseteq\cT_{k}^0$. Moreover, with $\chi_{k}^0$ denoting the characteristic function of $D_k^0$, the mesh-size function $h_k$ has the property \cite{MorinSiebertVeeser:2008} \cite{Siebert:2011}
		\begin{equation}\label{mesh-size->zero}
			\lim_{k\to\infty}\|h_k\chi_{k}^0\|_{L^\infty(D)} = 0.
		\end{equation}
		
		By the above convergence and the auxiliary convergence in Theorem \ref{thm:alg_conv_limit}, we can investigate the behavior of each maximal error indicator among $\mathcal{M}_{k}^j$ of the marked set in Algorithm \ref{alg_afem_eigenvalue}.
		
\begin{lemma}\label{lem:max-error-ind->zero}
Let $\left\{\phi_{k_n}^\ast, \left(\lambda_{k_n,i_j}^{\eps,\ast},w^\ast_{k_n,i_j}\right)_{j=1}^l\right\}_{n\geq 0}$ be the convergent subsequence given by Theorem \ref{thm:alg_conv_limit}. Over the sequence $\left\{\mathcal{M}_{k_n}\right\}_{n \geq 0}$ of corresponding marked sets determined in Algorithm \ref{alg_afem_eigenvalue}, there holds
	\begin{equation}\label{max-error-ind-zero}
				\lim_{n\to\infty} \max_{T\in\mathcal{M}_{k_n}}\eta_{k_n,j}(T) = 0\quad 0\leq j\leq l.
			\end{equation}
		\end{lemma}
		
\begin{proof}
As \eqref{marking} suggests, $T_{k_n}^j$ $(0 \leq j \leq l)$ is also the element with the largest error indicator among all $\eta_{k_n,j}(T)$ over $\mathcal{M}_{k_n}$. Due to the local quasi-uniformity of $\cT_{k_n}$, $T_{k_n}^j\in\mathcal{M}_{k_n}\subseteq\cT_{k_n}^0$ and \eqref{mesh-size->zero}, there holds
			\begin{equation}\label{lem:max-error-ind->zero_pf01}
    |\omega_{k_n}(T_{k_n}^j)|\leq c |T_{k_n}^j| \leq c \|h_{k_n}\chi_{k_n}^0\|^d_{L^\infty(D)}\to 0\quad \text{as}~ n\to\infty,~0\leq j\leq l.
\end{equation}
By the definition of $\eta_{k_n,j}(T_{k_n}^j)$ in the module \textsf{MARK} of Algorithm \ref{alg_afem_eigenvalue}, we use the inverse estimate and the scaled trace theorem to proceed with two different cases separately
\begin{equation*}
	\begin{aligned}
     \eta^2_{k_n,0}(T_{k_n}^0) & \leq c \Bigg( h_{T_{k_n}^0}^2\sum_{j=1}^l \left|\frac{\partial \Psi(\lambda_{k_n,i_1}^{\eps,\ast},\cdots,\lambda_{k_n,i_l}^{\eps,\ast})}{\partial\lambda_{i_{j}}^\eps}\right|^2\|w_{k_n,i_{j}}^\ast\|^4_{L^4(T_{k_n}^0)} + \|\bold{\nabla}\phi_{k_n}^\ast\|_{L^2(\omega_{k_n}(T_{k_n}^0))}^2 \\
				& \qquad\qquad  + h_{T^0_{k_n}}^2 \|f'(\phi_{k_n}^\ast)\|^2_{L^2(T^0_{k_n})}\Bigg),\\
        \eta^2_{k_n,j}(T_{k_n}^j)  &\leq c \left( \|\bold{\nabla}w_{k_n,i_j}^{\ast}\|^2_{L^2(\omega_{k_n}(T^j_{k_n}))} +  h_{T_{k_n}^j}^2\|(\alpha\phi_{k_n}^\ast-\lambda_{k_n,i_j}^{\eps,\ast})w_{k_n,i_j}^\ast\|_{L^2(T_{k_n}^j)}^2
				\right)\quad 1\leq j\leq l.
	\end{aligned}
\end{equation*}
Since $\Psi$ is $C^1$ and each sequence $\{\lambda_{k_n,i_j}^{\eps,\ast}\}_{n\geq0}$ ($1\leq j\leq l$) is bounded due to Lemma \ref{lem:eigenvalue_disc_estimate}, then each  $\left|\frac{\partial \Psi(\lambda_{k_n,i_1}^{\eps,\ast},\cdots,\lambda_{k_n,i_l}^{\eps,\ast})}{\partial\lambda_{i_{j}}^\eps}\right|$ ($1\leq j \leq l$) is also bounded. The right hand side of the above first inequality is further estimated by
\[
	\begin{aligned}
        \eta^2_{k_n,0}(T_{k_n}^0) & \leq c \bigg(\|h_{k_n}\chi_{k_n}^0\|^2_{L^\infty(D)}\sum_{j=1}^l\left(\|w_{k_n,i_j}^\ast-w_{\infty,i_j}^\ast\|^4_{L^4(D)}+\|w^\ast_{\infty,i_j}\|^4_{L^4(T_{k_n}^0)}\right)          \\
        & \qquad + \|\bold{\nabla} (\phi_{k_n}^\ast - \phi_{\infty}^\ast)\|_{L^2(D)}^2  + \|\bold{\nabla} \phi_{\infty}^\ast\|_{L^2(\omega_{k_n}(T_{k_n}^0))}^2   \\
        & \qquad + \|h_{k_n}\chi_{k_n}^0\|^2_{L^\infty(D)} \left(\|f'(\phi_{k_n}^\ast)-f'(\phi_{\infty}^\ast)\|_{L^2(D)}^2+\|f'(\phi_{\infty}^\ast)\|_{L^2(T_{k_n}^0)}\right)
		\bigg).
	\end{aligned}
\]
In view of Sobolev embedding theorem and Theorem \ref{thm:alg_conv_limit}, $\sum_{j=1}^l\|w_{k_n,i_j}^\ast-w_{\infty,i_j}^\ast\|^4_{L^4(D)}$ and $\|\bold{\nabla} (\phi_{k_n}^\ast- \phi_{\infty}^\ast)\|_{L^2(D)}^2$ both tend to zero as $n\to\infty$. Thanks to the pointwise convergence of $\left\{\phi_{k_n}^\ast\right\}_{n\geq0}$ as already pointed out in Remark \ref{rem:alg_conv_limit} and the continuity of $f'$ over $[0,1]$, Lebesgue's dominated convergence theorem implies that $\|f'(\phi_{k_n}^\ast)-f'(\phi_{\infty}^\ast)\|_{L^2(D)}^2$ also tends to zero as $n\to\infty$. The remaining three terms all tend to zero by \eqref{lem:max-error-ind->zero_pf01} and the absolute continuity of norms $\|\cdot\|_{L^4(D)}$ and $\|\cdot\|_{L^2(D)}$ with respect to the Lebesgue measure. On the other hand, the convergent sequence  $\left\{\lambda_{k_n,i_j}^{\eps,\ast}\right\}_{n\geq0}$ ($1\leq j\leq l$) is bounded. Noting $\phi_{k_n}^\ast\in [0,1]$, we further obtain for $1\leq j \leq l$,
\[
    \begin{aligned}
        \eta^2_{k_n,j}(T_{k_n}^j)  & \leq c \left( \|\bold{\nabla}(w_{k_n,i_j}^{\ast} - w_{\infty,i_j}^{\ast})\|^2_{L^2(D)} + \|\bold{\nabla}w_{\infty,i_j}^{\ast}\|^2_{L^2(\omega_{k_n}(T^j_{k_n}))}\right.\\
        & \qquad \left. + \|h_{k_n}\chi_{k_n}^0\|^2_{L^\infty(D)} \left(\|w_{k_n,i_j}^{\ast} - w_{\infty,i_j}^{\ast}\|^2_{L^2(D)}+\|w_{\infty,i_j}^{\ast}\|^2_{L^2(T^j_{k_n})}\right) \right).
	\end{aligned}
\]
As before, the right hand side vanishes by Theorem \ref{thm:alg_conv_limit}, \eqref{lem:max-error-ind->zero_pf01} and the absolute continuity of $\|\cdot\|_{L^2(D)}$ with respect to the Lebesgue measure when passing to the limit $n\to\infty$. In conclusion, the desired result is proved.
\end{proof}
		
Next we shall borrow the approach from \cite{GanterPraetorius:2022} to establish $\eta_{k_n,j}\to 0$ as $n\to\infty$ for each $j\in \{0,1,\cdots,l\}$. Instead of reliability and efficiency, key ingredients in the argument are stability on non-refined elements and reduction on refined elements of the sequence $\{\eta_{k,j}\}_{k\geq0}$ ($0\leq j\leq l$) given by Algorithm \ref{alg_afem_eigenvalue}, which are stated in the following two lemmas.
		\begin{lemma}\label{lem:est_stability}
			For the sequence $\left\{ \phi_{k}^\ast, \left(\lambda_{k,i_j}^{\eps,\ast},w_{k,i_j}^\ast\right)_{j=1}^l \right\}_{k\geq0}$ of discrete solutions given by Algorithm \ref{alg_afem_eigenvalue}, the corresponding estimators $\{\eta_{k,j}\}_{k\geq0}$ $(0\leq j\leq l)$ are stable in the sense that for $m<n$,
			\begin{equation}\label{est-0_stab}
				\begin{aligned}
					\eta_{n,0}(\cT_{m}\cap\cT_{n}) & \leq \eta_{m,0}(\cT_{m}\cap\cT_{n}) + c \Bigg ( \|\bold{\nabla}(\phi_{n}^\ast - \phi_{m}^\ast)\|_{L^2(D)} + \|f'(\phi_{n}^\ast) - f'(\phi_m^\ast)\|_{L^2(D)}  \\
					& \qquad  + \sum_{j=1}^{l}\left\|\frac{\Psi(\lambda_{n,i_1}^{\eps,\ast},\cdots,\lambda_{n,i_l}^{\eps,\ast})}{\partial \lambda_{i_j}^{\eps}}(w_{m,i_j}^\ast)^2-\frac{\Psi(\lambda_{m,i_1}^{\eps,\ast},\cdots,\lambda_{m,i_l}^{\eps,\ast})}{\partial \lambda_{i_j}^{\eps}}(w_{n,i_j}^\ast)^2\right\|_{L^2(D)}
					\Bigg),
				\end{aligned}
			\end{equation}
			\begin{equation}\label{est-j_stab}
				\begin{aligned}
					\eta_{n,j} (\cT_{m}\cap\cT_{n}) & \leq \eta_{m,j} (\cT_{m}\cap\cT_{n}) + c \Big(\|\bold{\nabla}(w_{n,i_j}^\ast-w_{m,i_j}^\ast)\|_{L^2(D)} \\
					& \qquad  + \|(\alpha\phi_{n}^\ast-\lambda_{n,i_j}^{\eps,\ast})w_{n,i_j}^\ast - (\alpha\phi_{m}^\ast-\lambda_{m,i_j}^{\eps,\ast})w_{m,i_j}^\ast\|_{L^2(D)}\Big)\quad 1\leq j\leq l.
				\end{aligned}
			\end{equation}
		\end{lemma}
		
		\begin{proof}
			By the triangle inequality and the inverse estimate, the proof is straightforward.
		\end{proof}

		\begin{lemma}\label{lem:est_reduction}
			For the sequence $\left\{ \phi_{k}^\ast, \left(\lambda_{k,i_j}^{\eps,\ast},w_{k,i_j}^\ast\right)_{j=1}^l \right\}_{k\geq0}$ of discrete solutions given by Algorithm \ref{alg_afem_eigenvalue}, the corresponding estimators $\{\eta_{k,j}\}_{k\geq0}$ $(0\leq j\leq l)$ are contractive in the sense that there exist $q_{j}\in(0,1)$ $(0\leq j\leq l)$ such that for $m<n$,
			\begin{equation}\label{est-0_reduction}
				\begin{aligned}
					\eta^2_{n,0}(\cT_{n}\setminus\cT_m) & \leq q_{0} \eta^2_{m,0}(\cT_{m}\setminus\cT_n) + c \Bigg ( \|\bold{\nabla}(\phi_{n}^\ast - \phi_{m}^\ast)\|_{L^2(D)}^2 + \|f'(\phi_{n}^\ast) - f'(\phi_m^\ast)\|_{L^2(D)}^2  \\
					& \qquad  + \sum_{j=1}^{l}\left\|\frac{\Psi(\lambda_{n,i_1}^{\eps,\ast},\cdots,\lambda_{n,i_l}^{\eps,\ast})}{\partial \lambda_{i_j}^{\eps}}(w_{n,i_j}^\ast)^2-\frac{\Psi(\lambda_{m,i_1}^{\eps,\ast},\cdots,\lambda_{m,i_l}^{\eps,\ast})}{\partial \lambda_{i_j}^{\eps}}(w_{m,i_j}^\ast)^2\right\|^2_{L^2(D)}
					\Bigg),
				\end{aligned}
			\end{equation}
			\begin{equation}\label{est-j_reduction}
				\begin{aligned}
					\eta_{n,j}^2 (\cT_{n}\setminus\cT_{m}) & \leq q_{j}\eta_{m,j}^2 (\cT_{m}\setminus\cT_{n}) + c \Big(\|\bold{\nabla}(w_{n,i_j}^\ast-w_{m,i_j}^\ast)\|^2_{L^2(D)}\\
					& \qquad  + \|(\alpha\phi_{n}^\ast-\lambda_{n,i_j}^{\eps,\ast})w_{n,i_j}^\ast - (\alpha\phi_{m}^\ast-\lambda_{m,i_j}^{\eps,\ast})w_{m,i_j}^\ast\|_{L^2(D)}^2\Big)\quad 1\leq j\leq l.
				\end{aligned}
			\end{equation}
		\end{lemma}
		
		\begin{proof}
			From the triangle inequality and the inverse estimate, it follows that on each $T\in\cT_{n}\setminus\cT_m$,
			\[
			\begin{aligned}
				\eta_{n,0}(\phi_{n}^\ast;T) &\leq \eta_{n,0}(\phi_{m}^\ast;T) + c \left ( \|\bold{\nabla}(\phi_{n}^\ast - \phi_{m}^\ast)\|_{L^2(\omega_{\cT_n}(T))} + \|f'(\phi_{n}^\ast) - f'(\phi_m^\ast)\|_{L^2(T)} \right. \\
				& \qquad \left. + \sum_{j=1}^{l}\left\|\frac{\Psi(\lambda_{n,i_1}^{\eps,\ast},\cdots,\lambda_{n,i_l}^{\eps,\ast})}{\partial \lambda_{i_j}^{\eps}}(w_{n,i_j}^\ast)^2-\frac{\Psi(\lambda_{m,i_1}^{\eps,\ast},\cdots,\lambda_{m,i_l}^{\eps,\ast})}{\partial \lambda_{i_j}^{\eps}}(w_{m,i_j}^\ast)^2\right\|_{L^2(T)}
				\right),
			\end{aligned}
			\]
			which, together with Young inequality, implies
			\[
			\begin{aligned}
				\eta_{n,0}^2(\phi_{n}^\ast;T) &\leq (1+\delta)\eta_{n,0}^2(\phi_{m}^\ast;T) + c (1+\delta^{-1}) \Bigg( \|\bold{\nabla}(\phi_{n}^\ast - \phi_{m}^\ast)\|^2_{L^2(\omega_{\cT_n}(T))} + \|f'(\phi_{n}^\ast) - f'(\phi_m^\ast)\|^2_{L^2(T)} \\
				& \qquad  + \sum_{j=1}^{l}\left\|\frac{\Psi(\lambda_{n,i_1}^{\eps,\ast},\cdots,\lambda_{n,i_l}^{\eps,\ast})}{\partial \lambda_{i_j}^{\eps}}(w_{n,i_j}^\ast)^2-\frac{\Psi(\lambda_{m,i_1}^{\eps,\ast},\cdots,\lambda_{m,i_l}^{\eps,\ast})}{\partial \lambda_{i_j}^{\eps}}(w_{m,i_j}^\ast)^2\right\|^2_{L^2(T)}
				\Bigg)\quad \forall \delta >0.
			\end{aligned}
			\]
			Now summing over all elements $T\in\cT_{n}\setminus\cT_m$ and using the finite overlapping property of $\omega_{\cT_{n}}(T)$, we further obtain
\[
    \begin{aligned}
        \eta_{n,0}^2(\phi_{n}^\ast;\cT_n\setminus\cT_m) &\leq (1+\delta)\eta_{n,0}^2(\phi_{m}^\ast;\cT_n\setminus\cT_m) + c (1+\delta^{-1}) \Bigg( \|\bold{\nabla}(\phi_{n}^\ast - \phi_{m}^\ast)\|^2_{L^2(D)} + \|f'(\phi_{n}^\ast) - f'(\phi_m^\ast)\|^2_{L^2(D)}  \\
		& \qquad  +\sum_{j=1}^{l}\left\|\frac{\Psi(\lambda_{n,i_1}^{\eps,\ast},\cdots,\lambda_{n,i_l}^{\eps,\ast})}{\partial \lambda_{i_j}^{\eps}}(w_{n,i_j}^\ast)^2-\frac{\Psi(\lambda_{m,i_1}^{\eps,\ast},\cdots,\lambda_{m,i_l}^{\eps,\ast})}{\partial \lambda_{i_j}^{\eps}}(w_{m,i_j}^\ast)^2\right\|^2_{L^2(D)}
		\Bigg).
	\end{aligned}
\]
Since $\cT_n$ is a refinement of $\cT_m$, each $T\in \cT_{n}\setminus\cT_m$ is generated from some $T'\in \cT_{m}\setminus\cT_n$ by possible successive bisections. Noting $\phi_m^\ast\in V_m$ and the bisection rule, we find $J_{0}(\phi_{m}^\ast)=0$ across all $F\in\mathcal{F}_{n}(D)$, which are also in the interior of $T'$, and $h_T^d=|T|\leq \frac{1}{2}|T'| = \frac{1}{2}h_{T'}^d$. Thus it gives
\[
    \eta_{n,0}^2(\phi_{m}^\ast;\cT_n\setminus\cT_m) \leq 2^{-1/d} \eta_{m,0}^2(\phi_{m}^\ast;\cT_m\setminus\cT_n).
\]
We conclude \eqref{est-0_reduction} with $q_0=(1+\delta)2^{-1/d}\in (0,1)$ for sufficiently small $\delta>0$. The other estimate can be proved in a similar manner.
\end{proof}
		
Now we are in a position to present the first main result of this paper.
		
\begin{theorem}\label{thm:estimator->zero}
For the convergent subsequence $\left\{\phi_{k_n}^\ast, \left(\lambda_{k_n,i_j}^{\eps,\ast},w_{k_n,i_j}^\ast\right)_{j=1}^l\right\}_{n\geq 0}$ in Theorem \ref{thm:alg_conv_limit}, the corresponding $l+1$ sequences of estimators $\{\eta_{k_n,j}\}_{ n \geq 0}$ $(0\leq j\leq l)$ all converge to zero as $n\to\infty$.
\end{theorem}
		
\begin{proof}
    The proof is inspired by arguments for Theorem 3.1 in \cite{GanterPraetorius:2022}, Theorem 2.1 in \cite{MorinSiebertVeeser:2008} and Theorem 2.2 in \cite{Siebert:2011} with some modifications. As $\cT_{k_n}$ is a refinement of $\cT_{k_m}$ for $m<n$, there holds $\cT_{k_m}^+:=\bigcap_{j\geq k_{m}}\cT_{j}\subset\cT_{k_m}\cap\cT_{k_n}$. This allows us to define an enlarged set  $\cT^+_{k_m\to k_n}:=\{T\in \cT_{k_n}\cap\cT_{k_m}~|~T\cap D_{k_m}^+\neq \emptyset\}\supseteq\cT_{k_m}^+$. The uniform shape regularity of $\{\cT_k\}_{k\geq0}$ guarantees a more important observation on $\cT^+_{k_m\to k_n}$
	\begin{equation}\label{thm:estimator->zero_pf01}
		\# \cT^+_{k_m\to k_n}\leq c \# \cT_{k_m}^+
	\end{equation}
	with $c$ depending only on $\cT_0$.
	Now we  split $\cT_{k_n}$ and $\eta_{k_n,i_j}^2$ $(0\leq j\leq l)$ respectively as %($0\leq j\leq l$)
\[
	\begin{aligned}
\cT_{k_n}&=(\cT_{k_n}\setminus\cT_{k_m})\cup(\cT_{k_n}\cap\cT_{k_m})
				=(\cT_{k_n}\setminus\cT_{k_m})\cup \left(\cT_{k_n}\cap\left(\left(\cT_{k_m}\setminus\cT^+_{k_m\to k_n}\right)\cup \cT^+_{k_m\to k_n} \right)\right)\\
				& = (\cT_{k_n}\setminus\cT_{k_m})\cup \left(\cT_{k_n}\cap (\cT_{k_m}\setminus\cT^+_{k_m\to k_n})\right)  \cup \left(\cT_{k_n}\cap\cT^+_{k_m\to k_n}\right)\\
				& =
				(\cT_{k_n}\setminus\cT_{k_m})\cup \left(\cT_{k_n}\cap (\cT_{k_m}\setminus\cT^+_{k_m\to k_n})\right)  \cup \cT^+_{k_m\to k_n},
				%\\
				%& \subseteq
				%(\cT_{k_n}\setminus\cT_{k_m})\cup %\left(\cT_{k_n}\cap %(\cT_{k_m}\setminus\cT^+_{k_m})\right)  \cup %\cT^+_{k_m\to k_n}
				%                = (\cT_{k_n}\setminus\cT_{k_m})\cup %\left(\cT_{k_n}\cap\cT_{k_m}^0\right)  \cup \cT^+_{k_m\to %k_n},
			\end{aligned}
			\]
\begin{equation}\label{thm:estimator->zero_pf02}
    \eta_{k_n,j}^2 = \eta_{k_n,j}^2(\cT_{k_n}\setminus\cT_{k_m}) + \eta_{k_n,j}^2(\cT_{k_n}\cap (\cT_{k_m}\setminus\cT^+_{k_m\to k_n})) + \eta_{k_n,j}^2(\cT^+_{k_m\to k_n})\quad 0\leq j\leq l.
\end{equation}
Next we estimate the three summands in the right hand of \eqref{thm:estimator->zero_pf02} respectively.
			
			\noindent\textit{Step 1.} First, reduction properties \eqref{est-0_reduction} and \eqref{est-j_reduction} in Lemma \ref{lem:est_reduction} imply that
			\begin{equation}\label{thm:estimator->zero_pf03}
				\begin{aligned}
					\eta_{k_n,0}^2(\cT_{k_n}\setminus\cT_{k_m}) & \leq q_{0} \eta_{k_m,0}^2(\cT_{k_m}\setminus\cT_{k_n}) + c \Bigg ( \|\bold{\nabla}(\phi_{k_n}^\ast - \phi_{k_m}^\ast)\|_{L^2(D)}^2 + \|f'(\phi_{k_n}^\ast) - f'(\phi_{k_m}^\ast)\|^2_{L^2(D)}  \\
					& \qquad   + \sum_{j=1}^{l}\left\|\frac{\Psi(\lambda_{k_n,i_1}^{\eps,\ast},\cdots,\lambda_{k_n,i_l}^{\eps,\ast})}{\partial \lambda_{i_j}^{\eps}}(w_{k_n,i_j}^\ast)^2-\frac{\Psi(\lambda_{k_m,i_1}^{\eps,\ast},\cdots,\lambda_{k_m,i_l}^{\eps,\ast})}{\partial \lambda_{i_j}^{\eps}}(w_{k_m,i_j}^\ast)^2\right\|_{L^2(D)}^2
					\Bigg)\\
					& := q_{0} \eta_{k_m,0}^2(\cT_{k_m}\setminus\cT_{k_n}) + \epsilon_0(k_n,k_m),
				\end{aligned}
			\end{equation}
			\begin{equation}\label{thm:estimator->zero_pf04}
				\begin{aligned}
					\eta_{k_n,j}^2 (\cT_{k_n}\setminus\cT_{k_m}) & \leq q_{j}\eta_{k_m,j}^2 (\cT_{k_m}\setminus\cT_{k_n}) + c \Big(\|\bold{\nabla}(w_{k_n,i_j}^\ast-w_{k_m,i_j}^\ast)\|^2_{L^2(D)}\\
					& \qquad  + \|(\phi_{k_n}^\ast-\lambda_{k_n,i_j}^{\eps,\ast})w_{k_n,i_j}^\ast - (\phi_{k_m}^\ast-\lambda_{k_m,i_j}^{\eps,\ast})w_{k_m,i_j}^\ast\|_{L^2(D)}^2\Big)\\
					&:= q_{j} \eta_{k_m,j}^2(\cT_{k_m}\setminus\cT_{k_n}) + \epsilon_j(k_n,k_m),\quad 1\leq j\leq l.
				\end{aligned}
			\end{equation}
			From Theorem \ref{thm:alg_conv_limit}, we know $\{\phi_{k_n}^\ast\}_{n\geq0}$ and $\{w_{k_n,i_j}^\ast\}_{n\geq0}$ ($1\leq j\leq j$) are both Cauchy sequences in $H^1(D)$ and $H_0^1(D)$ respectively. In the proof of Lemma \ref{lem:max-error-ind->zero}, it has been proved that $\{f'(\phi_{k_n}^\ast)\}_{n\geq0}$ is also a Cauchy sequence in $L^2(D)$. For each $j\in \{1,2,\cdots,l\}$, we can derive
			\[
			\begin{aligned}
				&\quad \left\|\frac{\Psi(\lambda_{k_n,i_1}^{\eps,\ast},\cdots,\lambda_{k_n,i_l}^{\eps,\ast})}{\partial \lambda_{i_j}^{\eps}}(w_{k_n,i_j}^\ast)^2-\frac{\Psi(\lambda_{k_m,i_1}^{\eps,\ast},\cdots,\lambda_{k_m,i_l}^{\eps,\ast})}{\partial \lambda_{i_j}^{\eps}}(w_{k_m,i_j}^\ast)^2\right\|_{L^2(D)}^2\\
				& \leq 2\|w_{k_n,i_j}^\ast\|_{L^4(D)}^4\left|\frac{\Psi(\lambda_{k_n,i_1}^{\eps,\ast},\cdots,\lambda_{k_n,i_l}^{\eps,\ast})}{\partial \lambda_{i_j}^{\eps}} - \frac{\Psi(\lambda_{k_m,i_1}^{\eps,\ast},\cdots,\lambda_{k_m,i_l}^{\eps,\ast})}{\partial \lambda_{i_j}^{\eps}}\right|^2 \\
				&\quad + 2\left|\frac{\Psi(\lambda_{k_m,i_1}^{\eps,\ast},\cdots,\lambda_{k_m,i_l}^{\eps,\ast})}{\partial \lambda_{i_j}^{\eps}}\right|^2\|(w_{k_n,i_j}^\ast)^2-(w_{k_m,i_j}^\ast)^2\|_{L^2(D)}^2.
			\end{aligned}
			\]
			From Theorem \ref{thm:alg_conv_limit} again, Sobolev embedding theorem and $C^1$ regularity of $\Psi$, it follows  that for any given $\epsilon>0$,
\[		
    2\|w_{k_n,i_j}^\ast\|_{L^4(D)}^4\left|\frac{\Psi(\lambda_{k_n,i_1}^{\eps,\ast},\cdots,\lambda_{k_n,i_l}^{\eps,\ast})}{\partial \lambda_{i_j}^{\eps}} - \frac{\Psi(\lambda_{k_m,i_1}^{\eps,\ast},\cdots,\lambda_{k_m,i_l}^{\eps,\ast})}{\partial \lambda_{i_j}^{\eps}}\right|^2 < \frac{\epsilon}{2}
\]
when $m$ and $n$ are sufficiently large. Noting
\[		
    \|(w_{k_n,i_j}^\ast)^2-(w_{k_m,i_j}^\ast)^2\|_{L^2(D)}^2\leq \|w_{k_n,i_j}^\ast + w_{k_m,i_j}^\ast\|_{L^4(D)}^2 \|w_{k_n,i_j}^\ast - w_{k_m,i_j}^\ast\|_{L^4(D)}^2,
\]
we deduce in a similar manner
\[			
    2\left|\frac{\Psi(\lambda_{k_m,i_1}^{\eps,\ast},\cdots,\lambda_{k_m,i_l}^{\eps,\ast})}{\partial \lambda_{i_j}^{\eps}}\right|^2\|(w_{k_n,i_j}^\ast)^2-(w_{k_m,i_j}^\ast)^2\|_{L^2(D)}^2 < \frac{\epsilon}{2}.
\]
Therefore, $\left\{\frac{\Psi(\lambda_{k_n,i_1}^{\eps,\ast},\cdots,\lambda_{k_n,i_l}^{\eps,\ast})}{\partial \lambda_{i_j}^{\eps}}(w_{k_n,i_j}^\ast)^2\right\}_{n\geq0}$ ($1\leq j\leq l$) is a Cauchy sequence in $L^2(D)$. Finally, we bound the third term in the right hand side of \eqref{thm:estimator->zero_pf04} as
\[
	\begin{aligned}
    &\quad \|(\alpha\phi_{k_n}^\ast-\lambda_{k_n,i_j}^{\eps,\ast})w_{k_n,i_j}^\ast - (\alpha\phi_{k_m}^\ast-\lambda_{k_m,i_j}^{\eps,\ast})w_{k_m,i_j}^\ast\|_{L^2(D)}^2\quad 1\leq j\leq l\\
	&\leq
    2\|(\alpha\phi_{k_n}^\ast-\alpha\phi_{k_m}^\ast+\lambda_{k_m,i_j}^{\eps,\ast}-\lambda_{k_n,i_j}^{\eps,\ast})w_{k_n,i_j}^\ast\|_{L^2(D)}^2 + 2\|(\alpha\phi_{k_m}^\ast-\lambda_{k_m,i_j}^{\eps,\ast})(w_{k_n,i_j}^\ast-w_{k_m,i_j}^\ast)\|_{L^2(D)}^2\\
	&\leq
    4\|\alpha(\phi_{k_n}^\ast-\phi_{k_m}^\ast)w_{k_n,i_j}^\ast\|_{L^2(D)}^2 + 4\left|\lambda_{k_m,i_j}^{\eps,\ast}-\lambda_{k_n,i_j}^{\eps,\ast}\right|^2\|w_{k_n,i_j}^\ast\|_{L^2(D)}^2 \\
    &\quad + 2(\alpha\|\phi_{k_m}^\ast\|_{L^\infty(D)}+\lambda_{k_m,i_j}^{\eps,\ast})^2\|w_{k_n,i_j}^\ast-w_{k_m,i_j}^\ast\|_{L^2(D)}^2.
    \end{aligned}
\]
In view of H\"{o}lder inequality, the Sobolev embedding theorem, Theorem \ref{thm:alg_conv_limit},  $\|w_{k_n,i_j}^{\ast}\|_{L^2(D)}=1$ and $\|\phi_{k_m}^\ast\|_{L^\infty(D)}\leq 1$, the three terms in the right hand side may be estimated further respectively by
\[
    4\|\alpha(\phi_{k_n}^\ast-\phi_{k_m}^\ast)w_{k_n,i_j}^\ast\|_{L^2(D)}^2
    \leq 4 \alpha^2 \|\phi_{k_n}^\ast-\phi_{k_m}^\ast\|^2_{L^4(D)}\|w_{k_n,i_j}^\ast\|_{L^4(D)}^2
    \leq c \|\phi_{k_n}^\ast-\phi_{k_m}^\ast\|^2_{H^1(D)},
\]
\[
    4\left|\lambda_{k_m,i_j}^{\eps,\ast}-\lambda_{k_n,i_j}^{\eps,\ast}\right|^2\|w_{k_n,i_j}^\ast\|_{L^2(D)}^2
    =4\left|\lambda_{k_m,i_j}^{\eps,\ast}-\lambda_{k_n,i_j}^{\eps,\ast}\right|^2,
\]
\[
    2(\alpha\|\phi_{k_m}^\ast\|_{L^\infty(D)}+\lambda_{k_m,i_j}^{\eps,\ast})^2\|w_{k_n,i_j}^\ast-w_{k_m,i_j}^\ast\|_{L^2(D)}^2\leq
    c\|w_{k_n,i_j}^\ast-w_{k_m,i_j}^\ast\|_{L^2(D)}^2.
\]
So Theorem \ref{thm:alg_conv_limit} implies that $\|(\alpha\phi_{k_n}^\ast-\lambda_{k_n,i_j}^{\eps,\ast})w_{k_n,i_j}^\ast -(\alpha\phi_{k_m}^\ast-\lambda_{k_m,i_j}^{\eps,\ast})w_{k_m,i_j}^\ast\|_{L^2(D)}^2$ tends to zero as $m$ and $n\to \infty$. Thus we conclude from \eqref{thm:estimator->zero_pf03} and \eqref{thm:estimator->zero_pf04} that each of the $l+1$ subsequences of estimators $\{\eta_{k_m,j}\}_{m\geq0}$ ($0\leq j\leq l$) on refined elements is contractive up to a quantity $\epsilon_j(k_n,k_m)$ vanishing in the limit as Algorithm \ref{alg_afem_eigenvalue} goes on.
			
			\noindent\textit{Step 2.} From the easy observation $\cT_{k_m}^+\subseteq\cT_{k_m\to k_n}^+$, it follows that  $\cT_{k_m}\setminus \cT_{k_m\to k_n}^+ \subseteq \cT_{k_m}\setminus \cT_{k_m}^+ = \cT_{k_m}^0$. We note that $\cT_{k_n}\cap\cT_{k_m}^0 $ with $n>m$ is the set of elements  that will eventually be refined after the $k_m$-th iteration and belong to $\cT_{k_n}$ at the same time. This indicates that any $T\in \cT_{k_n}\cap\cT_{k_m}^0$ is not refined until after the $k_n$-th iteration. Thus by the definition of $\cT_{k_n}^0$, we obtain $\cT_{k_n}\cap\cT_{k_m}^0=\cT_{k_n}^0\cap\cT_{k_m}^0$. When $n$ is sufficiently large, due to $\eqref{mesh-size->zero}$ the mesh size of each element in $\cT_{k_n}^0$ is smaller than that of each one in $\cT_{k_m}^0$ for all fixed $m\in \mathbb{N}_0$. Therefore, letting $D\left(\cT^0_{k_m}\cap\cT^0_{k_n}\right):=\bigcup_{T\in\cT^0_{k_m}\cap\cT^0_{k_n}}T$ we may infer that $\cT_{k_m}^0\cap\cT_{k_n}^0=\emptyset$ for sufficiently large $n$ and \begin{equation}\label{thm:estimator->zero_pf05}
				\lim_{n\to\infty}\left|D\left(\cT^0_{k_m}\cap\cT^0_{k_n}\right)\right|=0
			\end{equation}
			for all fixed $m\in \mathbb{N}_0$. Now the scaled trace theorem, the inverse estimate, the fact $\cT_{k_n}\cap(\cT_{k_m}\setminus \cT_{k_m\to k_n}^+)\subseteq \cT_{k_n}^0\cap\cT_{k_m}^0$ and the triangle inequality yield
			\[
			\begin{aligned}
				&\quad\eta_{k_n,0}^2(\cT_{k_n}\cap (\cT_{k_m}\setminus\cT^+_{k_m\to k_n})) \\
				& \leq c \sum_{T\in \cT_{k_n}\cap \left(\cT_{k_m}\setminus\cT^+_{k_m\to k_n}\right)} \Bigg( h_{T}^2\sum_{j=1}^l \left|\frac{\partial \Psi(\lambda_{k_n,i_1}^{\eps,\ast},\cdots,\lambda_{k_n,i_l}^{\eps,\ast})}{\lambda_{i_j}^\eps}\right|^2\|w_{k_n,i_j}^\ast\|^4_{L^4(T)} + \|\bold{\nabla}\phi_{k_n}^\ast\|_{L^2(\omega_{k_n}(T))}^2 \\
				&\hspace{4.5cm}+h_{T}^2 \|f'(\phi_{k_n}^\ast)\|^2_{L^2(T)} \Bigg)\\
				&\leq c \left(\sum_{j=1}^{l}\|w_{k_n,i_j}^\ast\|^4_{L^4\left(D(\cT^0_{k_m}\cap\cT^0_{k_n})\right)}+\|\bold{\nabla}\phi_{k_n}^\ast\|_{L^2\left(D(\cT^0_{k_m}\cap\cT^0_{k_n})\right)}^2 +\|f'(\phi_{k_n}^\ast)\|^2_{L^2\left(D(\cT^0_{k_m}\cap\cT^0_{k_n})\right)}\right)\\
				&\leq c \bigg (\sum_{j=1}^{l}\|w_{k_n,i_j}^\ast-w_{\infty,i_j}^\ast\|^4_{L^4(D)}+\|\bold{\nabla}(\phi_{k_n}^\ast-\phi_\infty^\ast)\|_{L^2(D)}^2 +
				\|f'(\phi_{k_n}^\ast)-f'(\phi_\infty^\ast)\|^2_{L^2(D)}\\
				& \qquad \qquad
				+
				\sum_{j=1}^{l}\|w_{\infty,i_j}^\ast\|^4_{L^4\left(D(\cT^0_{k_m}\cap\cT^0_{k_n})\right)}
				+ \|\bold{\nabla}\phi_\infty^\ast\|_{L^2\left(D(\cT^0_{k_m}\cap\cT^0_{k_n})\right)}^2+\|f'(\phi_\infty^\ast)\|^2_{L^2\left(D(\cT^0_{k_m}\cap\cT^0_{k_n})\right)}\bigg).
			\end{aligned}
			\]
			Here the $C^1$ regularity of $\Psi$ and the boundedness of each sequence $\{\lambda_{k_n,i_j}^{\eps,\ast}\}_{n\geq0}$ ($1\leq j\leq l$) are also invoked. For each fixed $m\in \mathbb{N}_0$, passing to the limit $n\to\infty$ in the right hand side of the above estimate, the first three terms tend to zero by the argument as in the proof of Lemma \ref{lem:max-error-ind->zero} while \eqref{thm:estimator->zero_pf05} and the absolute continuity of norms $\|\cdot\|_{L^4(D)}$ and $\|\cdot\|_{L^2(D)}$ with respect to the Lebesgue measure  ensure that the remaining three ones vanishes in the limit. Hence we obtain for each fixed $m\in \mathbb{N}_0$,      \begin{equation}\label{thm:estimator->zero_pf06}
				\eta_{k_n,0}^2(\cT_{k_n}\cap (\cT_{k_m}\setminus\cT^+_{k_m\to k_n}))\to 0 \quad \text{as}~n\to\infty.
			\end{equation}
			Likewise, by the fact $\|\phi_{k_n}^\ast\|_{L^\infty(D)}\leq 1$ and the boundedness of each sequence $\{\lambda_{k_n,i_j}^{\eps,\ast}\}_{n\geq0}$ again there holds for $1\leq j\leq l$ and each fixed $m\in\mathbb{N}_0$,
			\begin{align}
				\eta_{k_n,j}^2(\cT_{k_n}\cap (\cT_{k_m}\setminus\cT^+_{k_m\to k_n}))  & \leq c\left( \|\bold{\nabla}w_{k_n,i_j}^{\ast}\|^2_{L^2\left(D(\cT^0_{k_m}\cap\cT^0_{k_n})\right)} +  \|(\alpha\phi_{k_n}^\ast-\lambda_{k_n,i_j}^{\eps,\ast})w_{k_n,i_j}^\ast\|_{L^2\left(D(\cT^0_{k_m}\cap\cT^0_{k_n})\right)}^2\right)\nonumber\\
				& \leq c \Big( \|\bold{\nabla}(w_{k_n,i_j}^{\ast}-w_{\infty,i_j}^{\ast})\|^2_{L^2(D)} + \|\bold{\nabla}w_{\infty,i_j}^{\ast}\|^2_{L^2\left(D(\cT^0_{k_m}\cap\cT^0_{k_n})\right)}\nonumber\\
				& \qquad + \|w_{k_n,i_j}^\ast-w_{\infty,i_j}^\ast\|_{L^2(D)} + \|w_{\infty,i_j}^\ast\|_{L^2\left(D(\cT^0_{k_m}\cap\cT^0_{k_n})\right)}
				\Big)\to 0  \quad \text{as}~n\to\infty. \label{thm:estimator->zero_pf07}
			\end{align}
			
\noindent{\textit{Step 3.}} The module \textsf{MARK} in Algorithm \ref{alg_afem_eigenvalue} employs a separate marking for each $\eta_{k_n,j}$ $(0\leq j\leq l)$. Then by    \eqref{thm:estimator->zero_pf01}, \eqref{marking} and \eqref{max-error-ind-zero} in Lemma \ref{lem:max-error-ind->zero} we obtain for each fixed $m\in \mathbb{N}_0$,
\begin{equation}\label{thm:estimator->zero_pf08}
    \eta_{k_n,j}^2(\cT_{k_m\to k_n}^+) \leq \#\cT_{k_m\to k_n}^+ \max_{ T \in \cT_{k_m\to k_n}^+ }\eta_{k_n,j}^2(T) \leq c \#\cT_{k_m}^+ \max_{ T \in \mathcal{M}_{k_n} }\eta_{k_n,j}^2(T)\to 0\quad \text{as}~ n\to \infty,~0\leq j\leq l.
\end{equation}
%If $\widetilde{\eta}_{k_n} := \eta_{k_n,j}$ ($1\leq j\leq l$), a similar argument shows that for each fixed $m\in \mathbb{N}_0$
%			\begin{equation}\label{thm:estimator->zero_pf09}
%				\widetilde{\eta}_{k_n}^2(\cT_{k_m\to k_n}^+)=\eta_{k_n,j}^2(\cT_{k_m\to k_n}^+) \to 0\quad n\to \infty. %\leq \#\cT_{k_m\to k_n}^+ \max_{ T \in \cT_{k_m\to k_n}^+ %}\eta_{k_n,j}^2(T) \leq c \#\cT_{k_m}^+ \max_{ T \in %\cT_{k_m}^+ }\eta_{k_n,j}^2(T)\to 0\quad n\to \infty.
%			\end{equation}

			\noindent\textit{Step 4.} From  \eqref{thm:estimator->zero_pf06}-\eqref{thm:estimator->zero_pf08}, it follows that for each $m\in \mathbb{N}_0$, there exists some $N(m)>m$ such that for each $n>N(m)$,
\begin{equation}\label{thm:estimator->zero_pf10}
    \eta_{k_{n},j}^2(\cT_{k_n}\cap\cT_{k_m}) \leq
    (q' - q_j) \eta_{k_m,j}^2
\end{equation}
with each $q_j\in (0,1)$ in \eqref{thm:estimator->zero_pf03}-\eqref{thm:estimator->zero_pf04}  and some $q' \in (\max_{0\leq j\leq l}q_j,1)$. Using  \eqref{thm:estimator->zero_pf03}-\eqref{thm:estimator->zero_pf04} again and \eqref{thm:estimator->zero_pf10},
we may extract a further subsequence $\{k_{n_p}\}_{p\geq0}$ such that
\[
    \eta_{k_{n_{p+1}},j}^2 \leq q' \eta_{k_{n_{p}},j}^2 + \epsilon_{k_{n_{p}},j}\quad 0\leq j\leq l,
\]
where $\epsilon_{k_{n_p},j} = \epsilon_j(k_{n_{p+1}},k_{n_{p}})\to 0$ as $p\to\infty$ due to the conclusion at the end of Step 1.
			%\[
			%    \begin{aligned}
				%    \eta_{k_{n_{p+1}},0}^2 & \leq q'_0 \eta_{k_{n_{p}},0}^2 + %c_0 \Bigg ( \|\bold{\nabla}(\phi_{k_{n_{p+1}}}^\ast - %\phi_{k_{n_p}}^\ast)\|_{L^2(D)}^2 + %\|f'(\phi_{k_{n_{p+1}}}^\ast) - %f'(\phi_{k_{k_n}}^\ast)\|^2_{L^2(D)} +  \\
				%            & \quad   %\sum_{j=1}^{l}\left\|\frac{\Psi(\lambda_{k_{n_{p+1}},i_1}^{\eps,\ast},\cdots,\lambda_{k_{n_{p+1}},i_l}^{\eps,\ast})}{\partial %\lambda_{k_{n_{p+1}},i_j}^{\eps}}(w_{k_{n_{p+1}},i_j}^\ast)^2-\frac{\Psi(\lambda_{k_{n_p},i_1}^{\eps,\ast},\cdots,\lambda_{k_{n_p},i_l}^{\eps,\ast})}{\partial %\lambda_{k_{n_p},i_j}^{\eps}}(w_{k_{n_p},i_j}^\ast)^2\right\|_{L^2(D)}^2
				%            \Bigg)\\
				%            &: =  q'_0 \eta_{k_{n_{p}},0}^2 + %\epsilon_{k_{n_p},0},
				%    \end{aligned}
			%\]
			%\[
			%    \begin{aligned}
				%        \eta_{k_{n_{p+1}},j}^2 & \leq %q_{j}'\eta_{k_{n_p},j}^2 +  c_j %\Big(\|\bold{\nabla}(w_{k_{n_{p+1}},i_j}^\ast-w_{k_{n_p},i_j}^\ast)\|^2_{L^2(D)}
				%            + \\
				%            &\qquad\qquad %\|(\phi_{k_{n_{p+1}}}^\ast-\lambda_{k_{n_{p+1}},i_j}^{\eps,\ast})w_{k_{n_{p+1}},i_j}^\ast %-             %(\phi_{k_{n_{p}}}^\ast-\lambda_{k_{n_p},i_j}^{\eps,\ast})w_{k_{n_p},i_j}^\ast\|_{L^2(D)}^2\Big)\quad %1\leq j\leq l \\
				%            &: = q_{j}'\eta_{k_{n_p},j}^2 +             %\epsilon_{k_{n_p},j}.
				%    \end{aligned}
			%\]
			%Summing up the above estimates, we arrive at
			%$
			%    \sum_{j=0}^{l}\eta_{k_{n_{p+1}},j}^2 \leq \max_{0\leq %j\leq l} q_{j}' \sum_{j=0}^l \eta_{k_{n_p},j}^2 +  %\sum_{j=0}^l \epsilon_{k_{n_{p}},j}
			%$
			%with $\sum_{j=0}^l \epsilon_{k_{n_{p}},j} \to 0$ as %$p\to\infty$ due to the conclusion at the end of Step 1.
			So it follows from Lemma 3.2 in \cite{GanterPraetorius:2022} that %$\lim_{p\to\infty}\widetilde{\eta}_{k_{n_{p}}}^2 = 0$, which immediately implies that
\begin{equation}\label{thm:estimator->zero_pf11}
    \lim_{p\to\infty}\eta_{k_{n_p},j}^2 = 0 \quad 0\leq j\leq l.
\end{equation}
By \eqref{est-0_stab}-\eqref{est-j_reduction}, we argue as in Step 1 to obtain  $\eta_{k_n,j}^2\leq 2\eta_{k_{n_p},j}^2 + \xi_j(k_n,k_{n_p})$ ($0\leq j\leq l$), where $\xi_j(k_n,k_{n_p})$ tends to zero as $k_n>k_{n_p}\to \infty$ and $p\to\infty$. This and \eqref{thm:estimator->zero_pf11} conclude the proof.
			%Then for any $\epsilon>0$, there exists a $n_{p}$ such that %$2\eta_{k_{n_{p}}^2<\frac{\epsilon}{2}$ and %$\xi(k_n,k_{n_p})<\frac{\epsilon}{2}$ for any $n>n_{p}$.
				%$\max_{0\leq j\leq j}c_j $
				%    The definition of $\cT_{k_m}^+$ implies that %$\cT_{k_m}\cap \cT_{k_n}$ for $m<n$ By %\eqref{est-0_reduction} and \eqref{est-j_reduction} in Lemma %\ref{lem:est_reduction}, we split for $m<n$
				%    \[
				%        \eta_{k_n,0} = %\eta_{k_n,0}(\cT_{k_n}\setminus\cT_{k_m}) +
				%    \]
			\end{proof}
			
			\begin{theorem}\label{thm:alg_conv}
				Let $\left\{\phi_{k}^\ast,\left(\lambda_{k,i_j}^{\eps,\ast},w_{k,i_j}^\ast\right)_{j=1}^l\right\}_{k\geq0}$ be the sequence of discrete minimizers and the $l$ sequences of corresponding eigenpairs solving \eqref{eigen-opt_phasefield_disc}-\eqref{eigen_vp_phasefield_disc} over $\{\cT_k\}_{k\geq0}$ generated by Algorithm \ref{alg_afem_eigenvalue}. Assuming all eigenvalues $\left\{\lambda_{k,i_j}^{\eps,\ast}\right\}_{k\geq0}$ $(1\leq j\leq l)$ are simple, then there exists a subsequence $\{\phi_{k_n}^\ast\}_{n\geq0}$ and $l$ subsequences $\left\{\left(\lambda_{k_n,i_j}^{\eps,\ast},w_{k_n,i_j}^\ast\right)_{j=1}^l\right\}_{n\geq0}$ converging strongly to some $\overline{\phi}$ in the $H^1(D)$-topology and $l$ corresponding eigenpairs in the $\mathbb{R}\times H_0^1(D)$-topology, which solve \eqref{optsys_phase-field}.
				%\begin{equation}\label{alg_conv}
				%        \phi_{k_n}^\ast\to\overline{\phi}\quad %\text{in}~H^1(D),\quad \lambda_{k_n,i_j}^{\eps,\ast}\to %\overline{\lambda}_{i_j}^{\eps},\quad  w_{k_n,i_j}^{\ast} \to %\overline{w}_{i_j}\quad  \text{in}~H^1_0(D).
				%    \end{equation}
		\end{theorem}
		
		\begin{proof}
			As mentioned as the beginning of this section, it is sufficient to prove that the limit in Theorem \ref{thm:alg_conv_limit}
			\[
			\left\{\phi_\infty^\ast, \left(\lambda_{\infty,i_1}^{\eps,\ast},w_{\infty,i_1}\right), \cdots, \left(\lambda_{\infty,i_l}^{\eps,\ast},w_{\infty,i_l}\right)\right\}\in \U_\infty\times (\mathbb{R}_{>0}\times V_\infty^0)^l
			\]
			solve \eqref{optsys_phase-field}. Invoking the $H^1(D)$ strong convergent subsequence $\{\phi_{k_n}^\ast\}_{n\geq 0}$ from Theorem \ref{thm:alg_conv_limit}, we find 
			\begin{equation}\label{thm:alg_conv_pf01}
				\lim_{n\to\infty}\int_{D} \bold{\nabla}\phi_{k_n}^\ast \cdot \bold{\nabla}(\phi-\phi_{k_n}^\ast) \dx = \int_{D} \bold{\nabla}\phi_{\infty}^\ast \cdot \bold{\nabla}(\phi-\phi_{\infty}^\ast) \dx \quad \forall \phi\in \mathcal{U}\textcolor{red}{.}
			\end{equation}
			By the pointwise convergence of $\{\phi_{k_n}^\ast\}_{n\geq 0}$ to $\phi_\infty^\ast$ in Remark \ref{rem:alg_conv_limit} and Lebesgue's dominated convergence theorem,
			\begin{equation}\label{thm:alg_conv_pf02}
				\lim_{n\to\infty}\int_D f'(\phi_{k_n}^\ast)(\phi - \phi_{k_n}^\ast) \dx  \textcolor{blue}{=} \int_D f'(\phi_{\infty}^\ast) (\phi - \phi_\infty^\ast) \dx \quad \forall \phi\in\mathcal{U}.
			\end{equation}
			Due to the fact $\left|(w_{\infty,i_j}^\ast)^2\left(\phi_{k_n}^\ast-\phi_\infty^\ast\right)\right|\leq (w_{\infty,i_j}^\ast)^2\in L^1(D)$ almost everywhere in $D$, H\"{o}lder inequality, the $H_0^1(D)$ strong convergence of $\{w_{k_n,i_j}^\ast\}_{n\geq0}$ from Theorem \ref{thm:alg_conv_limit}, the pointwise convergence of $\left\{\phi_{k_n}^\ast\right\}_{n\geq 0}$ and Lebesgue's dominated convergence theorem further imply that passing to the limit $n\to\infty$ for each $j\in\{1,2,\cdots,l\}$,
			\[
			\begin{aligned}
				\left| \int_{D} \left((w_{k_n,i_j}^\ast)^2 - (w_{\infty,i_j}^\ast)^2 \right) \phi  \dx \right| & \leq \|w_{k_n,i_j}^\ast + w_{\infty,i_j}^\ast\|_{L^2(D)}
				\|w_{k_n,i_j}^\ast - w_{\infty,i_j}^\ast\|_{L^2(D)}\|\phi\|_{L^\infty(D)}\\
				&\leq 2\|w_{k_n,i_j}^\ast - w_{\infty,i_j}^\ast\|_{L^2(D)}\to 0 \quad \forall \phi\in\mathcal{U},
			\end{aligned}
			\]
			\[
			\begin{aligned}
				 \left|\int_{D}(w_{k_n,i_j}^\ast)^2\phi_{k_n}^\ast \dx - \int_{D} (w_{\infty,i_j}^\ast)^2\phi_{\infty}^\ast \dx \right|
				&\leq \left|\int_{D}\left((w_{k_n,i_j}^\ast)^2 - (w_{\infty,i_j}^\ast)^2\right) \phi_{k_n}^\ast \dx\right|+\left|\int_{D}(w_{\infty,i_j}^\ast)^2\left(\phi_{k_n}^\ast-\phi_\infty^\ast\right) \dx\right|\\
				&\leq 2\|w_{k_n,i_j}^\ast - w_{k_n,\infty}^\ast\|_{L^2(D)} + \int_{D}\left|(w_{\infty,i_j}^\ast)^2\left(\phi_{k_n}^\ast-\phi_\infty^\ast\right)\right| \dx\to 0 \quad \forall \phi\in\mathcal{U}.
			\end{aligned}
			\]
			A collection of the above two estimates, along with the convergence of $\left\{\lambda_{k_n,i_j}^{\eps,\ast}\right\}_{n\geq 0}$ ($1\leq j\leq l$) from Theorem \ref{thm:alg_conv_limit} and the $C^1$ regularity of $\Psi$, yields that as $n\to\infty$,
			\begin{equation}\label{thm:alg_conv_pf03}
				\begin{aligned}
					\alpha\sum_{j=1}^l\frac{\partial \Psi(\lambda_{k_n,i_1}^{\eps,\ast},\cdots,\lambda_{k_n,i_l}^{\eps,\ast})}{\lambda_{i_j}^\eps}
					&\int_D (w_{k_n,i_j}^\ast)^2   \left( \phi - \phi_{k_n}^\ast\right) \dx \\
					& \to \alpha\sum_{j=1}^l\frac{\partial \Psi(\lambda_{\infty,i_1}^{\eps,\ast},\cdots,\lambda_{\infty,i_l}^{\eps,\ast})}{\lambda_{i_j}^\eps} \int_D(w_{\infty,i_j}^\ast)^2 \left( \phi - \phi_{\infty}^\ast\right) \dx \quad \forall \phi\in \mathcal{U}.
				\end{aligned}
			\end{equation}
			On the other hand, recalling \eqref{residual_0}, \eqref{eqn:est_derive01}-\eqref{eqn:est_derive02} and noting $\eta_{k_n,0}\to 0$ as $n\to\infty$ given by Theorem \ref{thm:estimator->zero}, we obtain
			\begin{equation}\label{thm:alg_conv_pf04}
				\liminf_{n\to\infty}\langle \mathcal{R}_0(\phi_{k_n}^\ast), \phi - \phi_{k_n}^\ast \rangle  \geq 0\quad \forall \phi\in \mathcal{U}.
			\end{equation}
			It follows from \eqref{thm:alg_conv_pf01}-\eqref{thm:alg_conv_pf04}  that $\overline{\phi}:=\phi_{\infty}^\ast$ and $\left\{\left(\lambda_{\infty,i_1}^{\eps,\ast},w_{\infty,i_1}\right), \cdots, \left(\lambda_{\infty,i_l}^{\eps,\ast},w_{\infty,i_l}\right)\right\}$ solve the first variational inequality of \eqref{optsys_phase-field}. By Theorem \ref{thm:alg_conv_limit} again, we have for $1\leq j\leq l$,
			\[
			\int_D \bold{\nabla} w_{k_n,i_j}^\ast \cdot \bold{\nabla} v \dx \to \int_D \bold{\nabla} w_{\infty,i_j}^\ast \cdot \bold{\nabla} v \dx,\quad \lambda_{k_n,i_j}^{\eps,\ast}\int_D w_{k_n,i_j}^\ast v \dx\to \lambda_{\infty,i_j}^{\eps,\ast} \int_{D}w_{\infty,i_j}^\ast v \dx \quad \forall v\in H_0^1(D).
			\]
			As $\phi_{\infty}^\ast\in [0,1]$ and  $|(\phi_{k_n}^\ast-\phi_\infty^\ast)v|^2 \leq |v|^2\in L^1(D)$ almost everywhere in $D$ for all $v\in H_0^1(D)$, the pointwise convergence of $\{\phi_{k_n}^\ast\}_{n\geq0}$, Lebesgue's dominated convergence theorem and Theorem \ref{thm:alg_conv_limit} imply for each $j\in\{1,2,\cdots,l\}$,
\[
    \begin{aligned}
        \left|\int_{D} \phi_{k_n}^\ast w_{k_n,i_j}^\ast v  \dx - \int_{D} \phi_{\infty}^\ast w_{\infty,i_j}^\ast v  \dx \right| & \leq \left| \int_{D} \left(\phi_{k_n}^\ast-\phi_\infty^\ast\right) w_{k_n,i_j}^\ast v \dx \right| + \left| \int_{D} \phi_{\infty}^\ast \left( w_{k_n,i_j}^\ast - w_{\infty,i_j}^\ast \right) v\dx \right| \\
				& \leq \|(\phi_{k_n}^\ast-\phi_\infty^\ast)v\|_{L^2(D)} +
				\|w_{k_n,i_j}^\ast - w_{\infty,i_j}^\ast\|_{L^2(D)}\|v\|_{L^2(D)} \to 0 \quad \text{as}~ n\to\infty.
\end{aligned}
\]
Therefore, in view of the definition \eqref{residual_j} there holds 
\[
    \lim_{n\to\infty}\langle \mathcal{R}_j(\phi_{k_n}^{\ast},\lambda_{k_n,i_j}^{\eps,\ast},w_{k_n,i_j}^{\ast}), v \rangle \textcolor{blue}{=} \int_{D} \bold{\nabla} w_{\infty,i_j}^\ast \cdot \bold{\nabla} v \dx + \alpha \int_{D} \phi_\infty^\ast w_{\infty,i_j}^\ast v  \dx - \lambda_{\infty,i_j}^{\eps,\ast}\int_{D}  w_{\infty,i_j}^\ast v \dx\quad \forall v\in H_0^1(D).
\]
			As $\eta_{k_n,j}\to 0$ as $n\to\infty$ by Theorem \ref{thm:estimator->zero}, we invoke \eqref{eqn:est_derive03} to find
			$\langle \mathcal{R}_j(\phi_{k_n}^{\ast},\lambda_{k_n,i_j}^{\eps,\ast},w_{k_n,i_j}^{\ast}), v \rangle\to 0$ as $n\to\infty$ for any $v\in H_0^1(D)$ and $j\in \{1,2,\cdots,l\}$. As a consequence, $\overline{\phi}$ and each of $\left\{\left(\lambda_{\infty,i_j}^{\eps,\ast},w_{\infty,i_j}^\ast\right)\right\}_{j=1}^l$ solve the second variational equation of \eqref{optsys_phase-field}.
	\end{proof}
		
\begin{remark}\label{rem:alg_conv}
We have proved that the $l$ sequences of discrete eigenpairs, associated with adaptively-generated discrete minimizers $\{\phi^\ast_k\}_{k\geq0}$,  contain $l$ convergent subsequences, whose limits   $\left(\lambda_{\infty,i_j}^{\eps,\ast},w_{\infty,i_j}^\ast\right)$ each solve \eqref{eigen_vp_phasefield} with $\phi=\overline{\phi}$. Now that the $i_j$-th $(1\leq j\leq l)$ discrete eigenpair over $\cT_k$ is chosen as the approximation, one may naturally ask a question whether each limit $\left(\lambda_{\infty,i_j}^{\eps,\ast},w_{\infty,i_j}^\ast\right)$ is the $i_j$-th eigenpair of \eqref{eigen_vp_phasefield}. To answer this question seems highly nontrivial and may be a topic of our future research.
\end{remark}

\paragraph{Acknowledgement} The research of Yifeng Xu was in part supported by National Natural Science Foundation of China (Projects  12250013, 12261160361 and 12271367), Science and Technology Commission of Shanghai Municipality (Projects 20JC1413800 and 22ZR1445400) and General Research Fund from Shanghai Normal University (KF202318 and KF202468). The work of Shengfeng Zhu was supported in part by National Key Basic Research Program under grant 2022YFA1004402, National Natural Science Foundation of China under grant 12471377 and Science and Technology Commission of Shanghai Municipality (Projects 22ZR1421900 and 22DZ2229014).

\end{document}